\documentclass[english]{article}
\usepackage{graphicx} 
\graphicspath{{images/}}
\usepackage{geometry}
\usepackage[dvipsnames]{xcolor}
\usepackage{stmaryrd}
\usepackage{comment}
\usepackage{tikz-cd}
\usepackage{quiver}
\usepackage{dsfont}
\usepackage[colorlinks=true, linkcolor=mypurple, urlcolor= black, citecolor=Orange]{hyperref}
\usepackage[all]{xy}
\usepackage[utf8]{inputenc}
\usepackage{amsmath,amsthm,amssymb,amsfonts,mathtools}

\usepackage[maxbibnames=99,backend=biber,style=alphabetic]{biblatex}
\usetikzlibrary{matrix}
\usepackage{lmodern}
\usepackage[T1]{fontenc}

\definecolor{mypurple}{RGB}{155,0,255}
\definecolor{myorange}{RGB}{255,144,0}
\definecolor{mygreen}{RGB}{0,115,85}
\usepackage{float}

\usepackage{caption}
\usepackage[labelfont=bf]{caption}
\DeclareCaptionFormat{myformat}{\fontsize{9}{10}\selectfont#1#2#3}
\usepackage{placeins}
\usepackage{pinlabel}

\usepackage{aliascnt}
\usepackage[capitalize]{cleveref}

\newtheorem{thm}{Theorem}[section]

\newaliascnt{lem}{thm}
\newtheorem{lem}[lem]{Lemma}
\aliascntresetthe{lem}

\newaliascnt{prop}{thm}
\newtheorem{prop}[prop]{Proposition}
\aliascntresetthe{prop}

\newaliascnt{cor}{thm}
\newtheorem{cor}[cor]{Corollary}
\aliascntresetthe{cor}

\newaliascnt{rk}{thm}
\newtheorem{rk}[rk]{Remark}
\aliascntresetthe{rk}

\newaliascnt{notation}{thm}
\newtheorem{notation}[notation]{Notation}
\aliascntresetthe{notation}

\newaliascnt{ex}{thm}

\aliascntresetthe{ex}

\newaliascnt{constr}{thm}
\newtheorem{constr}[constr]{Construction}
\aliascntresetthe{constr}

\newaliascnt{propdef}{thm}

\aliascntresetthe{propdef}

\newaliascnt{definition}{thm}
\newtheorem{definition}[definition]{Definition}
\aliascntresetthe{definition}

\newaliascnt{cl}{thm}

\aliascntresetthe{cl}

\newtheorem{mainthm}{Theorem}

\crefname{thm}{Theorem}{Theorems}
\Crefname{thm}{Theorem}{Theorems}

\crefname{lem}{Lemma}{Lemmas}
\Crefname{lem}{Lemma}{Lemmas}

\crefname{prop}{Proposition}{Propositions}
\Crefname{prop}{Proposition}{Propositions}

\crefname{cor}{Corollary}{Corollaries}
\Crefname{cor}{Corollary}{Corollaries}

\crefname{rk}{Remark}{Remarks}
\Crefname{rk}{Remark}{Remarks}

\crefname{notation}{Notation}{Notations}
\Crefname{notation}{Notation}{Notations}

\crefname{ex}{Example}{Examples}
\Crefname{ex}{Example}{Examples}

\crefname{constr}{Construction}{Constructions}
\Crefname{constr}{Construction}{Constructions}

\crefname{propdef}{Proposition/Definition}{Propositions/Definitions}
\Crefname{propdef}{Proposition/Definition}{Propositions/Definitions}

\crefname{definition}{Definition}{Definitions}
\Crefname{definition}{Definition}{Definitions}

\crefname{cl}{Claim}{Claims}
\Crefname{cl}{Claim}{Claims}

\newcommand{\bR}{\mathbb R}
\newcommand{\bN}{\mathbb N}
\newcommand{\bS}{\mathbb S}
\newcommand{\bZ}{\mathbb Z}
\newcommand{\sus}{\Sigma^{\infty}}

\newcommand{\susplus}{\Sigma^{\infty}_+}
\newcommand{\widesus}{\widetilde{\Sigma}^{\infty}}
\newcommand{\lop}{\Omega^{\infty}}
\newcommand{\Ztwo}{\mathbb Z / 2 \mathbb Z}
\newcommand{\Ftwo}{\mathbb F_2}

\newcommand{\nupar}{\nu_{\mathrm{P}}^{\haut}}

\newcommand{\nuparStwo}{\nu_{S^2}^{\haut}}
\newcommand{\CobdG}{\mathrm{Cob}_{\mathrm{d}}^{\mathrm{G}}}
\newcommand{\PCobdG}{\mathrm{Cob}_{\mathrm{d}}^{\mathrm{(S)G}}}
\newcommand{\PCobdO}{\mathrm{Cob}_{\mathrm{d}}^{\mathrm{(S)O}}}

\newcommand{\CobdSG}{\mathrm{Cob}_{\mathrm{d}}^{\mathrm{SG}}}
\newcommand{\sCobtwo}{\mathrm{Cob}_2^{\mathrm{SG}}}
\newcommand{\sCobtwooo}{\mathrm{Cob}_2^{\mathrm{SO}}}
\newcommand{\sCobtwoH}{\mathrm{Cob}_2^{\mathrm{H}}}
\newcommand{\cobred}{\mathrm{Cob}_2^{\mathrm{SG,red}}}
\newcommand{\Bhaut}{\mathrm{Bhaut}}
\newcommand{\haut}{\mathrm{haut}}
\newcommand{\Emb}{\mathrm{Emb}}
\newcommand{\Map}{\mathrm{Map}}
\newcommand{\Diff}{\mathrm{Diff}}
\newcommand{\B}{\mathrm{B}}
\newcommand{\cSG}{\mathrm{(S)G}}
\newcommand{\cSO}{\mathrm{(S)O}}
\newcommand{\Pone}{\mathrm{P}_1}
\newcommand{\Ho}{\mathrm{H}} 

\newcommand{\sphG}{\mathrm{Sph^{SG}}}

\newcommand{\Th}{\mathrm{Th}}
\newcommand{\PTPH}{\mathrm{PT}_{\mathrm{P}}^{\mathrm{H}}}
\newcommand{\FPH}{\mathrm{F}_{\mathrm{P}}^{\mathrm{H}}}
\newcommand{\DPH}{\mathrm{D}_{\mathrm{P}}^{\mathrm{H}}}
\makeatletter 
\newcommand{\colim@}[2]{%
  \vtop{\m@th\ialign{##\cr
    \hfil$#1\operator@font colim$\hfil\cr
    }}%
}

\newcommand{\colim}{%
  \mathop{\mathpalette\colim@{\fill@\scriptscriptstyle}}\nmlimits@
}
\makeatother
\makeatletter
\newcommand{\hocolim@}[2]{%
  \vtop{\m@th\ialign{##\cr
    \hfil$#1\operator@font hocolim$\hfil\cr
    }}%
}
\newcommand{\hocolim}{%
  \mathop{\mathpalette\hocolim@{\fill@\scriptscriptstyle}}\nmlimits@
}
\makeatother

\usepackage{blindtext}

\usepackage{subfiles}

\title{The homotopy type of the Poincaré cobordism category for surfaces}

\begin{document}
\author{Azélie Picot}
\date{}
\maketitle

\begin{abstract}
    We define a version of the surface cobordism category $\sCobtwo(X)$ over a base space $X$ where surfaces are considered up to self homotopy equivalences instead of diffeomorphisms. We prove the induced functor $\B \sCobtwo(-) : \mathcal{S} \rightarrow \mathcal{S}$ is not $1$-excisive. We show its first derivative $\partial_1 B \sCobtwo(-)$ in the Goodwillie sense is equivalent to a Thom spectrum over $\Bhaut_*^+(S^2)$.
\end{abstract}

\section{Introduction}

\subsection{The main character} \label{1.1}

Smooth cobordism categories (potentially with tangential structures) have been extensively studied in \cite{GMTW}, \cite{grw9} and \cite{grw14}. They have been key to better understanding the cohomology of moduli spaces of manifolds. The homotopy type of variants of the smooth cobordism category, such as a cobordism category for topological manifolds (\cite{gomezlopez2022homotopytypetopologicalcobordism}) or a cobordism category for PL manifolds (\cite{lopez2024homotopytypeplcobordism}) 
have also been investigated. Cobordism categories for chain complexes have been studied in \cite{hebestreit2025stablemodulispaceshermitian} and \cite{calmès2021hermitianktheorystableinftycategorie} as well. In this paper, we define another variant: a cobordism category for Poincaré complexes. 

A Poincaré complex is a finite space which satisfies Poincaré duality with respect to some local coefficient system, while a Poincaré pair is a pair of spaces which satisfies relative Poincaré duality (see \cref{2.1}). Poincaré complexes and Poincaré pairs are a homotopy-theoretic analogue of manifolds and bordisms. In this paper, we define a topologically-enriched category $\CobdG$, the cobordism category of Poincaré complexes, whose objects are $(d-1)$-dimensional Poincaré complexes and whose morphism spaces are given by $$\mathrm{Mor}(P_0,P_1)\simeq \coprod_{[Q,P_0,P_1]} \Bhaut_{\partial}(Q,P_0\sqcup P_1),$$ where the disjoint union runs over $d$-dimensional Poincaré pairs $(Q,P_0\sqcup P_1)$ up to equivalence and $\haut_{\partial}(Q,P_0\sqcup P_1)$ is the space of self-equivalences of $Q$ that restrict to the identity on $P_0\sqcup P_1$. Composition is obtained by gluing morphisms along their common boundary. More generally, for each space $X$, we can define a category $\CobdG(X)$. Its objects and morphisms are respectively $(d-1)$-dimensional Poincaré complexes and $d$-dimensional Poincaré pairs equipped with a map to $X$, compatible with the boundary data. We also define an oriented version $\CobdSG(X)$ where we restrict to oriented Poincaré complexes and orientation-preserving self-homotopy equivalences.

\subsection{Main Results}

The main results of this paper concern the $2$-dimensional oriented Poincaré cobordism category $\sCobtwo(X)$. Taking the geometric realization of the nerve of $\sCobtwo(X)$ gives a functor $\B\sCobtwo(-) : \mathcal{S} \rightarrow \mathcal{S}$, where $\mathcal{S}$ denotes the $\infty$-category of spaces. Before explicitly describing $\B \sCobtwo(X)$, we first describe the first Goodwillie derivative of the functor $\B \sCobtwo(-)$. In other words, we compute the best approximation of $\B \sCobtwo(-)$ by a homology theory or an excisive functor.

In \cref{5}, we construct a spherical fibration $\nuparStwo$ over $\Bhaut_*(S^2)$, where $\haut_*(-)$ denotes self-homotopy equivalences preserving the basepoint. The spherical fibration $\nuparStwo$ is the family of the underlying spherical fibration of the stable normal bundle of $S^2$. We give a formal definition in \cref{dualizingcomplexfamilies} and \cref{notation-4.1.12}. Generalizing the parametrized Pontryagin-Thom construction for bundles from \cite{GMTW}, we define in \cref{constr-6.1-nattrans} a natural transformation $$\alpha(-) : \B \sCobtwo(-) \Rightarrow \lop (\Sigma \Th(\nuparStwo)\otimes (\susplus -)).$$ Our main statement determines the best approximation of the functor $\B \sCobtwo(-)$ at the point by an excisive functor. We denote this best approximation by $\Pone \B \sCobtwo(-)$ in the sense of Goodwillie calculus (see \cref{5.1} and \cref{P_1^X_F} for definitions).

\begin{mainthm} \label{thmA} 
The first approximation map $$p_1 \mathrm{B}\sCobtwo(-) : \mathrm{B}\sCobtwo(-) \Rightarrow \Pone \mathrm{B}\sCobtwo(-)$$ is equivalent to the natural transformation $$\gamma : \mathrm{B}\sCobtwo(-) \Rightarrow  P,$$ where $P$ is the pullback of the cospan in the following diagram and $\gamma$ is induced by $\alpha$: 

\[\begin{tikzcd}
	{} & \textcolor{mypurple}{{\mathrm{BCob}_2^{\mathrm{SG}}(-)}} \\
	&& \textcolor{mypurple}{P} & {\Omega^{\infty-1} (\mathrm{Th}(\nu_{S^2}^{\mathrm{haut}}) \wedge \Sigma^{\infty}_+-)} & {.} \\
	&& {\mathrm{BCob}_2^{\mathrm{SG}}} & {\mathrm{Th}(\nu_{S^2}^{\mathrm{haut}})}
	\arrow[draw={mypurple}, Rightarrow, from=1-2, to=2-3]
	\arrow["\alpha"{description}, curve={height=-12pt}, Rightarrow, from=1-2, to=2-4]
	\arrow["{(-\rightarrow \star)_*}"{description}, curve={height=6pt}, Rightarrow, from=1-2, to=3-3]
	\arrow[Rightarrow, from=2-3, to=2-4]
	\arrow[Rightarrow, from=2-3, to=3-3]
	\arrow["{(-\rightarrow \star)_*}"{description}, Rightarrow, from=2-4, to=3-4]
	\arrow[""{name=0, anchor=center, inner sep=0}, "{\alpha(\star)}"{description}, Rightarrow, from=3-3, to=3-4]
	\arrow["\lrcorner"{anchor=center, pos=0.125}, draw=none, from=2-3, to=0]
\end{tikzcd}\]
\end{mainthm}

In order to prove \cref{thmA}, we establish a pushout formula for the classifying space $\B \sCobtwo(\mathrm{X})$ which involves $\B \mathrm{Cob}_2^{\mathrm{SO}}$. The cobordism category $\sCobtwooo(X)$ is the topologically-enriched category with objects $1$-dimensional closed oriented smooth manifolds equipped with a map to $X$, while mapping spaces are equivalent to: $$\sCobtwooo((M_0,f_0),(M_1,f_1)) \simeq \coprod_{[\Sigma]} \Map_{\partial}(\Sigma,X)\sslash\Diff^+_{\partial}(\Sigma).$$ The disjoint union runs over $2$-dimensional oriented cobordisms $(\Sigma,M_0,M_1)$ and $\Map_{\partial}(\Sigma,X)$ denotes the space of maps restricting to $f_0,f_1$ on the boundary of $\Sigma$. Taking the classifying space gives a functor $\B\mathrm{Cob}_{\mathrm{2}}^{\mathrm{SO}}(-) : \mathcal{S} \rightarrow \mathcal{S}$. The celebrated result of \cite{GMTW} describes the homotopy type of $\B \sCobtwooo(X)$ as follows: 
 $$B \mathrm{Cob}_2^{\mathrm{SO}} (\mathrm{X}) \simeq \Omega^{\infty} (\Sigma \mathrm{MTSO}(2)\otimes \susplus X),$$ where $\mathrm{MTSO}(2)$ denotes the Thom spectrum of the stable inverse of the universal $2$-dimensional vector bundle over $\mathrm{BSO}(2)$.

The forgetful maps $\Diff_{\partial}(\Sigma) \rightarrow \haut_{\partial}(\Sigma)$ induce a natural transformation $$\B\sCobtwooo(-) \Rightarrow \B \sCobtwo(-).$$ On the other hand, the sphere $S^2$ is an endomorphism of the empty manifold in $\sCobtwooo(X)$, hence we have a natural map $$\Map(S^2,\mathrm{X})\sslash\Diff^+(S^2) \rightarrow \Omega_{\emptyset} \mathrm{BCob}_2^{SO}(\mathrm{X}),$$ which extends to $$\mathrm{Q}_+( \Map(S^2,\mathrm{X})\sslash\mathrm{SO}(3)) \rightarrow \Omega_{\emptyset} \mathrm{BCob}_2^{SO}(\mathrm{X}),$$ where we use the equivalence $\Diff^+(S^2) \simeq \mathrm{SO}(3)$ from \cite{Smale1959DIFFEOMORPHISMSOT} and $\mathrm{Q}_+=\Omega^{\infty} \Sigma^{\infty}_+$ is the free infinite loop-space functor. Similarly, we have an induced map $$\mathrm{Q}_+(\Map(S^2,\mathrm{X})\sslash\haut^+(S^2)) \rightarrow \Omega_{\emptyset} \mathrm{BCob}_2^{SG}(\mathrm{X}),$$ where $\haut^+(S^2)$ is the monoid of oriented self-homotopy equivalences of $S^2$. These maps assemble into a homotopy commutative square \begin{equation} \label{square-thmB}
\begin{tikzcd}
	{\mathrm{Q}_+ (\mathrm{Map}(S^2,X)\sslash\mathrm{Diff}^+(S^2))} & {\Omega_{\emptyset}\mathrm{BCob}_2^{\mathrm{SO}}(X)} \\
	{\mathrm{Q}_+ (\mathrm{Map}(S^2,X)\sslash\mathrm{haut}^+(S^2))} & {\Omega_{\emptyset} \mathrm{BCob}_2^{\mathrm{SG}}(X)}
	\arrow[from=1-1, to=1-2]
	\arrow[from=1-1, to=2-1]
	\arrow[from=1-2, to=2-2]
	\arrow[from=2-1, to=2-2]
\end{tikzcd}.
\end{equation}

The top horizontal map of the square \eqref{square-thmB} deloops to a map $\susplus \Map(S^2,X)\sslash\Diff^+(S^2) \rightarrow \mathrm{MTSO}(2) \otimes \susplus X$. Let $\mathrm{PH}(2,X)$ denote the pushout of the cospan

\begin{equation} \label{delooping-square}
\begin{tikzcd}
	{\susplus (\mathrm{Map}(S^2,X)\sslash\mathrm{Diff}^+(S^2))} & {\mathrm{MTSO}(2)\otimes \susplus X} \\
	{\susplus (\mathrm{Map}(S^2,X)\sslash\mathrm{haut}^+(S^2))} & {\mathrm{PH}(2,X)}
	\arrow[from=1-1, to=1-2]
	\arrow[from=1-1, to=2-1]
	\arrow[from=1-2, to=2-2]
	\arrow[from=2-1, to=2-2]
\end{tikzcd}.
\end{equation}
 Our next result determines that $\mathrm{PH}(2,-)$ is actually a delooping through the category of spectra $\mathrm{Sp}$ of $\B \sCobtwo(-)$.
 
\begin{mainthm}\label{thmB} 
    For any space $X$, the square \eqref{square-thmB} is a homotopy pullback square. Moreover, a delooping of the square \eqref{square-thmB} in $\mathrm{Sp}$ is given by the square \eqref{delooping-square}.
\end{mainthm}

The computation of the homotopy type of $\Omega_{\emptyset} \B \sCobtwooo(-)$ was motivated in \cite{GMTW} by the study of mapping class groups of surfaces. The functor $\Omega_{\emptyset} \B \sCobtwooo(-)$ is equivalent to the functor $\lop ( \mathrm{MTSO}(2) \otimes (\susplus - ))$. As a corollary, the functor $\Omega_{\emptyset} \B \sCobtwooo(-)$ is excisive. This is surprising because the simplicial levels of the nerve are not excisive. Even though this was not the initial motivation for computing $\Omega_{\emptyset} \B \sCobtwooo(-)$, the functor $\Omega_{\emptyset} \B \sCobtwooo$ being excisive follows from a computation.

We can now ask whether the functor $\Omega_{\emptyset} \B \sCobtwooo(-)$ being excisive depends on the smooth nature of its objects and morphism spaces, or if the functor $\Omega_{\emptyset} \B \sCobtwo(-)$ is also excisive.

To do so, we construct a non-zero obstruction to the map $\Omega \alpha(\star) : \Omega_{\emptyset} \B \sCobtwo \rightarrow \lop \Th(\nuparStwo)$ being an equivalence. We deduce that the functor $\B \sCobtwo(-)$ and $\Pone \B \sCobtwo(-)$ are not equivalent, and then that the functor $\B \sCobtwo(-)$ is not excisive. 

There exists a class $\epsilon \in \Ho^3(\Bhaut_*^+(S^2),\Ztwo)$, which can be interpreted as the first obstruction to lifting a spherical fibration to a vector bundle. The Thom class $U$ of $\nuparStwo$ induces a class $\epsilon.U \in \Ho^1(\Th(\nuparStwo), \Ztwo)$. On cohomology, we have the following morphism: $$\sigma^* : \Ho^*(\Th(\nuparStwo), \Ztwo) \rightarrow \Ho^*(\lop_0 \Th(\nuparStwo), \Ztwo),$$ where $\lop_0$ denotes taking the connected component of the basepoint. Let $$\kappa_{\epsilon} \in \Ho^1(\lop_0\Th(\nuparStwo),\Ztwo)$$ denote the class $\sigma^*(\epsilon. U)$. The following theorem states that the class $\kappa_{\epsilon}$ is the first failure to the map $\Omega \alpha(\star)$ being an equivalence.

\begin{mainthm} \label{thmC}
    The map $$\Omega \alpha(\star) : \Omega_{\emptyset}\B \sCobtwo(\star) \rightarrow \Omega^{\infty} \Th(\nuparStwo)$$  \begin{enumerate}
        \item is a rational equivalence;
        \item induces an isomorphism on $\pi_0$;
    \end{enumerate}
however the class $\kappa_{\epsilon}$ is non-zero and is mapped to $0$ via the pullback morphism $$(\Omega \alpha)^1 : \Ho^1(\lop \Th(\nuparStwo), \Ztwo) \rightarrow \Ho^1(\Omega_{\emptyset} \B \sCobtwo, \Ztwo) .$$ 
\end{mainthm}

As discussed above, as a direct corollary of \cref{thmC} and \cref{crit-ex} we have: 

\begin{cor} \label{COR-thmC}
    The functor $\B \sCobtwo(-)$ is not excisive. 
\end{cor}

\subsection{Outline of the paper}

First, \cref{2} of this paper is devoted to the definition and the construction of a model of the $d$-dimensional Poincaré cobordism category $\PCobdG$ (both non-oriented and oriented). We start with some recollections on Poincaré complexes in Subsection \ref{2.1}. In Subsection \ref{2.2}, we describe a simplicial set model of $\B \haut(P)$, for $P$ a finite space. In Subsection \ref{2.3}, we define the cobordism category $\PCobdG$ as a category enriched in the category $\mathrm{sSet}$ of simplicial sets and write a functor from the smooth cobordism category $\PCobdO$ to $\PCobdG$. In Subsection \ref{2.4}, we prove that the nerve of the Poincaré cobordism category $\B \PCobdG$ is actually an infinite loop space.

In Section \ref{3}, we give a proof of Theorem \ref{thmB}. We begin with comparing diffeomorphisms and self-homotopy equivalences of surfaces in Subsection \ref{3.1bis}. In Subsection \ref{subsec32}, we introduce a reduced cobordism category $\cobred$, obtained from $\sCobtwo$ by deleting spherical components in morphisms. We show that the proof of Theorem \ref{thmB} is equivalent to identifying the fiber of a certain reduction functor $\B \mathrm{red^{SG}} : \B\sCobtwo \rightarrow \B \cobred$. We do the latter in Subsection \ref{subsec33} by using a version of Quillen's Theorem B established in \cite{Steinebrunner_2020} for enriched categories. 

In \cref{4}, we define for $P$ a Poincaré complex and the universal fibration $P\sslash\haut(P) \rightarrow \B \haut(P)$ a spherical fibration $\nupar$ over the total space $P\sslash\haut(P)$, as well as a map $\mathrm{PT}_P^{\haut}: \susplus \B \haut(P) \rightarrow \Th(\nupar)$. Any Poincaré complex $P$ admits a spherical fibration $\nu_P$ (called the Spivak fibration) and a Pontryagin-Thom collapse map $\bS \rightarrow \Th(\nu_P)$. Intuitively, the spherical fibration $\nupar$ and the map $\mathrm{PT}_P^{\haut}$ are families of Spivak fibrations and Pontryagin-Thom collapse maps of the fibers of the universal fibration. Applying these constructions to $S^2$ gives the spherical fibration $\nuparStwo$ appearing in the Statement of \cref{thmA}.

\cref{5} aims to recall the sufficient amount of Goodwillie calculus to prove \cref{thmA} in \cref{6.2}. In Subsection \ref{5.1}, we discuss the classification of excisive functors from the $\infty$-category of spaces $\mathcal{S}$ to spectra $\mathrm{Sp}$. In Subsection \ref{5.2}, we give a recipe to compute the first polynomial approximation $\Pone^*\mathrm{F}$ of a functor $\mathrm{F}: \mathcal{S} \rightarrow \mathrm{Sp}$. Lastly in Subsection \ref{5.3}, we compute the first Goodwillie approximation of the functor $\mathrm{F}_P^{\haut(P)} : \mathrm{X} \rightarrow \Sigma^{\infty}_+ \Map(P,\mathrm{X})\sslash\haut(P)$, where $P$ is a Poincaré complex of dimension $d$. The main upshot, given by \cref{main-prop-P1}, is that the first derivative of $\mathrm{F}_P^{\haut(P)}$ is equivalent to the Thom spectrum $\Th(\nupar)$ of \cref{4}.

Finally in the last \cref{6}, we prove the last two theorems: \cref{thmA} and \cref{thmC}. In Subsection \ref{6.2}, we use the results of \cref{5} and the pushout decomposition given in Theorem \ref{thmB} to determine the first polynomial approximation of $\B \sCobtwo(-)$. We finally prove \cref{thmC} in Subsection \ref{6.3}  by computing the ranks of $\Ztwo$-cohomology groups of $\mathrm{PH}(2,*)$ and $\Th(\nuparStwo)$.

\subsection*{Acknowledgements}

I am extremely grateful to my advisor Søren Galatius for suggesting this project, the many enlightening discussions, suggesting the model of $\B\haut(P)$ proved in \cref{2.2}, as well as reading manifold drafts. I would also like to thank Isaac Moselle for helpful discussions which led to \cref{5}. I would like to thank Carlos Andrés Alvarado Álvarez, Jonathan Clivio, Fadi Mezher, Isaac Moselle, Jonathan Sejr Pedersen and Jan Steinebrunner for reading some parts of the paper. Finally, I thank Siddhi Krishna for the wonderful Inkscape tutorial as well as Florian Riedel for helping me with \LaTeX.

I was supported by the Danish National Research Foundation DNRF151 through the GeoTop center in Copenhagen. Part of this work was written during visits at Columbia University. I would like to thank the mathematics department of Columbia for their hospitality.

\tableofcontents

\section{Constructing a Poincaré Cobordism Category} \label{2}

In this section, we aim to define the cobordism category $\CobdG(X)$ announced in the introduction. To simplify the discussion, assume $X$ is a point. We could try to define it as a topologically enriched category, with objects homotopy classes of $(d-1)$-dimensional Poincaré complexes and morphisms spaces are \[\bigsqcup_{(Q,P_0,P_1)} \B\haut_{\partial}(Q,P_0 \sqcup P_1),\] where the disjoint union is taken over homotopy classes of Poincaré pairs $(Q,P_0\sqcup P_1)$. The composition would then be induced by the union of Poincaré pairs along the common boundary. Furthermore, since homotopy types of smooth manifolds are Poincaré complexes and diffeomorphisms are self-equivalences, we would like to write a functor $\mathrm{Cob}_d^{\mathrm{O}} \rightarrow \mathrm{Cob}_d^{\mathrm{G}}$, where $\mathrm{Cob}_d^{\mathrm{O}}$ is the smooth cobordism category. We could try to define it as a topologically enriched category as above. However, for both categories the composition may not be strictly associative. To circumvent this issue, in \cite{grw9}, the authors upgrade the sets of objects and morphisms to spaces such that both objects and morphisms are embedded subsets of a high-dimensional Euclidean space $\bR^n$. 

In this paper, we define $\CobdG$ as a category internal to simplicial sets, i.e. a category with a simplicial set of objects $(\CobdG)_0$ and a simplicial set of morphisms $(\CobdG)_1$. The elements of the set of $0$-simplices of its objects $(\CobdG)_0$ are certain subsets of $\bR^n$ which are homotopy equivalent to a $(d-1)$-dimensional Poincaré complex. More precisely, a subset $U \subset \bR^n$ is a $0$-simplex of $(\CobdG)_0$ if $U$ is an open submanifold of $\bR^n$, such that $U$ is diffeomorphic to the interior of a compact manifold and $U$ has the homotopy type of a $(d-1)$-dimensional Poincaré complex. Similarly, $0$-simplices of the simplicial set $(\CobdG)_1$ of morphisms are subsets of $\bR^n$ homotopy equivalent to $d$-dimensional Poincaré pairs, with prescribed boundary with respect to the objects. The passage from the Euclidean dimension $n$ to dimension $n+1$ is through crossing with $\bR$. Examples of $0$-simplices of objects and morphisms are illustrated on Figures \ref{fig:ex-simplices-objects-0} and \ref{fig:morphisms-cob-cat}.

The first subsection deals with Poincaré complexes. The second subsection aims at replacing $\B\haut_{\partial}(Q,P)$ by some equivalent simplicial set model. In subsection \ref{2.3}, we define the Poincaré cobordism category $\CobdG$ and write a map functor $ \mathrm{Cob}_d^{\mathrm{O}} \rightarrow \CobdG$. Finally, in subsection \ref{2.4}, we show the nerve $\B\CobdG$ has an infinite loopspace structure. 

\subsection{Recollections on Poincaré Complexes} \label{2.1}

In this subsection, we give the necessary background on Poincaré complexes and Poincaré pairs, see \cite{Wall1967PoincareCI} and \cite{land2021reducibilitylowdimensionalpoincare} for references.
\begin{definition} \label{definition-2.1.1}
    Let $P$ be a connected finite CW-complex. We say that $P$ is a Poincaré complex or Poincaré Duality space of dimension $d$ if there exists a local coefficient system $\mathcal{L}$ on $P$ and a fundamental class $[P]\in \Ho_d(P,\mathcal{L})$ such that $\mathcal{L}$ is pointwise isomorphic to $\bZ$ and the morphism $$- \cap [P] : \Ho^*(P,\mathcal{M}) \rightarrow \Ho_{d-*}(P, \mathcal{M \otimes L})$$ is an isomorphism for all local systems $\mathcal{M}$ on $P$.\\ 
    We say that $P$ is orientable if $\mathcal{L}$ is isomorphic to the constant local system $\underline{\bZ}$. An orientation of $P$ is the choice of an isomorphism $\mathcal{L} \rightarrow \underline{\bZ}$.\\
    If $P$ has a finite number of connected components, we say that $P$ is a Poincaré complex of dimension $d$ if each one of its connected component is a Poincaré complex of dimension $d$.
\end{definition}

As manifolds with boundary work as a relative notion of manifolds, we now introduce Poincaré pair as a relative notion of Poincaré complexes.

\begin{definition}
    Let $(Q,P)$ be a finite CW pair. We say $(Q,P)$ is a Poincaré pair of dimension $d$ if there exists a coefficient system $\mathcal{L}$ on $\mathrm{Q}$ and a class  $[Q] \in \Ho_d(Q,P, \mathcal{L})$ such that $\mathcal{L}$ is pointwise isomorphic to $\bZ$ and such that the morphism $$\cap [Q] : \Ho^*(Q; \mathcal{M}) \rightarrow \Ho_{d-*}(Q,P; \mathcal{L\otimes M})$$ is an isomorphism for all coefficient system $\mathcal{M}$. The induced class $\partial_*([Q])\in \Ho_{d-1}(P, i^*\mathcal{L})$ makes $P$ into a Poincaré complex of dimension $(d-1)$, where $i : \mathrm{P} \rightarrow \mathrm{Q}$ is the inclusion. We say that $(Q,P)$ is orientable if $\mathcal{L}$ is isomorphic to the trivial coefficient system, with an orientation being the choice of such an isomorphism.\\
    Let $P_0,P_1$ be two $(d-1)$-dimensional Poincaré duality spaces. We say that $P_0$ is Poincaré bordant to $P_1$ if there exists a Poincaré duality pair $(Q,P_0\sqcup P_1)$.\\
    Let $(P_j,\mathcal{L}_j, [P_j], o_j)_{j=0,1}$ be two oriented $(d-1)$-dimensional Poincaré Duality spaces where $o_i : \mathcal{L}_i \rightarrow \underline{\bZ}$ denote the choice of the orientation on $P_i$. We say $P_0$ is cobordant to $P_1$ if there exists an oriented Poincaré duality pair $((Q,P_0 \sqcup P_1),\mathcal{L}, [Q], o)$ such that $(P_0,i_0^*\mathcal{L}, i_0^*o)$ inherits the same orientation $o_0$ and $(P_1,i_1^*\mathcal{L}, i_1^*o)$ has the opposite orientation of $(P_1,\mathcal{L}_1,o_1)$. 
\end{definition}

Spivak showed in \cite{SPIVAK196777} any Poincaré duality space admits a canonical pair consisting of a stable spherical fibration which plays the role of a stable normal bundle of a manifold, and a map which plays the role of the Pontryagin-Thom construction: 

\begin{thm}[Spivak Normal Fibration] \label{thm-2.1.3}
    Let $\mathrm{P}$ be a Poincaré complex of dimension $\mathrm{d}$. There exists a couple $$(\xi, \mathrm{c})$$ where $\xi : \mathrm{P} \rightarrow \mathrm{Pic}(\bS)$ is a stable spherical fibration of rank $(-\mathrm{d})$ and $$\mathrm{c} : \mathbb{S} \rightarrow \Th(\xi)$$ is a collapse map such that the collapse map $\mathrm{c}\in \pi_0(\Th(\xi))$ is sent to the fundamental class $[\mathrm{P}]\in \Ho_d(\mathrm{P},\mathcal{L})$ through the map 
\[\begin{tikzcd}
	{\pi_0(\Th(\xi))} & {\Ho_0(\Th(\xi),\mathbb{Z})} & {\Ho_d(P,\mathcal{L})} & {.}
	\arrow["{h_0}", from=1-1, to=1-2]
	\arrow["\cong", from=1-2, to=1-3]
\end{tikzcd}\]
    Moreover, the couple $(\xi, c)$ is unique up to equivalence.\\
    If $P$ is oriented, then $\xi$ is an oriented spherical fibration. 
\end{thm}

In \cref{4}, we explain a reformulation of Poincaré duality in terms of parametrized spectra, as expounded by Land in \cite{land2021reducibilitylowdimensionalpoincare}, following from \cite{Klein_2007}. 

If $P_0, P_1$ are two Poincaré complexes of dimension $d$, then their disjoint union $P_0 \sqcup P_1$ is again a Poincaré complex of dimension $d$. If $(Q,P_0,P_1)$ and $(Q',P_1,P_2)$ are two Poincaré cobordisms then taking the union $(Q\cup_{P_1} Q', P_0,P_2)$ is again a Poincaré cobordism from $P_0$ to $P_2$. If $P$ is a Poincaré complex, the pair $(P\times I, P,P)$ is also Poincaré. Therefore, it makes sense to define Poincaré bordism groups:

\begin{definition} \label{bordism-groups-def}
    Let $\Omega_d^G$ be the bordism classes of $d$-dimensional Poincaré complexes. Disjoint union makes $\Omega_d^G$ into a group, where $\emptyset$ is the unit and every class $[P]$ is of order $2$.
    
    In the same way, we define $\Omega_d^{SG}$ to be the oriented bordism classes of $d$-dimensional oriented Poincaré complexes. In a similar fashion, we can define Poincaré cobordism groups over a space $X$, $\Omega_d^{\cSG}(X)$, of cobordism classes of $d$-dimensional Poincaré complexes equipped with a map to $X$.
\end{definition}

We end this subsection by citing some results on the classification of Poincaré complexes in dimension $1$ and $2$, proved by Wall in \cite{Wall1967PoincareCI} and Eckmann-Müller in \cite{Eckmann1980}. 
For $g,n \geq 0$, let $\Sigma_{g,n}$ denote the genus $g$ surface with $n$ boundary components. 

\begin{thm}[Theorem $4.2$ in \cite{Wall1967PoincareCI}, Corollary $3$ and Theorem $2$ in \cite{Eckmann1980}]
    Let $(Q,P)$ be a connected Poincaré pair of dimension $d$. \begin{enumerate}
        \item if $d=1$, then $(Q,P)$ is equivalent to $(S^1,\emptyset)$ or $(D^1,S^0)$;
        \item if $d=2$, $P=\emptyset$ and if $Q$ is orientable, then $Q$ is equivalent to $\Sigma_g$ for some $g\geq 0$;
        \item if $d=2$, $Q$ orientable, $P\neq \emptyset$ and $\pi_1(Q)$ is finite, then $(Q,P)$ is equivalent to $(D^2,S^1)$;
        \item if $d=2$, $Q$ orientable and $P\neq \emptyset$, then $(Q,P)$ is homotopy equivalent to $(\Sigma_{g,n},\partial \Sigma_{g,n}).$
    \end{enumerate}
\end{thm}

\subsection{A pointset model for \texorpdfstring{$\B\haut_{\partial}(Q,P)$}{B haut partial(Q,P)}} \label{2.2}

As already mentioned in the introduction to Section \ref{2}, we wish to describe a point-set model of the space $\Map_f(Q,X)\sslash\haut_{\partial}(Q,P)$, where $(Q,P)$ is a Poincaré pair. In order to facilitate the preliminary discussion, we first discuss the case without boundary and $X$ is a point. For a smooth closed manifold $M$, its moduli space $\B\Diff(M)$ is equivalent to the space of submanifolds of $\bR^{\infty}$ which are diffeomorphic to $M$. In a similar flavor, for $P$ a Poincaré complex, we wish to model $\B\haut(P)$ by a space whose points are subsets of $\bR^{\infty}$ homotopy equivalent to $P$. The idea is to replace $P$ by an $n$-dimensional open manifold $U$, such that $U$ is diffeomorphic to the interior of a compact smooth thickening $N \subset \bR^n$ of $P$. Then, we show $\B\haut(P)$ is equivalent to a space whose points are subsets of $A \subset \bR^{n+k}$ which are diffeomorphic to $U\times \bR^k$, where we identify $A \subset \bR^{n+k}$ with $A\times \bR \subset \bR^{n+k+1}$ after taking the direct limit $k \rightarrow \infty$.

Before going any further, let us introduce a few notations. Let $M,N$ be two manifolds. Let $$S_{\bullet}\mathrm{Sub}_{\partial}(M,N)$$ denote the simplicial set, whose $k$-simplices are $\Delta^k$-parametrized families of subsets of $N$ diffeomorphic to $M$ relative boundary. For $X$ a space, let $$S_{\bullet}\mathrm{Sub}_{\partial}(M,N)_{/X}$$ denote the simplicial set with $k$-simplices pairs $(A,f)$, where $A$ is a $k$-simplex of $S_{\bullet} \mathrm{Sub}_{\partial} (M,N)$ and $f$ is a $\Delta^k$-parametrized map from $A$ to $X$. If $N$ is the Euclidean space $\bR^n$, we write $S_{\bullet} \mathrm{Sub}_{\partial}(M,n)_{/X}$ instead of $S_{\bullet} \mathrm{Sub}_{\partial}(M,\bR^n)_{/X}$. Crossing with $\bR$ induces a map $$S_{\bullet} \mathrm{Sub}_{\partial}(M,n)_{/X} \rightarrow S_{\bullet}\mathrm{Sub}_{\partial}(M\times \bR,n+1)_{/X}.$$ Let $S_{\bullet} \mathrm{Sub}_{\partial}(M,\infty)_{/X}$ denote the direct limit of the $S_{\bullet} \mathrm{Sub}_{\partial}(M,n)_{/X}$.

In this subsection, we aim to prove the following proposition. 

\begin{prop} \label{prop-approx-Bhaut}
    Let $(Q,P)$ be a Poincaré pair of dimension $d$ and let $f : P \rightarrow X$ be a map. Let $N$ be a smooth compact manifold of dimension $n$, such that $N$ can be embedded in $\bR^n$. Let $N_0$ be a codimension $0$ compact submanifold of $\partial N$, such that $N_0$ can be embedded in $\bR^{n-1}$. Assume $(N,N_0)$ deformation retracts to the pair $(Q,P)$. Let $(U_1,U_0)$ be the open manifold with boundary $(N - (\partial N - \mathrm{int}(N_0)), \mathrm{int}(N_0))$. Then, the simplicial set $$S_{\bullet} \mathrm{Sub}_{\partial}(U,\infty)_{/X}$$ is equivalent after geometric realization to $$\Map_f(Q,X)\sslash\haut_{\partial}(Q,P).$$
\end{prop}

Before proving Proposition \ref{prop-approx-Bhaut}, we recall suitable simplical set models of spaces of embeddings and diffeomorphisms of manifolds, following \cite{alma991086073420406532}. In a second phase, we discuss thickenings of Poincaré complexes, before finally giving a proof of Proposition \ref{prop-approx-Bhaut}. The latter involves embedding calculus to compare the $S_{\bullet} \mathrm{Sub}_{\partial}(U,n)$ with quotients of spaces of bundle maps, which have an entirely homotopy-theoretic description.

Let $P$ be a Poincaré Complex. Let $\haut(P)$ denote the grouplike monoid of self-homotopy equivalences of $P$, topologized as a subspace of $\Map(P,P)$. If $(Q,P)$ is a Poincaré duality pair, we consider the group-like monoid $\haut_{\partial}(Q,P)$ of self-equivalences $f$ such that $f_{|P}$ coincides with the inclusion $P \subset Q$. For $(Q,P)$ an oriented Poincaré Duality pair, we denote by $\haut_{\partial}^+(Q,P)$ the orientation-preserving self-homotopy equivalences of $Q$ relative $P$. Let $f\in \Map(P,X)$. Let $\Map_f(Q,X)$ denote the space of maps $\phi\in \Map(Q,X)$ such that $\phi_{|P}$ coincides with $f$. For any space $P$, the classifying space $\Bhaut(P)$ classifies Serre fibrations with fiber equivalent to $P$. For a pair $(Q,P)$, $\Bhaut_{\partial}(Q,P)$ classifies relative Serre fibrations with fiber equivalent to $(Q,P)$.

Let $M,N$ be smooth manifolds and $e : \partial M \hookrightarrow N$ an embedding. Let $(Q,P)$ be a pair of spaces and $f : P \rightarrow X$ be a map to a space $X$.
We topologize the diffeomorphism group $\Diff_{\partial}(M)$ as the geometric realization of a simplicial group $\mathrm{S}_{\bullet}\Diff_{\partial}(M)$ with set of $k$-simplices:

\[
\mathrm{S}_k \Diff_{\partial}(M) = \Biggl\{ \begin{tikzcd}[row sep=tiny]
	{M\times \Delta^k} & {} & {M\times \Delta^k} \\
	& {\Delta^k}
	\arrow["\cong", from=1-1, to=1-3]
	\arrow["\pi"', from=1-1, to=2-2]
	\arrow[from=1-2, to=1-3]
	\arrow["\pi", from=1-3, to=2-2]
\end{tikzcd},\ \text{fixing}\  \partial M \times \Delta^k \Biggr\}.
\]
We topologize the embedding space $\Emb_e(M,N)$ as the geometric realization of the simplicial set $\mathrm{S}_{\bullet} \Emb_e(M,N)$ with $k$-simplices as follows: 
\[
\mathrm{S}_k \Emb_e(M,N) = \Biggl\{ \begin{tikzcd}
	{M\times \Delta^k} & {} & {N\times \Delta^k} \\
	& {\Delta^k}
	\arrow[hook, from=1-1, to=1-3]
	\arrow["\pi"', from=1-1, to=2-2]
	\arrow[from=1-2, to=1-3]
	\arrow["\pi", from=1-3, to=2-2]
\end{tikzcd},\ \text{restricting to}\ e\times \Delta^k\  \text{on}\ \partial M\times \Delta^k \Biggr\}. \]
Let $S_{\bullet}\mathrm{Emb}_e^{\simeq}(M,N)$ denote the subsimplicial set of embeddings $S_{\bullet} \mathrm{Emb}_e(M,N)$ which are equivalences.
The simplicial group $S_{\bullet} \Diff_{\partial}(M)$ acts levelwise and freely on $S_{\bullet} \mathrm{Emb}_e(M,N)$. We can then define the simplicial set $S_{\bullet} \mathrm{Sub}_{\partial}(M,N)$ to be the levelwise quotient simplicial set $$S_{\bullet} \mathrm{Emb}_e(M,N)/S_{\bullet} \Diff_{\partial}(M).$$

Let $\mathrm{S}_{\bullet} \Map_f(Q,X)$ be the simplicial set with $k$-simplices:
\[ \mathrm{S}_k \Map_f(Q,X) =\Biggl\{ \begin{tikzcd}
	{Q\times \Delta^k} & {} & {X\times \Delta^k} \\
	& {\Delta^k}
	\arrow[from=1-1, to=1-3]
	\arrow["\pi"', from=1-1, to=2-2]
	\arrow[from=1-2, to=1-3]
	\arrow["\pi", from=1-3, to=2-2]
\end{tikzcd},\ \text{restricting to}\ f\times \Delta^k\ \text{on}\ P\times \Delta^k \Biggl\}. \] We observe it is the singular complex of $\Map_f(Q,X)$, hence is a Kan complex.

The following is shown in \cite{alma991086073420406532} and \cite{maysimp}:

\begin{lem}[Proposition $2.5$ in \cite{alma991086073420406532}, Theorem $17.1$ in \cite{maysimp}] \label{lem-2.2-0.2}
The simplicial sets $$S_{\bullet} \Emb_{e}(M,N)$$ and $$S_{\bullet} \Diff_{\partial}(M)$$ are Kan complexes.     
\end{lem}

In what follows, we make extensive use of the expression "level-preserving" or "level-wise preserving" map. To avoid confusions, we clarify below what we mean:

\begin{definition}
    Let $U,V$ be two subsets of $\Delta^k \times \bR^d$. We say a map $\phi : U \rightarrow V$ is a level-preserving map/embedding/diffeomorphism if it is a map/embedding/diffeomorphism and if $\phi$ commutes with the projection $$\pi : \Delta^k \times \bR^n  \rightarrow \Delta^k.$$
\end{definition}

We now briefly discuss smooth thickenings of Poincaré pairs. Let $(Q,P)$ be a Poincaré duality pair. In particular, it is equivalent to a finite pair of CW-complexes. Following \cite{SPIVAK196777} or \cite{Browder1972}, for $k$ large enough, we can find an embedding $e_0 : P \hookrightarrow \bR^{k-1}$, i.e. an injective map which is a homeomorphism on its image. Up to replacing $Q$ with the mapping cylinder of the inclusion $P \hookrightarrow Q$, we can assume $P$ admits a collar neighborhood $P \times [0,1] \hookrightarrow Q$. We can then find an embedding $e : Q \hookrightarrow \bR^k$ such that $e$ restricts to $e_0 \times \text{id}$ on the collar $P \times [0,1]$. The space $Q$ embedded in $\bR^k$ can be thickened to a compact smooth submanifold with boundary $N \subset \bR^k$ such that $N_0= N \cap \bR^{k-1}$ is a compact thickening of $P$. Actually $N_0$ is a codimension $0$ submanifold of $\partial N$ and defines a manifold triad $(N,\partial N, N_0)$. We say such a triad is a \emph{relative thickening} of $(Q,P)$. On the other hand, we can take the pair $(N-N_1,\mathrm{int}(N_0))$, where $N_1$ is the manifold $\partial N \setminus \mathrm{int}(N_0)$ . We say it is a \emph{relative open thickening} of $(Q,P)$.

We can now begin the proof of Proposition \ref{prop-approx-Bhaut}. We break down the proof in several steps. For $M,N$ two manifolds and $e_0 : \partial M \hookrightarrow N$ an embedding, taking the derivative induces a map $$\Emb_{e_0}(M,N) \rightarrow \mathrm{Bun}_{Te_0}(TM,TN).$$  Here, $\mathrm{Bun}_{Te_0}(TM,TN)$ denotes the space of bundle maps from the tangent bundle $TM$ to $TN$ restricting to $Te_0$ on the boundary. In the following proposition, we use embedding calculus to show the derivative maps are highly-connected for thickenings.

\begin{prop} \label{emb-calculus}
   Let $(Q,P)$ be a Poincaré pair of dimension $d$. Let $(N,\partial N, N_0) \subset \bR^k$ be a relative thickening of $(Q,P)$ and $e_0 : N_0 \hookrightarrow \bR^k$ be a fixed embedding. Let $(\mathrm{U}_1,\mathrm{U}_0)$ be $(N-(\partial N-\partial_1 N), \text{int}(N_0))$.
  There exists a constant $c$ such that the derivative maps $$\Emb_{e_0}(\mathrm{U}_1,\bR^k) \rightarrow \mathrm{Bun}_{Te_0}(T\mathrm{U}_1,T\bR^k)$$ and $$\Emb^{\simeq}_{\partial}(\mathrm{U}_1,\mathrm{U}_1) \rightarrow \mathrm{Bun}^{\simeq}_{\partial}(T\mathrm{U}_1,T\mathrm{U}_1)$$ are $(\mathrm{k-d+c})$-connected.
\end{prop}

\begin{proof}  
The handle dimension of the pair $(\mathrm{U}_1,\mathrm{U}_0)$ is majorized by $d$. The handle dimension of a pair $(\mathrm{U}_1,\mathrm{U}_0)$ is invariant under crossing with $\bR$. In particular, up to replacing $(\mathrm{U}_1,\mathrm{U}_0)$ by the relative open thickening $(\mathrm{U}_1\times \bR^l, \mathrm{U}_0\times \bR^l)$ of $(Q,P)$ for $l$ sufficiently large, we can assume: $\text{hdim}(\mathrm{U}_1,\mathrm{U}_0)\leq k-3.$ 
We are in the situation for convergence of the embedding calculus tower. According to \cite{Goodwillie_1999}, there exists a constant $c$ such that the maps $$T_l\Emb_{e_0}(\mathrm{U}_1,\bR^k) \rightarrow T_{l-1}\Emb_{e_0}(\mathrm{U}_1,\bR^k)$$ are at least $l(k-\text{hdim}(\mathrm{U}_1,\mathrm{U}_0)+c)$-connected. In particular, since the approximation map $$\Emb_{e_0}(\mathrm{U}_1,\bR^k) \rightarrow T_{\infty}\Emb_{e_0}(\mathrm{U}_1,\bR^k)$$ is an equivalence, we deduce the map $$\Emb_{e_0}(\mathrm{U}_1,\bR^k) \rightarrow T_1\Emb_{e_0}(\mathrm{U}_1,\bR^k)$$ is at least $(k-d+c)$-connected. According to \cite{Goodwillie_1999}, the first stage in the embedding tower is given by $\mathrm{Bun}_{Te_0}(T\mathrm{U}_1,T\bR^k)$ and the derivative map is the approximation map.

Similarly, up to replacing $(\mathrm{U}_1,\mathrm{U}_0)$ by $(\mathrm{U}_1\times \bR^l, \mathrm{U}_0\times \bR^l)$ for $l$ large, we can assume $$\text{dim}(\mathrm{U}_1)-\text{hdim}(\mathrm{U}_1,\mathrm{U}_0) \geq k -d \geq 3.$$ Again, the embedding calculus tower for $\Emb_{\partial}(\mathrm{U}_1,\mathrm{U}_1)$ converges. There exists a constant $c$ independent of $k$ such that the map $d_{\mathrm{U}_1} : \Emb_{\partial}(\mathrm{U}_1,\mathrm{U}_1) \rightarrow \mathrm{Bun}_{\partial}(T\mathrm{U}_1,T\mathrm{U}_1)$ is $(k-d+c)$-connected. The spaces $\Emb^{\simeq}_{\partial}(\mathrm{U}_1,\mathrm{U}_1)$ and $\mathrm{Bun}_{\partial}^{\simeq}(T\mathrm{U}_1,T\mathrm{U}_1)$ are respectively obtained from $\Emb_{\partial}(\mathrm{U}_1,\mathrm{U}_1)$ and $\mathrm{Bun}_{\partial}(T\mathrm{U}_1,T\mathrm{U}_1)$ by restricting to path-components which are invertible in the monoids $\pi_0(\Emb_{\partial}(\mathrm{U}_1,\mathrm{U}_1))$ and $\pi_0(\mathrm{Bun}_{\partial}(T\mathrm{U}_1,T\mathrm{U}_1))$. For $k$ large enough, $k-d+c \geq 1$. Then, the map $d_{\mathrm{U}_1}$ induces an isomorphism on $\pi_0$. Thus, the map $\pi_0(\Emb_{\partial}(\mathrm{U}_1,\mathrm{U}_1))^{\times} \rightarrow \pi_0(\mathrm{Bun}_{\partial}(T\mathrm{U}_1,T\mathrm{U}_1))^{\times}$ is an isomorphism on invertible elements. We deduce the map $$\Emb_{\partial}^{\simeq}(\mathrm{U}_1,\mathrm{U}_1) \rightarrow \mathrm{Bun}_{\partial}^{\simeq}(T\mathrm{U}_1,T\mathrm{U}_1)$$ is also $(k-d+c)$-connected.  
\end{proof}

\begin{lem} \label{diffeomorphisms-lem 2.16}
    Let $(\mathrm{U}_1,\mathrm{U}_0) \subset (\bR^k,\bR^{k-1})$ be a relative open thickening. The monoid map $$\Diff_{\partial}(\mathrm{U}_1) \rightarrow \Emb_{\partial}^{\simeq}(\mathrm{U}_1,\mathrm{U}_1)$$ is a homotopy equivalence.
\end{lem}

\begin{figure}[H]
\labellist  
\pinlabel {\textcolor{mypurple}{$U_1$}} at 800 500 
\pinlabel {\textcolor{myorange}{$\partial_1N\times [0,\infty)$}} at 1800 400 
\pinlabel {\textcolor{myorange}{$\partial_1 N$}} at 1250 250 

\pinlabel {\textcolor{mypurple}{$U_0$}} at 290 250 
\endlabellist
\includegraphics[height=5cm, width=8cm]{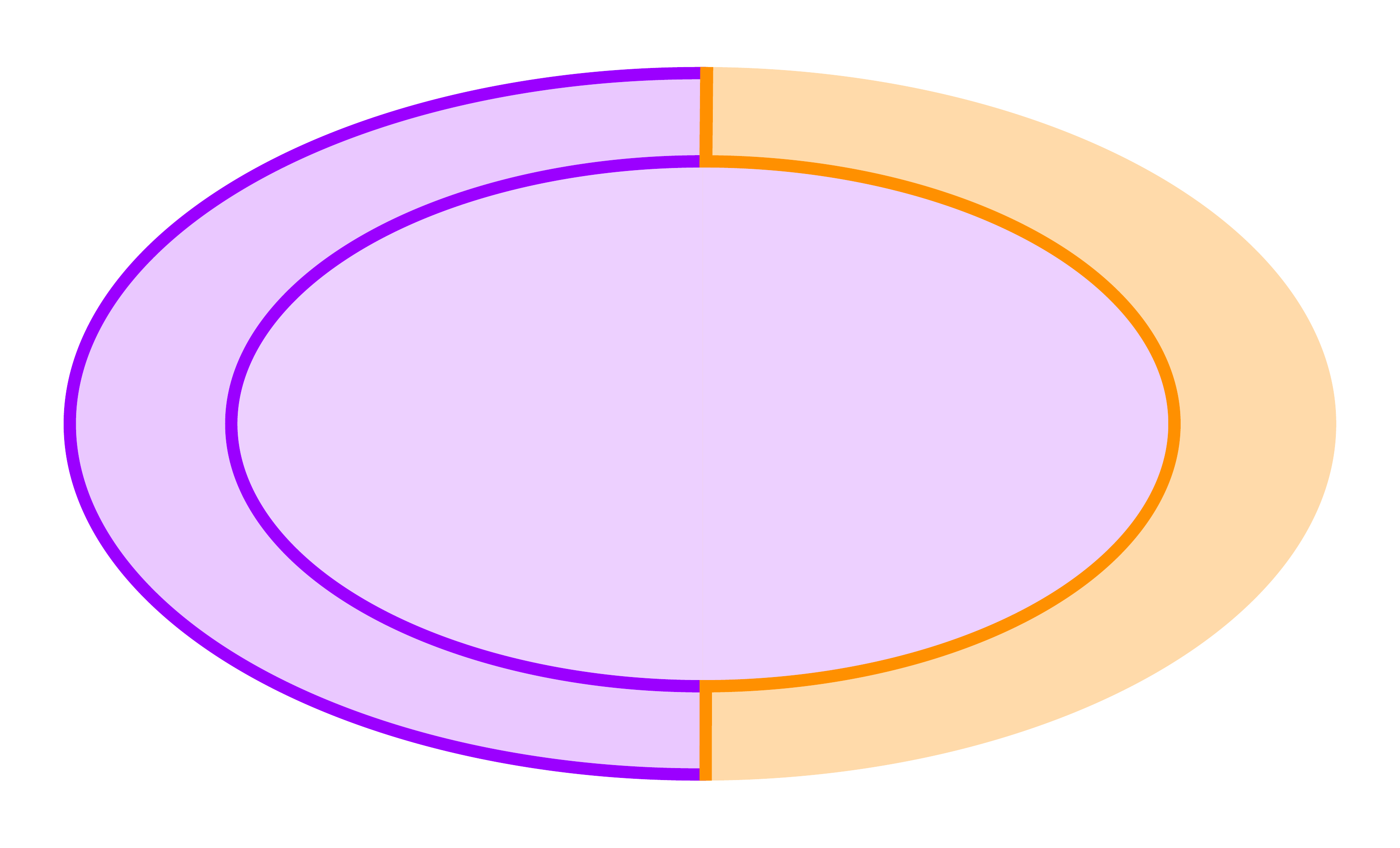}
\centering
\captionof{figure}{The relative thickening $(U_1,U_0)$}
\label{fig:for-proof-model-Bhaut}
\end{figure}
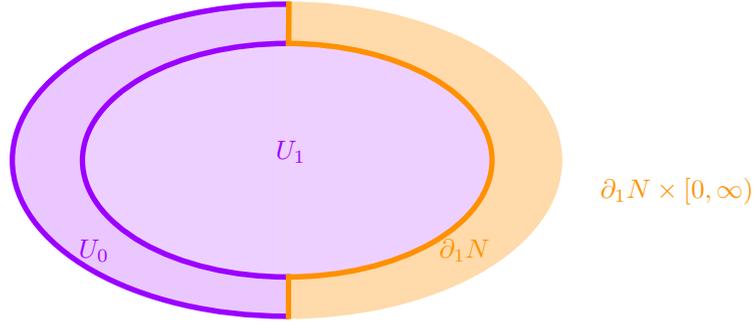

\begin{proof}
Let $(N,\partial N,\partial_0 N)$ be the manifold triad such that $\mathrm{U}_0$ is the interior of $\partial_0 N$ and $\mathrm{U}_1$ is obtained by taking $(N-\partial_1 N)$. Here, $\partial_1 N$ is such that $\partial N =\partial_0 N \cup_{\partial_{01} N} \partial_1 N$. By taking a collar of $\partial N$, we obtain inverse up to isotopy embeddings $j_1 : N \hookrightarrow \mathrm{U}_1$ and $\mathrm{U}_1 \hookrightarrow N$. By taking a collar of $\partial_0 N$, we also obtain inverse up to isotopy embeddings $\partial_0 N \hookrightarrow \mathrm{U}_0$ and $\mathrm{U}_0 \hookrightarrow \partial_0 N$. Consequently, the restriction map $$r_{|j_1} : \Emb_{\partial}(\mathrm{U}_1,\mathrm{U}_1) \rightarrow \Emb_{\partial_0 N}(N,\mathrm{U}_1)$$ is an equivalence.

The fiber at $j_1$ of the restriction map $$\Diff_{\partial}(\mathrm{U}_1) \rightarrow \Emb_{\partial_0 N}(N,\mathrm{U}_1)$$
is equivalent to $\Diff_{\partial}(\partial_1 N \times [0,\infty))$ which is contractible, as illustrated on Figure \ref{fig:for-proof-model-Bhaut}.

\end{proof}

We need the following lemma from \cite{bonatto2020decouplingdecorationsmodulispaces}. The statement was written originally for topological groups, but it is not too hard to adapt the proof to group-like monoids. 

\begin{lem}[Corollary 2.11 in \cite{bonatto2020decouplingdecorationsmodulispaces}] \label{lem-2.17}
Let $G_i$ be group-like monoids and $S_i$ be $G_i$-spaces for $i=1,2,3$. Assume there is a short exact sequence 
\[\begin{tikzcd}
	1 & {G_1} & {G_2} & {G_3} & 1
	\arrow[from=1-1, to=1-2]
	\arrow[from=1-2, to=1-3]
	\arrow["\phi", from=1-3, to=1-4]
	\arrow[from=1-4, to=1-5]
\end{tikzcd}\] such that $\phi$ is a principal $G_1$-bundle. Let $$S_1 \rightarrow S_2 \rightarrow S_3$$ be a fiber sequence of equivariant maps. Then the induced maps on quotients $$S_1\sslash G_1 \rightarrow S_2\sslash G_2 \rightarrow S_3\sslash G_3$$ form a fiber sequence.
\end{lem}

We now use \cref{lem-2.17} to identify $\Bhaut_{\partial}(Q,P)$ with a homotopy quotient of spaces of bundle maps. Let $(Q,P)$ be a Poincaré pair and let $(\mathrm{U}_1,\mathrm{U}_0) \subset (\bR^k,\bR^{k-1})$ be a relative open thickening of $(Q,P)$. The space of bundle maps $\mathrm{Bun}_{\partial}^{\simeq}(T\mathrm{U}_1,T\mathrm{U}_1)$ acts by precomposition on $\mathrm{Bun}_{Te_0}(T\mathrm{U}_1,T\bR^k)$. Furthermore, there is a forgetful monoid map $$\mathrm{Bun}_{\partial}^{\simeq}(T\mathrm{U}_1,T\mathrm{U}_1) \rightarrow \haut_{\partial}(\mathrm{U}_1,\mathrm{U}_0).$$ Sending $\mathrm{Bun}_{Te_0}(T\mathrm{U}_1,T\bR^k)$ to the point induces a map $$\mathrm{Bun}_{Te_0}(T\mathrm{U}_1,T\bR^k)\sslash\mathrm{Bun}_{\partial}^{\simeq}(T\mathrm{U}_1,T\mathrm{U}_1) \rightarrow \Bhaut_{\partial}(\mathrm{U}_1, U_0).$$ 
In the following lemma, we prove this map is an equivalence:

\begin{lem} \label{lem-2.18}
    If $(\mathrm{U}_1,\mathrm{U}_0)$ is a relative open thickening of a Poincaré pair, then the composite $$\mathrm{Bun}_{Te_0}(T\mathrm{U}_1,T\bR^k))\sslash\mathrm{Bun}_{\partial}^{\simeq}(T\mathrm{U}_1,T\mathrm{U}_1) \rightarrow \Bhaut_{\partial}(U_1,U_0)$$ is a weak equivalence.
\end{lem}

\begin{proof}
    Since $\mathrm{U}_1$ is a codimension $0$ submanifold of $\bR^k$, its tangent bundle $T\mathrm{U}_1$ is given by the projection $\mathrm{U}_1 \times \bR^k \rightarrow \mathrm{U}_1$. Let $\Map_{\partial}(\mathrm{U}_1,(\bR))$ denote the space of maps sending pointwise the boundary to $\text{id} \in GL_k(\bR)$. The space of bundle maps $\mathrm{Bun}_{Te_0}(T\mathrm{U}_1,T\bR^k)$ is then exactly given by $\Map_{\partial}(\mathrm{U}_1,\text{GL}_k(\bR)) \times \Map_:{e_0}(\mathrm{U}_1,\bR^k)$. The group structure on $GL_k(\bR)$ makes $\Map_{\partial}(\mathrm{U}_1,\text{GL}_k(\bR))$ into a topological group. We denote the multiplication by $\star$. On the other hand, $\mathrm{Bun}_{\partial}(T\mathrm{U}_1,T\mathrm{U}_1)^{\simeq}$ is given as a monoid by the semi-direct product $\Map_{\partial}(\mathrm{U}_1,GL_k(\bR)) \rtimes \haut_{\partial}(\mathrm{U}_1)$. 
    
    There is a monoid map $\iota : \Map_{\partial}(\mathrm{U}_1,\text{GL}_k(\bR)) \rightarrow \mathrm{Bun}_{\partial}(T\mathrm{U}_1,T\mathrm{U}_1)^{\simeq}$ given by sending a map $\phi$ to the couple $(\text{id}_{\mathrm{U}_1},\phi)$. The projection map $\pi : \mathrm{Bun}_{\partial}(T\mathrm{U}_1,T\mathrm{U}_1)^{\simeq} \rightarrow \haut_{\partial}(\mathrm{U}_1)$ is a monoid map and is a trivial $\Map_{\partial}(\mathrm{U}_1,\text{GL}_k(\bR))$-bundle. The following sequence 
\[\begin{tikzcd}
	1 & {\Map_{\partial}(\mathrm{U}_1,\text{GL}_k(\mathbb{R}))} & {\mathrm{Bun}^{\simeq}_{\partial}(T\mathrm{U}_1,T\mathrm{U}_1)} & {\haut_{\partial}(\mathrm{U}_1)} & 1
	\arrow[from=1-1, to=1-2]
	\arrow["\iota", from=1-2, to=1-3]
	\arrow["\pi", from=1-3, to=1-4]
	\arrow[from=1-4, to=1-5]
\end{tikzcd}\] is short exact. On the other hand, $\haut_{\partial}(\mathrm{U}_1)$ acts on $\Map_{e_0}(\mathrm{U}_1,\bR^k)$. The space of bundle maps $\mathrm{Bun}^{\simeq}(T\mathrm{U}_1,T\mathrm{U}_1)$ acts on $\mathrm{Bun}_{Te_0}(T\mathrm{U}_1,T\bR^k)$ by precomposition. Finally $\Map_{\partial}(\mathrm{U}_1,\text{GL}_k(\bR))$ acts on $\Map_{\partial}(\mathrm{U}_1,\text{GL}_k(\bR))$ by precomposition.

There is a projection map $p : \mathrm{Bun}_{Te_0}(T\mathrm{U}_1,T\bR^k) \rightarrow \Map_{e_0}(\mathrm{U}_1,\bR^k) $ equivariant under the action of $\mathrm{Bun}^{\simeq}_{\partial}(T\mathrm{U}_1,T\mathrm{U}_1)$. Fixing an embedding $e_1 : \mathrm{U}_1 \hookrightarrow \bR^k$, restricting to $e_0$ on the boundary, gives a map $$I : \Map_{\partial}(\mathrm{U}_1,\text{GL}_k(\bR)) \rightarrow \mathrm{Bun}_{Te_0}(T\mathrm{U}_1,T\bR^k)$$ given by sending $\phi \in \Map_{\partial}(\mathrm{U}_1,\text{GL}_k(\bR))$ to $(\phi,e_1)$. It is equivariant with respect to the map $$\iota : \Map_{\partial}(\mathrm{U}_1,GL_k(\bR)) \rightarrow \mathrm{Bun}_{\partial}^{\simeq}(T\mathrm{U}_1,T\mathrm{U}_1).$$ Indeed, for $(\phi,e_1) \in \mathrm{Bun}_{Te_0}(T\mathrm{U}_1,T\bR^k)$, $\psi \in \Map_{\partial}(\mathrm{U}_1,\text{GL}_k(\bR))$, $$I(\phi) \bullet \iota(\psi) = (e_1,\phi) \bullet (\text{id},\psi)= (e_1,\phi \star \psi) =I(\phi \star \psi).$$

Finally, there is a fiber sequence of equivariant maps 
\[\begin{tikzcd}
	{\Map_{\partial}(\mathrm{U}_1,\text{GL}_k(\mathbb{R}))} & {\mathrm{Bun}_{Te_0}(T\mathrm{U}_1,T\mathbb{R}^k)} & {\Map_{e_0}(\mathrm{U}_1,\mathbb{R}^k)}
	\arrow[from=1-1, to=1-2]
	\arrow[from=1-2, to=1-3]
\end{tikzcd}.\]
Note that $\Map_{\partial}(\mathrm{U}_1,GL_k(\bR))\sslash\Map_{\partial}(\mathrm{U}_1,GL_k(\bR))$ and $\Map_{e_0}(\mathrm{U}_1,\bR^k)$ are contractible.

Then, \cref{lem-2.17} allows to conclude $$ * \rightarrow \mathrm{Bun}_{Te_0}(T\mathrm{U}_1,T\bR^k)\sslash\mathrm{Bun}_{\partial}^{\simeq}(T\mathrm{U}_1,T\mathrm{U}_1) \rightarrow \Bhaut_{\partial}(\mathrm{U}_1)$$ is a fiber sequence.

\end{proof}

We may now conclude the proof of \cref{prop-approx-Bhaut}:

\begin{proof}[Proof of Proposition \ref{prop-approx-Bhaut}]
For $l\in \bN$, the pair $(\mathrm{U}_1 \times \bR^l, \mathrm{U}_0\times \bR^l) \subset (\bR^{k+l},\bR^{k+l-1})$ is a relative open thickening of $(Q,P)$.
Composing the maps from \cref{emb-calculus} and \cref{lem-2.18}, we obtain $(l+k-d+c)$-connected maps $$j_l : \Emb_{e_0}(\mathrm{U}_1\times \bR^l ,\bR^{k+l})/ \Diff_{\partial}(\mathrm{U}_1) \rightarrow \Bhaut_{\partial}(\mathrm{U}_1 \times \bR^l, \mathrm{U}_0 \times \bR^l).$$ The maps $j_l$ are compatible with crossing $(\mathrm{U}_1 \times \bR^l, \mathrm{U}_0 \times \bR^l)$ with $\bR$.
We then get a map 

\[ J : \hocolim_{l\rightarrow \infty}\Emb_{e_0}(\mathrm{U}_1\times \bR^l,\bR^{k+l})\sslash\Diff_{\partial}(\mathrm{U}_1) \rightarrow \hocolim_{l \rightarrow \infty} \B\haut_{\partial}(\mathrm{U}_1 \times \bR^l, \mathrm{U}_0 \times \bR^l). \]
We now remark the right handside of the map is equivalent to $\B\haut_{\partial}(Q,P)$. Indeed, $Q$ and $P$ are respectively retracts of $U_1$ and $U_0$, hence $\B\haut_{\partial}(U_1,U_0)$ is equivalent to $\B \haut_{\partial}(Q,P)$.

We now show the map $J$ induces isomorphism on homotopy groups. For $a\in \bN$, take a map  $\alpha : S^a \rightarrow \Bhaut_{\partial}(Q,P)$. For $l$ large enough the map $j_l$ induces an isomorphism on $\pi_a$ hence we can lift  $\alpha$ to a map in $\Emb_{e_0}(\mathrm{U}_1\times \bR^l ,\bR^{k+l})/ \Diff_{\partial}(\mathrm{U}_1)$. On the other hand, let 
\[\alpha : \mathrm{S}^a \rightarrow \hocolim_{l\rightarrow \infty}\Emb_{e_0}(\mathrm{U}_1\times \bR^l,\mathrm{U}_0\times \bR^l)\sslash\Diff_{\partial}(\mathrm{U}_1)\]  be a map such that $J(\alpha)$ is nullhomotopic. By compatness of $S^a$, it factors through some $\Emb_{e_0}(\mathrm{U}_1\times \bR^l ,\bR^{k+l})/ \Diff_{\partial}(\mathrm{U}_1)$. We can choose $l$ large enough so that the map $j_l$ is injective on $\pi_a$. Then $\alpha$ is null in $\pi_a(\Emb_{e_0}(\mathrm{U}_1\times \bR^l ,\bR^{k+l})\sslash \Diff_{\partial}(\mathrm{U}_1))$ hence 
\[\mathrm{J}_* : \pi_*(\hocolim_{l\rightarrow \infty} \Emb_{e_0}(\mathrm{U}_1\times \bR^l,\mathrm{U}_0\times \bR^l)\sslash\Diff_{\partial}(\mathrm{U}_1)) \rightarrow \pi_*(\Bhaut(P)) \] is an isomorphism.
Composition of diffeomorphisms makes $\mathrm{S}_{\bullet} \Diff_{\partial}(\mathrm{U}_1\times \bR^l)$ into a simplicial group. Moreover, it acts freely on $\mathrm{S}_{\bullet} \Emb_{e_0}(\mathrm{U}_1\times \bR^l, \bR^{l+k})$. According to Lemma $2.18$ in \cite{maysimp}, the quotient map $$\mathrm{S}_{\bullet} \Emb_{e_0}(\mathrm{U}_1\times \bR^l, \bR^{l+k}) \rightarrow  S_{\bullet} \mathrm{Sub}_{\partial}(U_1\times \bR^l,l+k)$$ is a Kan fibration. Thus the geometric realization of $S_{\bullet} \mathrm{Sub}_{\partial}(U_1 \times \bR^l, l+k)$  is equivalent to $\Emb_{e_0}(\mathrm{U}_1\times \bR^l, \bR^{k+l})\sslash\Diff_{\partial}(\mathrm{U}_1\times \bR^l)$. The maps $$S_{\bullet} \mathrm{Sub}_{\partial}(U_1\times \bR^l, l+k) \rightarrow S_{\bullet} \mathrm{Sub}_{\partial}(U_1\times \bR^{l+1},l+k+1)$$ are levelwise injective, hence induce cofibrations after geometric realization. Finally the geometric realization of 
 $S_{\bullet} \mathrm{Sub}_{\partial}(U_1,\infty)$
  is equivalent to  
  \[\hocolim_{l \rightarrow \infty} \Emb_{e_0}(\mathrm{U}_1\times \bR^l,\mathrm{U}_0\times \bR^l)\sslash\Diff_{\partial}(\mathrm{U}_1).\]
\end{proof}

We can derive a similar model for the classifying space $\B\haut_{\partial}(Q,P)$ of orientation-preserving self-equivalences of an oriented Poincaré pair $(Q,P)$. Let $S_{\bullet} \mathrm{Sub}_{\partial}^+(U_1,n)$ denote the quotient simplicial set $S_{\bullet} \mathrm{Emb}_{\partial}(U_1,\bR^n)/S_{\bullet} \Diff^+_{\partial}(U_1)$ where $U_1$ is an oriented manifold. Similarly, we define $S_{\bullet}\mathrm{Sub}_{\partial}^+(U_1,\infty)$ and $S_{\bullet}\mathrm{Sub}_{\partial}^+(U_1,\infty)_{/X}$.

\begin{prop} 
    Let $(Q,P)$ be an oriented Poincaré pair. Let $(N,\partial N, N_0)$ be an oriented relative thickening of $(Q,P)$. Let $(\mathrm{U}_1,\mathrm{U}_0)$ be the pair $(N\setminus (\partial N- N_0),\text{int}(N_0))$. Let $e_0 : N_0 \hookrightarrow \mathrm{S}^{k-1}$ be an orientation-preserving embedding extending $e_0 : P \hookrightarrow \mathrm{S}^{k-1}$ and let $f_0 : P \rightarrow X$.\\
    The simplicial set $$S_{\bullet} \mathrm{Sub}^+_{\partial}(U_1,\infty)_{/X}$$
    is equivalent after geometric realization to $$\Map_{f_0}(Q,X)\sslash\haut^+_{\partial}(Q,P).$$
\end{prop}

We now end this subsection with one remark:

\begin{rk}
    In the proof of \cref{prop-approx-Bhaut}, in order to get an actual quotient and not just a homotopy quotient $\Emb_{e_0}(\mathrm{U}_1,\bR^k)\sslash\Emb^{\simeq}_{\partial}(\mathrm{U}_1,\mathrm{U}_1),$ it was important to have a simplicial group $S_{\bullet} \Diff_{\partial}(U_1)$ act freely on $S_{\bullet}\Emb_{e_0}(U_1,\bR^k)$. If we had not restricted ourselves to the interiors of thickenings, we would have ended up with modding out by $\Emb_{\partial}^{\simeq}(N_1,N_1)$ where $N_1$ is the closure of $\mathrm{U}_1$. However, it is well-known that $\Diff(N_1) \rightarrow \Emb_{\partial}^{\simeq}(N_1,N_1)$ is not an equivalence when $N_1$ is compact. Since $N_1$ and its interior $\mathrm{U}_1$ are isotopic, we could have replaced $\Emb_{\partial_0}^{\simeq}(N_1,N_1)$ with $\Emb_{\partial}(\mathrm{U}_1,\mathrm{U}_1)$ which is equivalent to $\Diff_{\partial}(\mathrm{U}_1)$. However the action of $\Diff_{\partial}(\mathrm{U}_1)$ on $\Emb_{e_0}(N_1,\bR^k)$ is not free.
\end{rk}

\subsection{The Poincaré Cobordism Category as a category internal to simplicial sets} \label{2.3}

In this subsection, we aim to define the Poincaré cobordism category over a space $X$ mentioned in the introduction $\text{Cob}_d^{\mathrm{G}}(\mathrm{X})$. Informally, it is an $\infty$-category with objects pairs $(P,f : P \rightarrow \mathrm{X})$ (up to self-equivalences) where $P$ is a $(d-1)$-dimensional Poincaré complex. Its morphism spaces are homotopy equivalent to \[ \text{Cob}^G_d((P_0,f_0),(P_1,f_1))= \bigsqcup_{W} \Map_{f_0,f_1}(W,\mathrm{X})\sslash\haut_{\partial}(W),\] where the disjoint union runs over Poincaré cobordisms $(W,P_0,P_1)$ from $P_0$ to $P_1$. Composition is given by gluing Poincaré cobordisms along the common boundary. The symmetric monoidal structure is given by disjoint union of objects. Its homotopy category $\text{hCob}_d^{SG}(\mathrm{X})$ is the category with objects pairs $(P,f : P \rightarrow \mathrm{X})$ up to equivalence and morphisms are Poincaré pairs $(W,P_0,P_1,F,f_0,f_1)$ up to relative equivalence. We define along the way an oriented Poincaré cobordism category $\CobdSG(\mathrm{X})$, where its objects and morphisms are oriented Poincaré duality spaces and cobordisms. To avoid disjunctions on oriented versus non-oriented cases, let $\PCobdG(\mathrm{X})$ denote either the unoriented category $\mathrm{Cob}_d^{\mathrm{G}}$ or the oriented one $\CobdSG$.

We give give a model $\PCobdG(\mathrm{X})$ as a non-unital simplicial category or equivalently a category internal to simplicial sets. To avoid any confusion, we clarify what we mean by simplicial category below.

\begin{definition}[Simplicial Categories]
    A non-unital simplicial category $\mathcal{C}$ is a category internal to simplicial sets $\text{sSet}$. It is equivalent to the data of a simplicial set of objects $\mathrm{Ob}(\mathcal{C})$, a simplicial set of morphisms $\mathrm{Mor}(\mathcal{C})$, maps of simplicial sets $s,t : \mathrm{Mor}(\mathcal{C}) \rightarrow \mathrm{Ob}(\mathcal{C})$ which send a morphism to respectively its source and its target, as well as a composition morphism $m : \mathrm{Mor}(\mathcal{C}) \times_{\mathrm{Ob}(\mathcal{C})} \mathrm{Mor}(\mathcal{C}) \rightarrow \mathrm{Mor}(\mathcal{C})$ satisfying strict associativity: $\mathrm{m(m(f,g),h)=m(f,m(g,h))}$.\\
    Its nerve $\mathrm{N}_{\bullet}\mathcal{C}$ is the semi-simplicial object in $\text{sSet}$ such that $\mathrm{N}_0 \mathcal{C}=\mathrm{Ob}(\mathcal{C}), \mathrm{N}_1\mathcal{C}=\mathrm{Mor}(\mathcal{C})$ and  $$\mathrm{N}_k \mathcal{C} =\mathrm{Mor}(\mathcal{C})\times_{\mathrm{Ob}(\mathcal{C})}...\times_{\mathrm{Ob}(\mathcal{C})} \mathrm{Mor}(\mathcal{C}).$$
\end{definition}

\begin{rk}
    Let $\mathcal{C}$ be a non-unital simplicial category. Let $|\mathrm{N}_{\bullet} \mathcal{C}|$ denote the semi-simplicial space obtained by taking levelwise the geometric realization of $\mathrm{N}_{\bullet} \mathcal{C}$. Then, $|\mathrm{N}_{\bullet} \mathcal{C}|$ is a semi-Segal space (as in \cite[Definition $4.2$]{Steinebrunner_2022}) if $|\mathrm{N}_1 \mathcal{C}| \rightarrow |\mathrm{N}_0\mathcal{C}| \times |\mathrm{N}_0 \mathcal{C}|$ is a Serre fibration. 
    
    In particular, if $(s,t) : \mathrm{N}_1 \mathcal{C} \rightarrow \mathrm{N}_0 \mathcal{C}$ is a Kan fibration, then $|\mathrm{N}_{\bullet} \mathcal{C}|$ is a Segal semi-simplicial space, hence $\mathcal{C}$ is a model of an $\infty$-category.
\end{rk}

We start with objects. We define a simplicial set $\psi^{\cSG}_{d,\bullet}(n)$ of thickenings of $d$-dimensional (oriented) Poincaré complex, analoguous to the space of submanifolds $\psi_d(n,0)$ from \cite{grw9}.

\begin{definition} \label{def-setsofpoinc}
    Let $\psi_{d,k}^{\cSG}(n)$ be the set of open subsets $U\subset \Delta^{\mathrm{k}} \times \bR^n$ such that there exists a level (and orientation)-preserving diffeomorphism 
\[\begin{tikzcd}
	U&& {\Delta^k \times \mathrm{int}(M)} \\
	& {\Delta^k}
	\arrow["\phi"{description}, from=1-1, to=1-3]
	\arrow["\pi"{description}, from=1-1, to=2-2]
	\arrow["\pi"{description}, from=1-3, to=2-2]
\end{tikzcd},\] where $M \subset \bR^n$ is a compact (oriented) $n$-manifold which has the homotopy type of a (oriented) Poincaré complex of dimension $d$.

        If $\mathrm{X}$ is a space, let $\psi_{d,k}(n,\mathrm{X})$ be the set of pairs $(U,f)$ where $U\in \psi_{d,k}^{\cSG}(n)$ and $f: U\rightarrow \Delta^k \times \mathrm{X}$ is a level-preserving map.
        
        The face maps and degeneracies of $\Delta^k$ make the collection $(\psi_{d,k}^{\cSG}(n,\mathrm{X}))_k$ into a simplicial set $\psi_{d,\bullet}^{\cSG}(n,\mathrm{X})$.
        
        Let $\psi^{\cSG}_{d,\bullet}(\mathrm{X})$ be the colimit of $\psi^{\cSG}_{d,\bullet}(n,\mathrm{X})$ under the identification maps $(U,f) \rightarrow (U\times \bR, f\times \text{id}_{\bR})$.
\end{definition}

\begin{rk}  
    In the definition above, the elements of the set of $0$-simplices of $\psi^{\cSG}_{d,\bullet}(n,\mathrm{X})$ are pairs $(U,f)$ where $U\subset \bR^n$ is \begin{itemize}
        \item open;
        \item diffeomorphic to the interior of a compact (oriented) submanifold of $\bR^n$;
        \item a Poincaré Duality space of dimension $d$;
    \end{itemize} and $f : U \rightarrow \mathrm{X}$ is a map.
\end{rk}

\begin{figure}[H]
\labellist 
\pinlabel {\textcolor{mygreen}{$\Delta^1$}} at 1700 400 
\pinlabel {\textcolor{mygreen}{$0$}} at 1870 150 
\pinlabel {\textcolor{mygreen}{$1$}} at 2780 70 
\pinlabel {\textcolor{mygreen}{$\sigma$}} at 2420 70 
\endlabellist
\includegraphics[height=3cm, width=8cm]{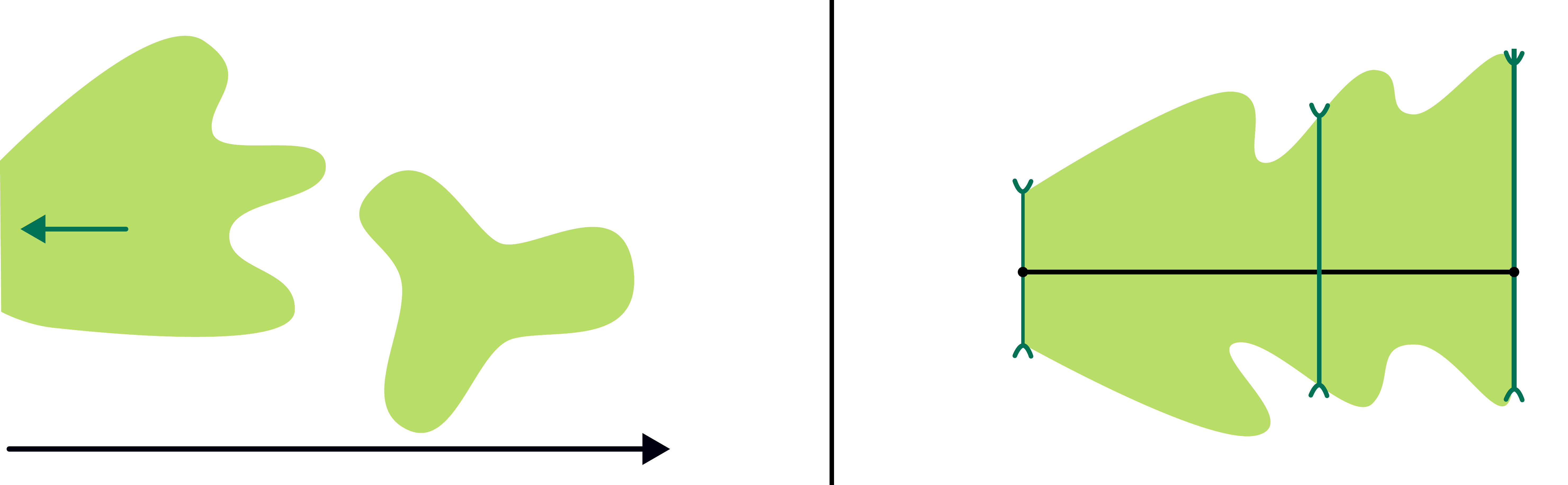}
\centering
\captionof{figure}{On the left: examples of $0$-simplices of $\psi_0^{\mathrm{SG}}(2)$. 
On the right: a $1$-simplex of $\psi_0^{\mathrm{SG}}(1)$.}
\label{fig:ex-simplices-objects-0}
\end{figure}
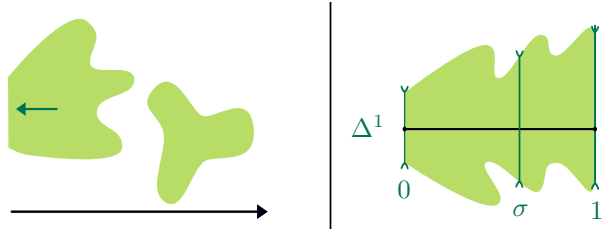

\begin{figure}[H]
\labellist 
\pinlabel {\Large\textcolor{mygreen}{$U_1$}} at 900 460 
\pinlabel {\Large\textcolor{mypurple}{$U_2$}} at 1300 1300 
\pinlabel {\Large\textcolor{myorange}{$U_3$}} at 2080 250
\endlabellist
\includegraphics[height=5cm, width=8cm]{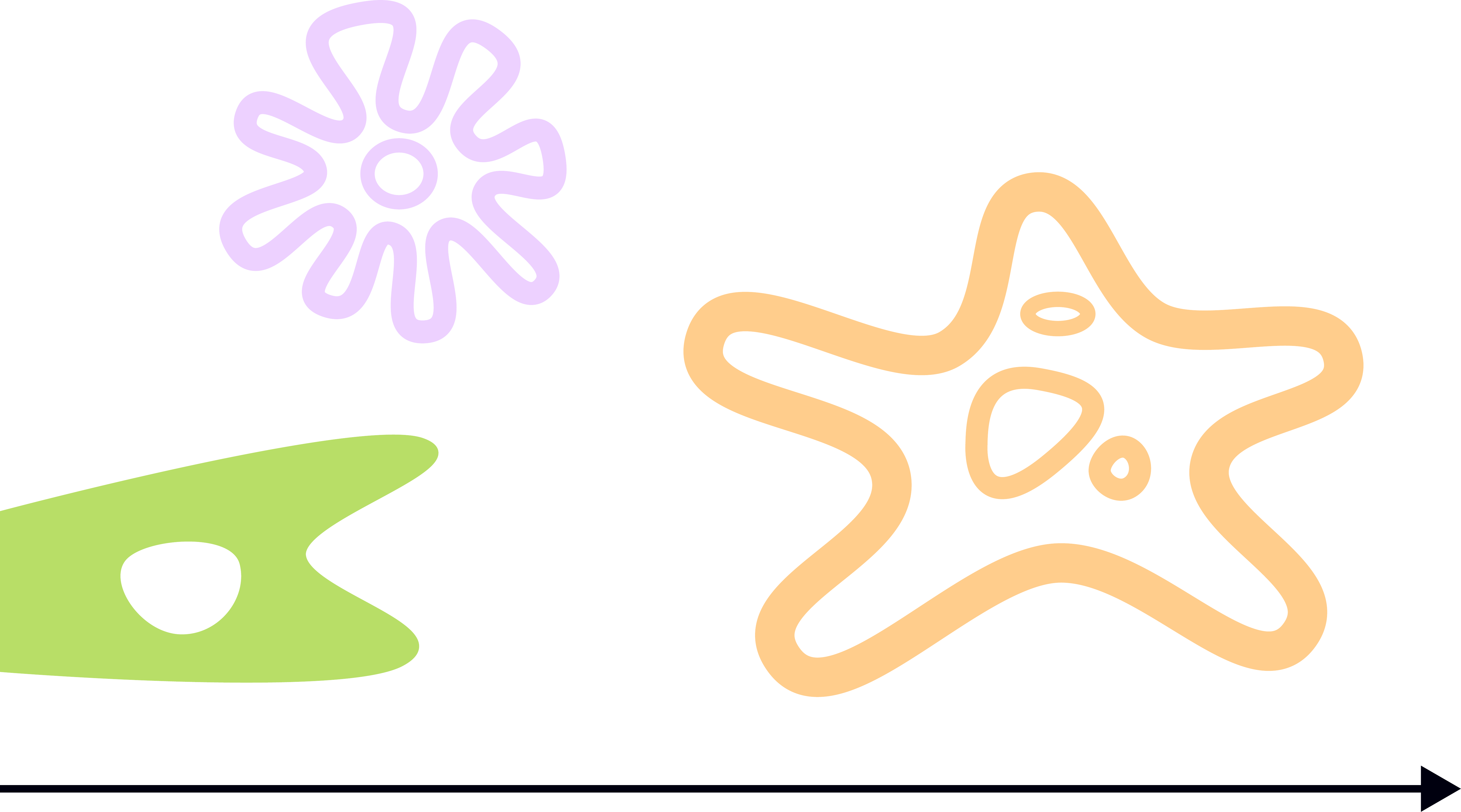}
\centering
\captionof{figure}{Three $0$-simplices $U_1,U_2$ and $U_3$ of $\psi_{1}^{\mathrm{SG}}(2)$}
\label{fig:ex-simplices-objects-1}
\end{figure}

Illustrations \eqref{fig:ex-simplices-objects-0} and \eqref{fig:ex-simplices-objects-1} give examples of simplices in $\psi_d^{\mathrm{SG}}(n)$ for low values of $d,n$.

In the lemma below, we describe the homotopy type of $|\psi_{d,\bullet}^{\cSG}(\mathrm{X})|$:

\begin{lem}
The space $|\psi_{d,\bullet}^{\cSG}(\mathrm{X})|$ is homotopy equivalent to $$\bigsqcup_{P} \Map(P,\mathrm{X})\sslash\haut^{(+)}(P)$$ where the disjoint union runs over (oriented) Poincaré duality spaces of dimension $d$.   
\end{lem}

\begin{proof}    
To begin with, we show that if $U_1$ and $U_2$ are two open $k$-dimensional thickenings of the same Poincaré duality space $P$, then for $l$ big enough, $U_1\times \bR^l$ and $U_2\times \bR^l$ are diffeomorphic. Let $\mathrm{P}$ be a Poincaré Duality space of dimension $\mathrm{d}$ and let $\mathrm{N}_1, \mathrm{N}_2 \subset \bR^n$ be codimension $0$ manifolds with boundary homotopy equivalent to $\mathrm{P}$. Let $U_{\mathrm{i}}$ be the interior of $\mathrm{N}_{\mathrm{i}}$. For $l>0$, $U_{\mathrm{i}} \times \bR^{l}$ is an open tubular neighborhood of $\mathrm{N}_{\mathrm{i}}$ in $\bR^{n+l}$. According to Corollary $2$ in \cite{Mazur1961StableEO}, for $l \geq k+2$, $U_1\times \bR^{l}$ and $U_2\times \bR^{l}$ are diffeomorphic. 

In particular no Poincaré Duality space is counted twice in $\psi_{\mathrm{d},\bullet}^{\cSG}(X)$. After stabilization, the set $\psi_{d,k}^{\cSG}(\mathrm{X})$ is exactly the set $$\bigsqcup_P \text{colim}_l \big(\mathrm{S}_k\Map(U_1\times \bR^l,\mathrm{X})\times \mathrm{S}_k\Emb(U_1\times \bR^l)\big)/\mathrm{S}_k\Diff(U_1\times \bR^l),$$ where the disjoint union runs over $P$ Poincaré complexes of dimension $d$ and an associated open thickening $U_1$. These identifications are compatible with the face maps and degeneracies. Since geometric realization commutes with all colimits, we can apply \cref{prop-approx-Bhaut} to conclude that $|\psi_{d,\bullet}^{\cSG}(\mathrm{X})|$ is equivalent to $\bigsqcup_P \Map(P,X)//\haut^{(+)}(P)$.
    
\end{proof}

\begin{notation}
    In what follows, let $\mathrm{x}_1 : \bR^n \rightarrow \bR$ denote the projection on the first coordinate. If $\mathrm{I} \subset \bR$ and $\mathrm{A}$ is a subset of $\bR^n$, let $$\mathrm{A}_\mathrm{I}$$ denote $$\mathrm{x}_1^{-1}(\mathrm{I})\cap \mathrm{A}.$$ Let $\mathrm{e}_1$ denote the norm $1$ vector defined in the $\mathrm{x}_1$-direction. If $\mathrm{A}\subset\bR^n$, we denote by $\mathrm{A}+\mathrm{t.e}_1$ the translation of $\mathrm{A}$ along the $\mathrm{x}_1$-axis.
\end{notation}

\begin{figure}[H]
\labellist 
\pinlabel {\Large{\textcolor{mygreen}{$\tau$}}} at 1890 35 
\pinlabel {\Large{\textcolor{mygreen}{$0$}}} at 380 35 
\endlabellist
\includegraphics[height=6cm, width=12cm]{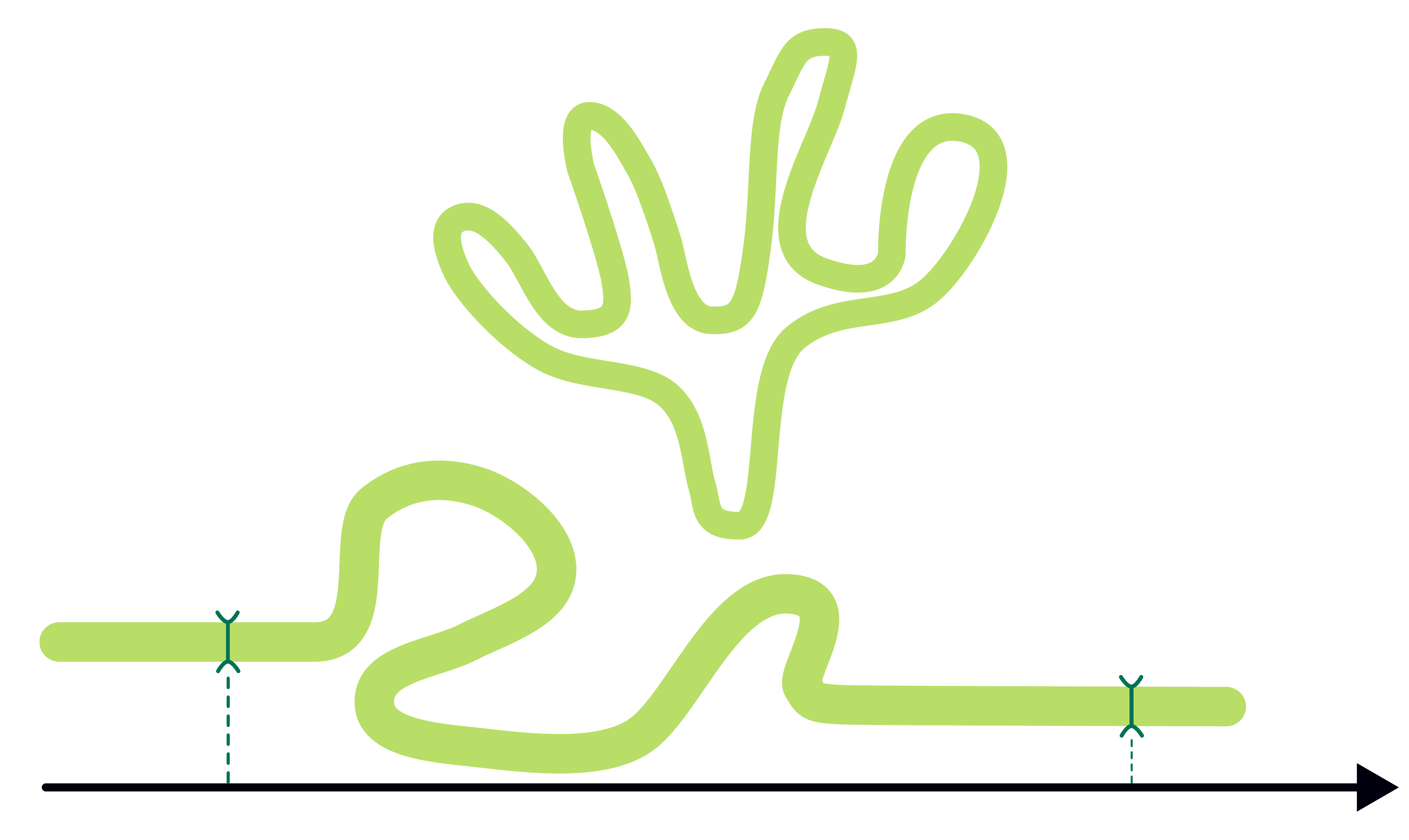}
\centering
\captionof{figure}{A $0$-simplex in $\psi_1^{\cSG}(2,1)$, or a $0$-simplex of the simplicial set of morphisms of $\mathrm{Cob}_1^{\cSG}(2)$}
\label{fig:morphisms-cob-cat}
\end{figure}

 Let $N$ be a compact manifold with boundary. If $\partial_0 N,\partial_1 N$ are two disjoint, compact, codimension $0$ submanifolds of $N$ (possibly with boundary), we denote by $\partial_2 N$ the complement $\partial N \setminus \mathrm{int}(\partial_0 N \sqcup \partial_1 N)$. We denote this data by the quadruple $(N,\partial N, \partial_0 N, \partial_1 N)$ .\\
 We give a definition of the simplicial set of $d$-dimensional Poincaré cobordisms $\psi_{d,\bullet}^{\cSG}(n,1)$, following the notation $\psi_d(n,1)$ from \cite{grw9}. See figure \ref{fig:morphisms-cob-cat} for examples of $0$-simplices in $\psi_{\mathrm{d}}^{\cSG}(n,1)$.

\begin{definition} \label{def2.24}
    Let $\psi_{d,k}^{\cSG}(n,1)$ be the set of pairs $(\mathrm{W},\tau)$ with $W \subset \Delta^k \times \bR^n$ and $\tau >0$ such that:
    \begin{enumerate}
        \item $W$ is open; 
        \item $W_{|(-\infty,0]} = W_0 \times (-\infty,0]$;
        \item $W_{[\tau,\infty)}=W_{\tau}\times [\tau,\infty)$.;
        \item there exists a (oriented) compact manifold quadruple $(N,\partial N, \partial_0 N, \partial_1 N)$ and a relative (orientation-preserving and) level-preserving diffeomorphism: \[\begin{tikzcd}
	{(\mathrm{W},\mathrm{W}_0, \mathrm{W}_t)} & {} & {(\Delta^k \times (N-\partial_2 N), \Delta^k\times \mathrm{int}(\partial_0 N), \Delta^k \times \mathrm{int}(\partial_1 N)) } \\
	& {\Delta^k}
	\arrow["\phi"{description}, from=1-1, to=1-3]
	\arrow["\pi"{description}, from=1-1, to=2-2]
	\arrow["\pi"{description}, from=1-3, to=2-2]
\end{tikzcd};\] 
        \item the triad $(N,\partial_0 N, \partial_1 N)$ is a homotopy equivalent to a $d$-dimensional Poincaré cobordism $\mathrm{(Q,P_0,P_1)}.$
    \end{enumerate}
If $\mathrm{X}$ is a space, let $\psi_{d,k}^{\cSG}(n,1,\mathrm{X})$ be the set of triples $(W,\tau,f)$ where $(W,\tau) \in \psi_{d,k}^{\cSG}(n)$ and $\mathrm{f}: \mathrm{W} \rightarrow \Delta^k \times \mathrm{X}$ is a level-preserving map.

        The face maps and degeneracies of $\Delta^k$ make the collection $(\psi_{d,k}^{\cSG}(n,1,\mathrm{X}))_k$ into a simplicial set $\psi_{d,\bullet}^{\cSG}(n,1,\mathrm{X})$.
        
        Let $\psi^{\cSG}_{d,\bullet}(1,\mathrm{X})$ be the colimit of $\psi^{\cSG}_{d,\bullet}(n,1,\mathrm{X})$ under the identification maps $(W,t,f) \rightarrow (W\times \bR, t, f\times \text{id}_{\bR})$.
\end{definition}

We can now define a simplicial category $\PCobdG(n,\mathrm{X})$ as follows: 

\begin{definition} \label{def-poinc-cat}
    Let $\PCobdG(n,\mathrm{X})$ be the simplicial category with:
    \begin{itemize}
        \item its simplicial set of objects $\mathrm{Ob}(\PCobdG(n,\mathrm{X}))$ is $\psi_{d-1,\bullet}^{\cSG}(n,\mathrm{X})$;
        \item its simplicial set of morphisms $\mathrm{Mor}(\PCobdG(n,\mathrm{X}))$ is $\psi_{d,\bullet}^{\cSG}(n,1,\mathrm{X})$;
        \item the source and target morphisms $s,t : \mathrm{Mor}(\PCobdG(n,\mathrm{X})) \rightarrow \mathrm{Ob}(\PCobdG(n,\mathrm{X}))$ are given by sending $\mathrm{(W,\tau,f)}$ to $\mathrm{(W_0,f_{|W_0})}$ and $\mathrm{(W,\tau,f)}$ to $\mathrm{(W_{\tau},f_{|W_{\tau}})}$ respectively; 
        \item the composition $(\mathrm{W,\tau,f})\circ (\mathrm{W',\tau',f'})$ is given by gluing along the common boundary:
        $$(\mathrm{W} \cup_{\mathrm{W}_{\mathrm{\tau}}} (\mathrm{W'}+ \mathrm{\tau}\mathrm{e}_1) , \mathrm{\tau+\tau'}, \mathrm{f}\cup_{\mathrm{f}_{\mathrm{W}_{\tau}}} \mathrm{f}').$$  
    \end{itemize}
\end{definition}

Crossing with $\bR$ defines functors  $$\PCobdG(n,\mathrm{X}) \rightarrow \mathrm{Cob}_{\mathrm{d}}^{\cSG}(\mathrm{n+1}, \mathrm{X}).$$
We can now define the $\mathrm{d}$-dimensional Poincaré Cobordism Category: 

\begin{definition}
    Let $\PCobdG(\mathrm{X})$ be the simplicial category \[\colim_{n\rightarrow \infty} \PCobdG(\mathrm{n,X}).\]
    Its simplicial set of objects is given by 
 \[\colim_{n\rightarrow \infty} \mathrm{Ob}(\PCobdG(n,\mathrm{X})).\]
 Its simplicial set of morphisms is \[\colim_{n\rightarrow \infty} \mathrm{Mor}(\PCobdG(n,\mathrm{X})).\]
 Moreover, postcomposing with a map $\mathrm{f} : \mathrm{X} \rightarrow \mathrm{Y}$ gives a functor $$\PCobdG(\mathrm{X}) \rightarrow \PCobdG(\mathrm{Y}).$$
\end{definition}

\begin{rk}
   The emptyset defines an object and a morphism in $\PCobdG(\mathrm{X})$.
\end{rk}

\begin{rk} \label{rk-non-unit}
    The category $\PCobdG$ is non-unital. In other words, it does not admit strict identity morphisms. However, it admits what are called weak units, or units up to equivalence. These are described in Subsection \ref{subsec33}. 
\end{rk}

According to \cref{prop-approx-Bhaut}, the above definition of the Poincaré cobordism category agrees with the hand-wavy one given in the introduction \ref{1.1}:

\begin{lem}
    Let $(U_0,f_0), (U_1,f_1)$ be two objects in $\PCobdG(\mathrm{X})$ such that $U_i$ is homotopy equivalent to a $(d-1)$-dimensional (oriented) Poincaré complex $\mathrm{P}_i$. Their simplicial set of morphisms $\PCobdG(\mathrm{X})((U_0,f_0),(U_1,f_1)) $ is equivalent after geometric realization to $$\bigsqcup_{\mathrm{Q}} \Map_{\mathrm{f}_0,\mathrm{f}_1}(\mathrm{Q,\mathrm{X}})\sslash\haut^{(+)}_{\partial}(\mathrm{Q,P_0,P_1)}),$$ where the disjoint union runs over $d$-dimensional (oriented) Poincaré cobordisms $(\mathrm{Q,P_0,P_1})$.
\end{lem}

We would like to write a functor from the smooth cobordism category $\PCobdO(X)$ to the Poincaré one $\PCobdG(X)$. However, for an embedded smooth closed manifold in $\bR^n$, there are many choices of thickenings possible. To do this, we replace  $\PCobdO(X)$ with an equivalent category $\mathrm{Cob}_d^{\cSO,\mathrm{tub}}(X)$, as done in \cite{Madsen-Tillmann}. In the end, we have a zigzag of functors: $$\PCobdO(X) \leftarrow \mathrm{Cob}_d^{\cSO,\mathrm{tub}}(X) \rightarrow \PCobdG(X).$$
Morally, the $0$-simplices of the objects of $\mathrm{Cob}_d^{\cSO,\mathrm{tub}}$ are pairs $(A,U)$, where $A$ is a closed subset of $\bR^n$ diffeomorphic to a closed manifold of dimension $(d-1)$ and $A\subset U\subset \bR^n$ is a tubular neighborhood of $A$. Before giving a clear definition of $\mathrm{Cob}_d^{\cSO,\mathrm{tub}}$, we discuss $\epsilon$-neighborhoods of compact manifolds: 

\begin{definition} \label{def-epsilon}
    Let $A\subset \bR^n$ be a smooth compact $(d-1)$-dimensional submanifold of $\bR^n$. The total space of its normal vector bundle $\nu_{A}$ is given by $$\nu_A:=\{(x,v)\in A\times \bR^n\ |\  v\in T_xA^{\perp}\}.$$ For $\epsilon>0$, we define $$D_{\epsilon}(\nu_A) :=\{(x,v)\in \nu_A\ |\  |v-x|<\epsilon\}.$$
    We define $$e(A):=\{\epsilon >0\ |D_{\epsilon}(\nu_A) \subset \bR^n\ \mathrm{is\  a\  tubular\ neighborhood\ of\ A}\}.$$
    
    If $f: A\rightarrow X$ is a map and $\epsilon \in e(A)$, we can define a map $f_{\epsilon} : D_{\epsilon}(\nu_A) \rightarrow X$ by letting $$f(x,v)=f(x).$$
\end{definition}

We now give a definition of the simplicial sets of smooth $(d-1)$-dimensional submanifolds of $\bR^n$:

\begin{definition}
    For $d,n\geq 0$, let $\psi_{d,k}^{\cSO}(n)$ be the set of closed subsets $A\subset \Delta^k \times \bR^n$ such that there exists a level and (orientation)-preserving diffeomorphism \[\begin{tikzcd}
	A && {\Delta^k \times M} \\
	& {\Delta^k}
	\arrow["\phi"{description}, from=1-1, to=1-3]
	\arrow["\pi"{description}, from=1-1, to=2-2]
	\arrow["\pi"{description}, from=1-3, to=2-2]
\end{tikzcd},\] where $M$ is compact closed (oriented) $d$-dimensional manifold.

Let $\psi_{d,k}^{\cSO, \mathrm{tub}}(n)$ be the set of pairs $(A,\epsilon)$ such that $A \in \psi_{d,k}^{\cSO}(n)$ and $\epsilon \in e(A)$.

For $X$ a space, let  $\psi_{d,k}^{\cSO}(n,X)$ be the set of pairs $(A,f)$ where $A \in \psi_{d,k}^{\cSO}(n)$ and $f \in \Map(A,X)$.

Let $\psi_{d,k}^{\cSO, \mathrm{tub}}(n,X)$ be the set of tuples $(A,\epsilon,f,f_{\epsilon})$ where $(A,\epsilon) \in \psi_{d,k}^{\cSO, \mathrm{tub}}(n)$, $f \in \Map(A,X)$ and $f_{\epsilon}$ is as in Definition \ref{def-epsilon}.

As in Definition \ref{def-setsofpoinc}, these sets form simplicial sets $\psi_{d,\bullet}^{\cSO}(n,X), \psi_{d,\bullet}^{\cSO, \mathrm{tub}}(n,X)$. Letting $n$ go to infinity, we get simplicial sets $\psi_{d,\bullet}^{\cSO}(X), \psi_{d,\bullet}^{\cSO,\mathrm{tub}}(X)$.

We can also define simplicial sets $\psi_{d,\bullet}^{\cSO}(n,1,X),\psi_{d,\bullet}^{\cSO,\mathrm{tub}}(n,1,X)$ of $d$-dimensional cobordisms with cylindrical ends as in Definition \ref{def2.24}.
\end{definition}

We can finally give a simplicial category model of the usual smooth cobordism category $\PCobdO(X)$ and its variation $\mathrm{Cob}_d^{\cSO,\mathrm{tub}}(X)$ as follows: 

\begin{definition}
    Let $\PCobdO(X)$ be the simplicial category with: 
    \begin{itemize}
        \item its simplicial set of objects is $\psi_{d-1}^{\cSO}(X)$;
        \item its simplicial set of morphisms is $\psi_{d}^{\cSO}(1,X)$;
        \item the source and target morphisms $s,t$ takes $(A,f,\tau) \in \psi_d^{\cSO}(1,X)$ to respectively $(A_0,f_{|A_0})$ and $(A_{\tau},f_{A_{\tau}})$;
        \item composition is given by union along the common boundary as in Definition \ref{def-poinc-cat}.
    \end{itemize}
    Let $\mathrm{Cob}_d^{\cSO,\mathrm{tub}}(X)$ be the simplicial category with:
    \begin{itemize}
        \item its simplicial set of objects is $\psi_{d-1}^{\cSO,\mathrm{tub}}(X)$;
        \item its simplicial set of morphisms is $\psi_{d}^{\cSO}(1,X)$;
        \item the source and target morphisms, as well as composition are as in Definition \ref{def-poinc-cat}.
    \end{itemize}
\end{definition}

There is a forgetful functor  $$\phi_X : \mathrm{Cob}_d^{\cSO,\mathrm{tub}}(X) \rightarrow \PCobdO(X),$$ 
which sends an object $(A,\epsilon,f,f_{\epsilon})$ of $\mathrm{Cob}_d^{\cSO,\mathrm{tub}}(X)$ to $(A,f)$. This construction is natural in $X$. Over each submanifold $A \subset \bR^n$, the space $e(A)$ of admissible $\epsilon$ is an interval, hence contractible. We deduce the following lemma, whose proof can be found in \cite{Madsen-Tillmann}:
\begin{lem}
    The forgetful functor $$ \phi_X :\mathrm{Cob}_d^{\cSO,\mathrm{tub}}(X) \rightarrow \PCobdO(X)$$ induces an equivalence  $$|\phi_X| :\B\mathrm{Cob}_d^{\cSO,\mathrm{tub}}(X) \rightarrow \B\PCobdO(X).$$
\end{lem}

\begin{figure}[H]
\labellist  
\pinlabel {\Large{\textcolor{mygreen}{$(A,\epsilon)$}}} at 3200 270 
\pinlabel {\Large{\textcolor{mygreen}{$A$}}} at 1350 270 
\pinlabel {\Large{\textcolor{mygreen}{$D_{\epsilon}(\nu_A)$}}} at 5120 270 
\endlabellist
\includegraphics[height=5cm, width=14cm]{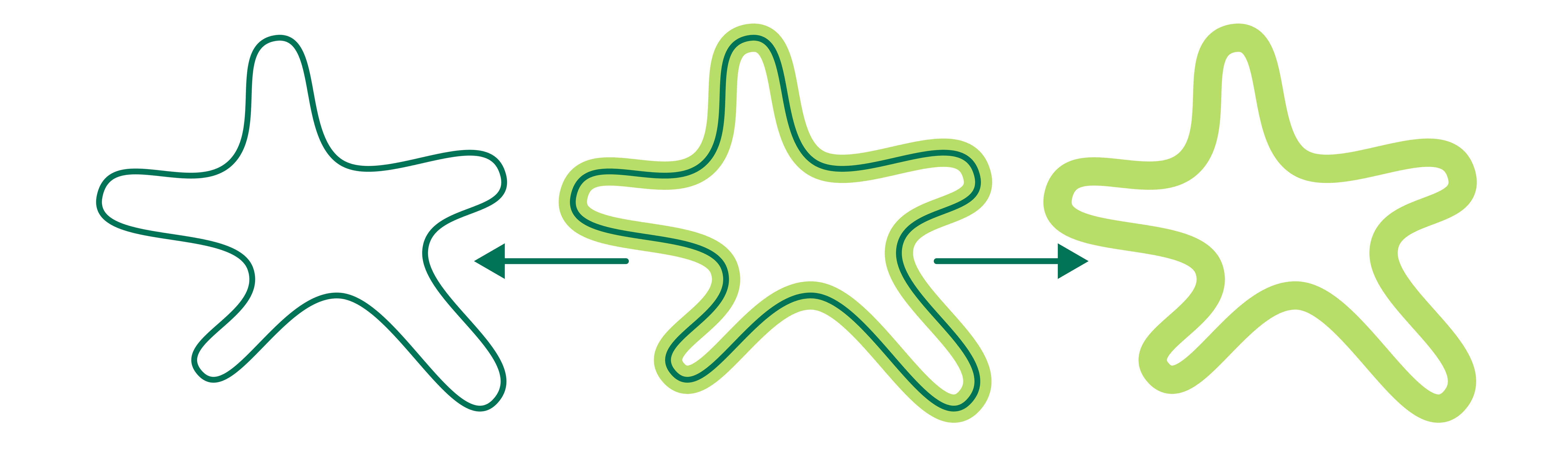}
\centering
\captionof{figure}{Action of the zigzag of functors $\PCobdO \leftarrow \mathrm{Cob}_d^{\mathrm{(S)O,tub}} \rightarrow \PCobdG$ on $0$-simplices of objects}
\label{fig:From-SO-to-SG}
\end{figure}

On the other hand, there is a forgetful functor $$\mathrm{tub}_X : \mathrm{Cob}_d^{\cSO,\mathrm{tub}}(X) \rightarrow \PCobdG(X) $$ sending a triple $(A,\epsilon, f, f_{\epsilon})$ to $(D_{\epsilon}(\nu_A),f_{\epsilon})$. This construction is again natural in $X$. The action of the resulting zigzag $$\PCobdO(X) \leftarrow \mathrm{Cob}_d^{\cSO,\mathrm{tub}}(X) \rightarrow \PCobdG(X)$$ on objects is illustrated on Figure \ref{fig:From-SO-to-SG}.

Taking the geometric realization of the nerve of the different cobordism categories gives functors $$\B\PCobdO(-), \B \PCobdG(-) : \mathcal{S} \rightarrow \mathcal{S}.$$ Here $\mathcal{S}$ denotes the $\infty$-category of spaces. By the previous lemma, the functors $\B\PCobdO(-)$ and $\B\mathrm{Cob}_d^{\cSO,\mathrm{tub}}(-)$ are equivalent, thus we can construct a functor as in the following lemma:

\begin{lem} \label{lemme-smoothtopoinc}
 The functor $\mathrm{tub}_X : \mathrm{Cob}_d^{\cSO,\mathrm{tub}}(X) \rightarrow \PCobdG(X) $ induces a natural transformation
 $$u(-) : \B \PCobdO(-) \rightarrow \B \PCobdG(-).$$
\end{lem}

We end on a remark on connected components of cobordism categories.

\begin{rk} \label{rk-pi_0-cob-cat}
    The group of connected components $\pi_0 (\B \PCobdO(X))$ is the $(d-1)$-dimensional cobordism group $\Omega_{d-1}^{\cSO}(X)$ over $X$. The group of connected components $\pi_0(\B\PCobdG(X))$ is the $(d-1)$-dimensional cobordism group $\Omega_{d-1}^{\cSG}(X)$.
\end{rk}

\begin{rk}
    In \cite{grw9}, the authors study moduli spaces of manifolds equipped with a tangential structure. For Poincaré complexes, we could make sense of a "Spivak structure". Let $B$ be a connected space and $\sigma : B \rightarrow \mathrm{Pic}(\bS)$ be a map. We say that a Poincaré complex $P$ admits a $\sigma$-structure if the classifying map $P \rightarrow \mathrm{Pic}(\bS)$ of its Spivak fibration $\nu_P$ admits a lift $s$ to $\sigma$ as in the following diagram:
\[\begin{tikzcd}
	& B \\
	P & {\mathrm{Pic}(\mathbb{S})}
	\arrow["\sigma", from=1-2, to=2-2]
	\arrow["s", dashed, from=2-1, to=1-2]
	\arrow["{\nu_P}"', from=2-1, to=2-2]
\end{tikzcd}.\] We say that the couple $(P,s)$ is Poincaré complex equipped with a Spivak $\sigma$-structure. We could now consider the self-homotopy equivalences $\haut^{\sigma}(P)$ preserving the $\sigma$-structure. One could wonder what could the analoguous simplicial model of $\B\haut^{\sigma}(P)$ be. Similarly, we could define a cobordism category of Poincaré complexes equipped with a $\sigma$-structure $\mathrm{Cob}_d^{\mathrm{SG},\sigma}$.
\end{rk}

\subsection{Delooping the functor \texorpdfstring{$\B\PCobdG(-)$}{B  PCobdG(-)}} \label{2.4}

In this subsection, we aim to show that the functor $\B \PCobdG(-)$ factors through spectra. More precisely, we construct a connective delooping of $\B \PCobdG(-)$ compatible with the connective delooping of $\B\PCobdO(-)$ from \cite{GMTW}. Let $\mathrm{Sp}^{\geq 0}$ denote the category of connective spectra. Let $\tau_{\geq 0} : \mathrm{Sp} \rightarrow \mathrm{Sp}^{\geq 0}$ denote the connective truncation functor. 

The main Theorem of \cite{GMTW} says there is a natural equivalence $$\B \PCobdO(X) \simeq \lop (\Sigma \mathrm{MT(S)O}(d) \otimes \susplus X).$$ We recall that the (oriented) Madsen-Tillmann spectrum $\mathrm{MT(S)O}(d)$ is the Thom spectrum of the stable inverse of the universal (oriented) $d$-bundle over $\mathrm{B(S)O}(d)$. In particular, $$\tau_{\geq 0}(\Sigma \mathrm{MT(S)O}(d) \otimes \susplus -)$$ defines a connective delooping of the functor $\B \PCobdO(-)$. We can now state the main result of this subsection:

\begin{prop} \label{main-prop-4.2}
    There exists a functor $$\mathrm{C}^{\mathrm{(S)G}}(d,-) :\mathcal{S} \rightarrow \mathrm{Sp}^{\geq 0}$$ and a natural equivalence $$\B \PCobdG(-) \Rightarrow \lop \mathrm{C}^{\mathrm{(S)G}}(d,-).$$    
Moreover, there exists a natural transformation $$\tilde{u} : \tau_{\geq 0}(\Sigma \mathrm{MT(S)O}(d) \otimes (\susplus -))\Rightarrow \mathrm{C^{(S)G}}(d,-)$$ such that $\Omega^{\infty} \circ \tilde{u} $ is equivalent to the natural transformation $u(-) :  \B \PCobdO(-)\rightarrow  \PCobdG(-)$ from Lemma \ref{lemme-smoothtopoinc}.
\end{prop}

The disjoint union of two Poincaré complexes of dimension $(d-1)$ is again a $(d-1)$-dimensional Poincaré complex. Then, we could show Proposition \ref{main-prop-4.2} by defining a symmetric monoidal structure on $\PCobdG(X)$, with the monoidal product being the disjoint union. However, in our model disjoint unions of objects may not be disjoint as subsets of the Euclidean space. Thus, defining a monoidal structure on $\PCobdG$ may be too rigid for our context. Instead, we define a $\Gamma$-space structure on $\B\PCobdG(-)$, as in \cite{Nguyen_2017}.

Let $\Gamma^{\mathrm{op}}$ be the Segal category with objects pointed sets $n_+: = \{\star,1,\ldots,n\}$ for $n\geq 0$ and the set of maps $\Gamma^{\mathrm{op}}(n,m)$ is the set of based maps of sets $\{\star, 1, \ldots, n\} \rightarrow \{\star, 1, \ldots, m\}$. Let $\rho_i : \{\star, 1, \ldots, n\} \rightarrow \{ \star, 1\}$ be the map sending every element but $i$ to $\star$. We recall the following definitions: 

\begin{definition}
    A \textcolor{black}{$\Gamma$-space} is a functor $X : \Gamma^{\mathrm{op}} \rightarrow \mathcal{S}$. It is a \emph{special} $\Gamma$-space if $X(0_+)$ is contractible and if it satisfies the \emph{Segal condition}, i.e. the Segal maps
\[\begin{tikzcd}
	{X(n_+)} && {X(1_+) \times \ldots X(1_+)}
	\arrow["{(\rho_1,\ldots, \rho_n)}", from=1-1, to=1-3]
\end{tikzcd}\] are equivalences for any $n$.

A map of very special $\Gamma$-spaces $f : X \rightarrow Y$ is a natural transformation from $X$  to $Y$.

Let $\Gamma^s\mathcal{S}$ denote the category of special $\Gamma$-spaces. Let $\Gamma^{vs}\mathcal{S}$ denote the category of very special $\Gamma$-spaceS.
\end{definition}

Let $X : \Gamma^{\mathrm{op}} \rightarrow \mathcal{S}$ be a special $\Gamma$-space. We have the following zigzag of maps: 
\[\begin{tikzcd}
	{X(1_+)\times X(1_+)} & {X(2_+)} & {X(1_+)}
	\arrow["\cong"', from=1-2, to=1-1]
	\arrow["\mu", from=1-2, to=1-3]
\end{tikzcd}.\] Here, the equivalence $X(2_+) \rightarrow X(1_+) \times X(1_+)$ is the Segal map, while $\mu$ is induced by the constant map $\{1,2\} \rightarrow \{1\}$. Picking a homotopy inverse of the Segal map gives a multiplication $X(1_+) \times X(1_+) \rightarrow X(1_+)$, making $X(1_+)$ into a $H$-space.

At the level of connected components, we get a multiplication $\pi_0(X(1_+)) \times \pi_0(X(1_+)) \rightarrow \pi_0(X(1_+))$. 

\begin{definition}
    A special $\Gamma$-space $X$ is \emph{grouplike} or \emph{very special} if the monoid $\pi_0(X(1))$ is actually a group. 
\end{definition}

The upshot of \cite{SEGAL1974293} is that the $H$-space $X(1_+)$ is a group-like $E_{\infty}$-space if $X : \Gamma^{\mathrm{op}} \rightarrow \mathcal{S}$ is a very special $\Gamma$-space. In \cite{SEGAL1974293}, Segal constructs a functor $$\mathbf{B}(-) : \Gamma^{s} \mathcal{S} \rightarrow \mathrm{Sp}^{\geq 0}.$$ 

Let $X$ be a very special $\Gamma$-space. The following theorem states that $\mathbf{B}(X)$ is a connective delooping of $X(1_+)$:

\begin{thm} \label{Segalstheorem}
    There is a natural equivalence from the functor \[
\begin{array}{rcl}
 \Gamma^s \mathcal{S} & \longrightarrow & \mathcal{S} \\
X & \longmapsto & X(1_+)
\end{array}
\] to the composite 
\[\begin{tikzcd}
	{\Gamma^{vs}\mathcal{S}} & {\mathrm{Sp}^{\geq 0}} & {\mathcal{S}}
	\arrow["{\mathbf{B}(-)}"', from=1-1, to=1-2]
	\arrow["\lop"', from=1-2, to=1-3]
\end{tikzcd}.\] 
\end{thm}

We now define the following $\Gamma$-space structure on $\B \PCobdG(X)$:

\begin{definition}
    Let $\B\PCobdG(X)\langle - \rangle : \Gamma^{\mathrm{op}} \rightarrow \mathcal{S}$ be the $\Gamma$-space defined by letting $$\B \PCobdG(X)\langle n_+ \rangle  = \B \PCobdG(X \times \{1,\ldots, n\}).$$ 
    A based map $\lambda : n_+ \rightarrow m_+$ induces a map $\lambda_X : X\times n_+ \rightarrow X\times m_+$. Let $$\PCobdG(X)\langle \lambda\rangle : \PCobdG(X)\langle n_+ \rangle  \rightarrow \PCobdG(X)\langle m_+ \rangle $$ be the functor sending an object $(U,f) \in \mathrm{Ob}(\PCobdG(X)\langle n_+ \rangle)$ to the object $(\lambda_X^{-1}(U),f_{|\lambda_X^{-1}(U)})$ of $\PCobdG(X)\langle m_+ \rangle$. A morphism $(W,\tau,F)\in \mathrm{Mor}(\PCobdG(X)\langle n_+ \rangle)$ is sent to the morphism $(\lambda_X^{-1}(W), \tau, F_{|\lambda_X^{-1}(W)})$. In other words, the functor $\PCobdG(X)\langle \lambda \rangle$ deletes the connected components of the morphisms and objects which are mapped to $X\times \{\star\}$ through $\lambda_X$.
\end{definition}

In the following proposition, we show $\B \PCobdG(X)\langle - \rangle$ indeed defines a very special $\Gamma$-space. 

\begin{prop} \label{veryspecialgamma}
    For $X$ a space, $\B \PCobdG(X)\langle - \rangle$ is a very special $\Gamma$-space.
\end{prop}

Along the way, we prove a lemma on the limits and colimits that the functor $\B\PCobdG(-)$ preserves.

\begin{lem} \label{coproductsproducts}
    The functor $\B \CobdG(-)$ preserves filtered colimit and sends finite coproducts to finite products. 
\end{lem}

\begin{proof}
    Poincaré complexes and Poincaré pairs are compact in $\mathcal{S}$, hence both $\text{Map}(P,-)$ and $\text{Map}_{f_0,f_1}(Q,-)$ commute with filtered colimits. Since finite products and colimits commute with filtered colimits, we deduce that $\text{Ob}(\PCobdG(-))$ and $\text{Mor}(\PCobdG(-))$ commute with filtered colimits. Finite products and geometric realization commute with filtered colimits as well. We deduce $\B\PCobdG$ commutes with filtered colimits.
    
    Secondly, we consider the functor $\PCobdG(X\sqcup Y) \rightarrow \PCobdG(X)$ taking an object $(U,f)$ in $\mathrm{Ob}(\PCobdG(X\sqcup Y))$ to $(\iota_X^{-1}(U), f_{|\iota_X^{-1}(U)}),$ where $\iota_X : X \rightarrow X\sqcup Y$ is the standard inclusion. Similarly, it deletes the connected components of the morphisms which map to $Y$. Using decompositions of the mapping spaces $\Map(P_1\sqcup ...\sqcup P_n,X\sqcup Y)$ as a disjoint union of products, for $P_i$ connected spaces, we can conclude the map of semi-simplicial spaces $$|N_{\bullet}\PCobdG(X\sqcup Y))| \rightarrow |N_{\bullet}\PCobdG(X)| \times |N_{\bullet}\PCobdG(Y)|$$ is a levelwise equivalence. Finally, geometric realization of semi-simplicial spaces commutes with finite products according to \cite[Theorem $7.2$]{Ebert_2019}, which concludes the proof.
\end{proof}

\begin{proof}[Proof of Proposition \ref{veryspecialgamma}]
    Firstly, $\PCobdG(X)\langle0_+ \rangle$ is the category with one object $\emptyset$ and one morphism, hence $\B \PCobdG(X)\langle0_+ \rangle$ is contractible.
    
    We then need to show  $\B \PCobdG(X)\langle - \rangle$ satisfies the Segal condition. The space $X\times \{1,\ldots,n\}$ is a finite disjoint union, hence the Segal map is an equivalence according to Lemma \ref{coproductsproducts}. 

    Secondly, we need to show the multiplication on $\pi_0(\B\PCobdG(X))$ defined by the $\Gamma$-structure is exactly the disjoint union. The objects of  $\PCobdG(X)\langle 2_+ \rangle$ are objects $(U,f)$ of $\PCobdG(X)$ such that the connected components of $U$ are labelled by either $1$ and $2$. The category $$\PCobdG(X)\langle 1_+\rangle$$ is simply the category $\PCobdG(X)$. The map $$\mu : \pi_0(\B\PCobdG(X)\langle 2_+ \rangle ) \rightarrow \pi_0(\B\PCobdG(X))$$ forgets the labeling on objects. On the other hand, for $i=1$ or $2$, the functor $\PCobdG(X)\langle \rho_i  \rangle : \PCobdG(X)\langle 2_+\rangle \rightarrow \PCobdG(X)\langle 1_+\rangle $ sends an object ($U,f)$  to $(U^i, f^i)$, where $U^i$ is obtained from $U$ by restraining to components of $U$ labeled by $i$. Then, the Segal map
    $$\pi_0(\B \PCobdG(X)\langle 2_+\rangle) \rightarrow \pi_0(\B\PCobdG(X)\langle 1_+\rangle) \times \pi_0(\B\PCobdG(X)\langle 1_+\rangle)$$ decomposes objects and morphisms into their components labelled by $1$ or $2$. Let $[(U^1,f^1)]$ and $[(U^2,f^2)]$ be connected components of $\B \PCobdG(X)$, such that the representants $U^1$ and $U^2$ are disjoint. An inverse of the Segal map is then given by sending $[(U^1,f^1)]$ and $[(U^2,f^2)]$ to the connected component $[(U^1\sqcup U^2, F)]$. Here, the map $F : U^1 \sqcup U^2 \rightarrow X \times \{1,2\}$ labels each component $U^ i$ by $i$. The multiplication on $\pi_0(\B\PCobdG(X))$ induced by the $\Gamma$-structure coincides then with the disjoint union. As observed in \cref{rk-pi_0-cob-cat}, $\pi_0(\B\PCobdG(X))$ is in bijection with the Poincaré cobordism group $\Omega_d^{\cSG}(X)$. Since the latter, equipped with the disjoint union, is a group, we deduce that  the $\Gamma$-space $\B\PCobdG(X)\langle - \rangle$ is very special. 
\end{proof}

We now conclude with the proof of Proposition \ref{main-prop-4.2}:

\begin{proof}[Proof of \cref{main-prop-4.2}]
According to Proposition \ref{veryspecialgamma}, for any space $X$, $\B \PCobdG(X) \langle - \rangle$ is a very special $\Gamma$-space. Let $C^{\cSG}(d,X)$ denote the connective spectrum $\textbf{B} (\B \PCobdG(X) \langle - \rangle)$. Moreover, any map $f : X \rightarrow Y$ induces a functor $\PCobdG(X) \rightarrow \PCobdG(Y)$, hence induces a map of $\Gamma$-spaces  $$f : \B \PCobdG(X) \langle - \rangle \rightarrow \B \PCobdG(Y) \langle - \rangle. $$ According to Theorem \ref{Segalstheorem}, the map $f$ induces a map of spectra $C^{\cSG}(d,X) \rightarrow C^{\cSG}(d,Y)$. Consequently, $C^{\cSG}(d,-)$ defines a functor from $\mathcal{S}$ to $\mathrm{Sp}^{\geq 0}$. According to Theorem \ref{Segalstheorem}, the functor $C^{\cSG}(d,-)$ deloops $\B \PCobdG(-)$.

For the second part, the $\Gamma$-space $\B\PCobdG(X)\langle - \rangle$ is essentially identical to the $\Gamma$-space $\B\PCobdO(X) \langle - \rangle$ defined in \cite[Definition $8$]{Nguyen_2017}. In particular, $\B\PCobdO(X)\langle n \rangle$ is equivalent to $\B \PCobdO(X\times \{1,\ldots, n\}).$ The natural transformation $u(-) :\B \PCobdO(-) \Rightarrow \B \PCobdG(-)$ then induces a map of $\Gamma$-spaces $B\PCobdO(X)\langle - \rangle  \rightarrow B\PCobdG(X)\langle - \rangle $.

Finally, taking for each $X$ the connective delooping $\textbf{B} (\B \PCobdO(X) \langle - \rangle)$ defines a functor $C^{\cSO}(d,-) : \mathcal{S} \rightarrow \mathrm{Sp}^{\geq 0}$ which deloops $\B\PCobdO(-)$. According to the main theorem of \cite{Nguyen_2017}, there are natural equivalences $$C^{\cSO}(d,X) \simeq \tau_{\geq 0} (\Sigma \mathrm{MT(S)O}(d) \otimes \susplus X).$$
\end{proof}

\section{A formula for \texorpdfstring{$\B \sCobtwo(X)$}{B sCobtwo(X)}}\label{3}

As announced in the introduction, we restrict ourselves to the two dimensional case. In subsection \ref{2.1}, we saw that every $2$-dimensional (resp $1$-dimensional) oriented Poincaré pair is homotopy equivalent to a $2$-dimensional manifold (resp $1$-dimensional). In the previous subsection, we constructed a Poincaré cobordism category $\sCobtwo(X)$ and a comparison map from the geometric realization of the nerves of the smooth cobordism category $\B \sCobtwooo(X) \rightarrow \B \sCobtwo(X).$ 

Let $\mathrm{H}(M)$ denote either the monoid of orientation-preserving diffeomorphisms or 
self-homotopy equivalences of a manifold $M$. 
Let $\sCobtwoH(X)(\emptyset,\emptyset)$ denote the space of endomorphisms of $\emptyset$ in the category $\sCobtwoH(X)$.  We can restrict to the submonoid \textcolor{black}{$\mathrm{Sph^H}(X)$} of endomorphisms of $\emptyset$ in $\sCobtwoH(X)$ which are homotopy equivalent to finite disjoint unions of spheres. There is an equivalence of monoids: $$\mathrm{Sph^H}(X) \simeq \bigsqcup_{n\geq 0} \Map(S^2\times \{1,...,n\},X)\sslash \mathrm{H}(S^2\times \{1,...,n\}),$$ where the multiplication is given by disjoint union. 

Since $\mathrm{Sph^H}(X)$ is in particular a subcategory of $\sCobtwoH(X)$, the passage from smooth to Poincaré $u(X) : \B\sCobtwooo(X) \rightarrow \B \sCobtwo(X)$ restricts naturally to a map $\B \mathrm{Sph^{SO}}(X) \rightarrow \B \mathrm{Sph^{SG}}(X).$ After taking loops at the empty object, we get a commutative square:

\begin{equation} \label{main-square-thmB}
\begin{tikzcd}
	\Omega_{\emptyset} {\mathrm{BSph^{SO}}(X)} & \Omega_{\emptyset} {\mathrm{BCob}_2^{\mathrm{SO}}(X)} \\
	\Omega_{\emptyset} {\mathrm{BSph^{SG}}(X)} & \Omega_{\emptyset} {\mathrm{BCob}_2^{\mathrm{SG}}(X)}
	\arrow[from=1-1, to=1-2]
	\arrow[from=1-1, to=2-1]
	\arrow[from=1-2, to=2-2]
	\arrow[from=2-1, to=2-2]
\end{tikzcd}.
\end{equation}

We aim to show in the following subsections \ref{subsec32} and \ref{subsec33} that the square \eqref{main-square-thmB} is homotopy cartesian. To conclude the proof of Theorem \ref{thmB}, we show that the group-completion $ \Omega_{\emptyset}\B \mathrm{Sph^H}(X)$ of $\mathrm{Sph^H}(X)$ is equivalent to the free infinite loop space $\lop \susplus (\Map(S^2,X) \sslash \mathrm{H}(S^2))$, where $\mathrm{H}(S^2)$ denotes either orientation-preserving diffeomorphisms or self-equivalences of $S^2$.

\begin{figure}
\labellist 
\pinlabel {\large{\textcolor{mypurple}{$\mathrm{red}(\Sigma)$}}} at 200 600 
\pinlabel {\large{\textcolor{myorange}{$\mathrm{sph}(\Sigma)$}}} at 600 500 
\endlabellist
\includegraphics[width=10cm, height=7cm]{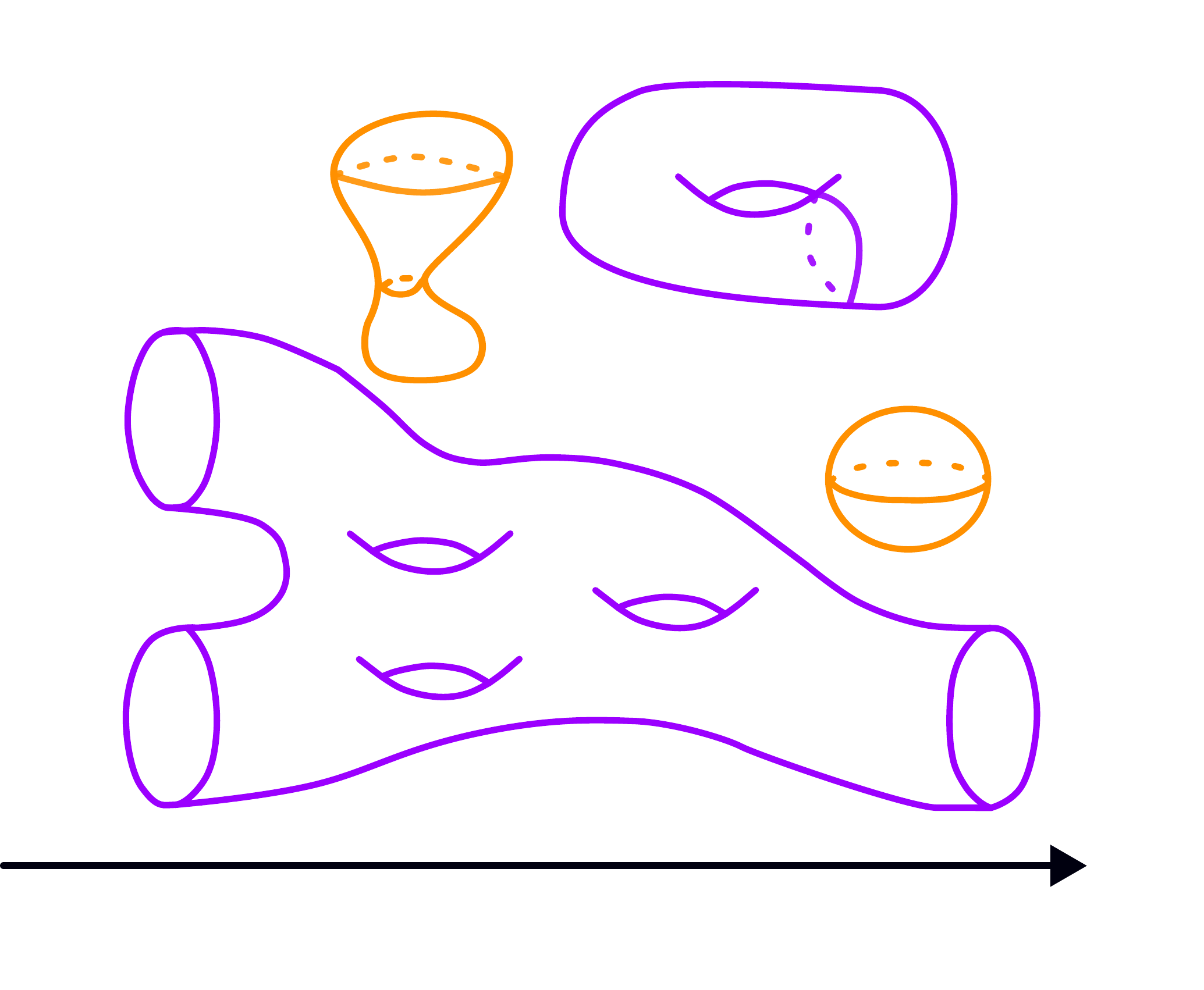}
\centering
\captionof{figure}{The reduced components $\textcolor{mypurple}{\mathrm{red}(\Sigma)}$ and spherical components $\textcolor{Orange}{\mathrm{sph}(\Sigma)}$ of a surface $\Sigma$}
\label{fig:closed-and-reduced}
\end{figure}

\begin{figure}[H]
\labellist  
\pinlabel {\large{\textcolor{mypurple}{$\Sigma_2$}}} at 1700 600 
\pinlabel {\large{\textcolor{mypurple}{$\Sigma_1$}}} at 500 600 
\pinlabel {\large{\textcolor{mypurple}{$M_1$}}} at 810 450 
\pinlabel {\large{\textcolor{mypurple}{$M_1$}}} at 1310 450 
\endlabellist
\includegraphics[height=5cm, width=14cm]{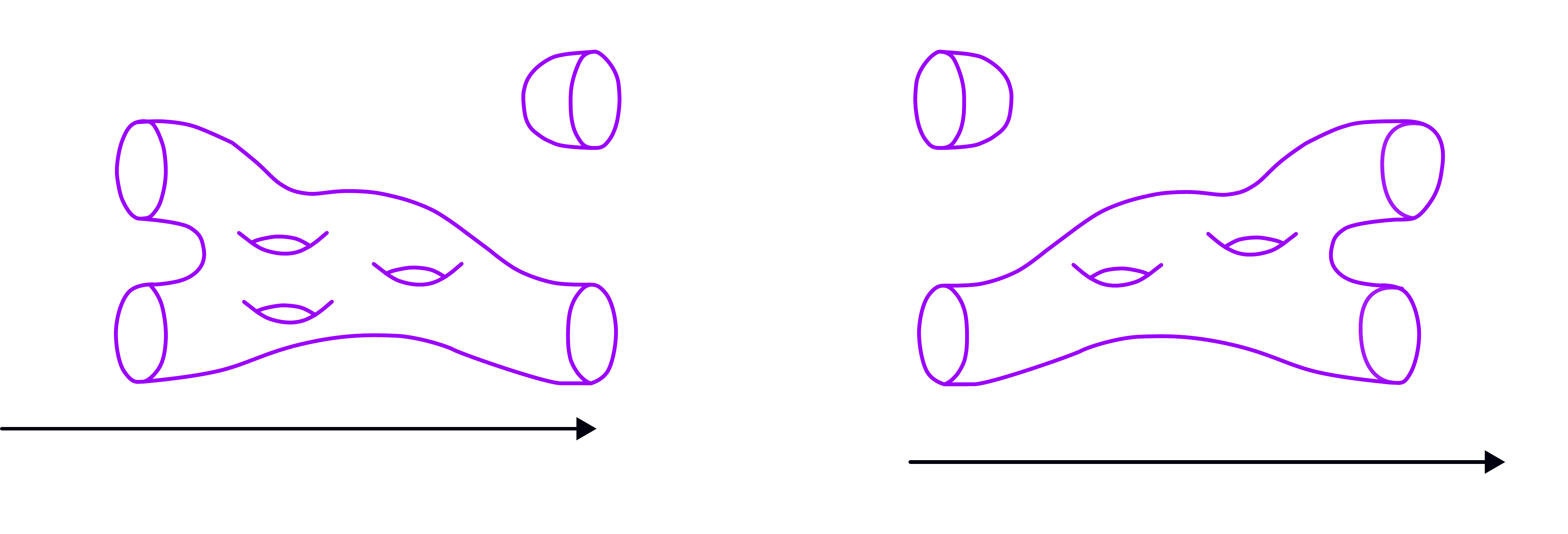}
\centering
\captionof{figure}{Two reduced surfaces $\Sigma_1, \Sigma_2$ such that $\Sigma_1 \cup_{M_1} \Sigma_2$ is not reduced}
\label{fig:two-red-compo}
\end{figure}
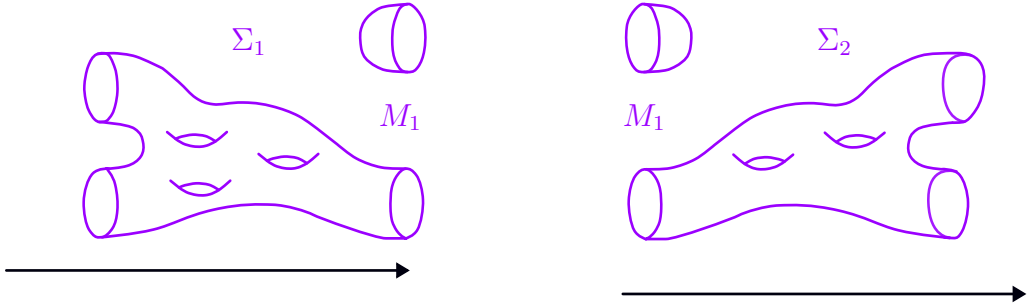

To begin with, we compare in subsection \ref{3.1bis} the homotopy types of the diffeomorphism group of surfaces and the monoid of self-homotopy equivalences of surfaces. In Theorem \ref{thm-hautvsdiff}, we see that the monoid map $\Diff_{\partial}(\Sigma) \rightarrow \haut_{\partial}^+(\Sigma)$ is a homotopy equivalence as long as $\Sigma$ has no connected component diffeomorphic to $S^2$. This suggests we define a \textcolor{black}{reduced} cobordism category \textcolor{black}{$\mathrm{Cob}_2^{\mathrm{H,red}}(X)$}, by restricting to morphisms having no connected component equivalent to $S^2$, see Definition \ref{def-3.7-red-cob} and Figures \ref{fig:two-red-compo} and \ref{fig:closed-and-reduced} in Subsection \ref{subsec32} for examples. Using the results from subsection \ref{3.1bis}, we show in \cref{same-red} that $\mathrm{Cob}_2^{\mathrm{SO,red}}(X)$ and $\mathrm{Cob}_2^{\mathrm{SG,red}}(X)$ have equivalent nerves. By deleting the spherical components of the morphisms in $\mathrm{Cob}_2^{\mathrm{H}}(X)$, we obtain a reduction functor:

\begin{equation} \label{red-functor}
\begin{tikzcd}
	{\mathrm{Cob}_2^{\mathrm{H}}(X) } & {\mathrm{Cob}_2^{\mathrm{H,red}}(X)}
	\arrow[from=1-1, to=1-2]
\end{tikzcd}.
\end{equation}
 Lastly, in Subsection \ref{subsec33}, we apply \textcolor{black}{Quillen's Theorem B} to the reduction functor \eqref{red-functor} to identify the homotopy fiber of the reduction map $\B \sCobtwoH(X) \rightarrow  \B \cobred(X)$ with $ \B \mathrm{Sph^H}(X)$, which concludes the proof of \cref{thmB}. 

\subsection{Diffeomorphisms versus self-homotopy equivalences of surfaces} \label{3.1bis}

In this subsection, we aim to compare the homotopy types of the diffeomorphism groups of surfaces $\Diff_{\partial}(\Sigma_{g,n})$ and their monoid of self-homotopy equivalences $\haut_{\partial}^+(\Sigma_{g,n})$. The main result is as follows: 

\begin{thm} \label{thm-hautvsdiff}
    Let $\Sigma_{g,n}$ be an surface such that $g+n>0$. Then the monoid map $$\Diff_{\partial}(\Sigma_{g,n}) \rightarrow \haut_{\partial}^+(\Sigma_{g,n})$$ is a homotopy equivalence.
\end{thm}

Since $\sCobtwoH$ has a simplicial set of objects, and not just a set, we compare in the following lemma the homotopy type of the spaces of objects of $\sCobtwooo$ and $\sCobtwo$:

\begin{lem} \label{lemS^1}
The monoid map $$\Diff^+(S^1) \rightarrow \haut^+(S^1)$$ is a homotopy equivalence.
\end{lem}

\begin{proof}
    The rotations induce an equivalence $S^1 \rightarrow \Diff^+(S^1)$. On the other hand, there is a fiber sequence $$(\haut^+_*(S^1))_{\mathrm{id}} \rightarrow (\haut^+(S^1))_{\mathrm{id}}  \rightarrow S^1,$$ hence the lemma. 
\end{proof}

We now give a sketch of proof of Theorem \ref{thm-hautvsdiff}:

\begin{proof}
    Let $\Sigma$ be an oriented compact connected surface. In the case $\chi(\Sigma)<0$, the equivalence follows from \cite{ee1} and \cite{pinotmcg}. It remains to deal with surfaces with nonnegative Euler characteristic, i.e. the disk, the annulus and the torus. The disk is contractible and Alexander trick applies in dimension $2$, hence both $\Diff_{\partial}(D^2)$ and $\haut_{\partial}^+(D^2)$ are contractible.
    
    Dehn twists generate $\Diff_{\partial}^+(S^1\times I)$, which has contractible connected components, according to \cite{ee} or \cite{pinotmcg}. On the other hand, Dehn twists also generate $\pi_0(\haut_{\partial}^+(S^1\times I))$. The components of $\haut_{\partial}^+(S^1\times I)$ are also contractible since the component of the identity $\haut_{\partial}^+(S^1\times I)$ retracts on $(\Omega\haut^+(S^1))_{\mathrm{id}}$.
    
    The torus is the Eilenberg-Maclane space $K(\bZ^2,1)$. It follows from a result of \cite{gottlieb} that $\haut(K(\bZ^2,1))$ is equivalent, as a monoid, to the semidirect product  
    $K(\bZ^2,1) \rtimes \mathrm{Aut}(\bZ^2).$ By restricting to the orientation-preserving components,  $\haut^+(K(\bZ^2,1))$ is equivalent to $K(\bZ^2,1) \rtimes \mathrm{SL}_2(\bZ).$ The equivalence with self-diffeomorphisms then follows from \cite{ee1} and \cite{pinotmcg}.
        
\end{proof}

We conclude with discussing the exceptional case $g+n=0$, i.e. the manifold is $S^2$. The following proposition, proved in \cite{10.1007/BFb0083829}, compares the spaces $\B\Diff^+(S^2)$ and $\B\haut^+(S^2)$:

\begin{prop} \label{HansenS^2}
    The homotopy fiber of the map $$\mathrm{BSO}(3) \rightarrow \Bhaut^+(S^2)$$ induced by the monoid map $\mathrm{SO}(3) \rightarrow \haut^+(S^2)$ is equivalent to $$\widetilde{\Omega^2S^3}.$$
\end{prop}

We note that $\widetilde{\Omega ^2 S^3}$ is rationally contractible. Hence for any oriented surface $\Sigma$, the monoid map $\Diff_{\partial}^+(\Sigma) \rightarrow \haut_{\partial}^+(\Sigma)$ is a rational homotopy equivalence.

Let $\iota : \text{BSO}(3) \rightarrow \Bhaut^+(\mathrm{S}^2)$ be the map induced by the inclusion morphism $$\text{SO}(3) \rightarrow \haut^+(\mathrm{S}^2).$$ Let $C$ be the homotopy cofiber of $\iota$. The commutative square $$\xymatrix{\text{BSO}(3) \ar[r] \ar[d] & \star \ar[d] \\ \Bhaut^+(\mathrm{S}^2) \ar[r] & C}$$ induces a map on homotopy fibers $c(\iota) : \text{hofib}(\iota) \rightarrow \Omega C.$ Similarly, let $D$ be the (homotopy) cofiber of the map $j : \text{BSO}(2) \rightarrow \Bhaut_*^+(\mathrm{S}^2)$ and let $c(j) : \text{hofib}(j) \rightarrow \Omega D$ be the comparison map.

\begin{lem} \label{comparison-lemma-cofibfib}
    The maps $c(\iota)$ and $c(j)$ are $3$-connected. In particular, both $C$ and $D$ are $2$-connected.
\end{lem}

\begin{proof}
    The maps $\text{BSO}(3) \rightarrow \star$ and $\text{BSO}(3) \rightarrow \Bhaut^+(\mathrm{S}^2)$ are $2$-connected. According to the homotopy excision theorem applied to the pushout square 
\[\begin{tikzcd}
	{\mathrm{BSO}(3)} & \star \\
	{\mathrm{Bhaut}^+(S^2)} & C
	\arrow[from=1-1, to=1-2]
	\arrow["\iota", from=1-1, to=2-1]
	\arrow[from=1-2, to=2-2]
	\arrow[from=2-1, to=2-2]
\end{tikzcd},\] the map $\text{hofib}(\iota) \rightarrow \Omega C$ is $3$-connected. A similar argument proves the map $c(j)$ is $3$-connected.\\
    It follows from \cref{HansenS^2} that $\text{hofib}(\iota)$ is equivalent to $\widetilde{\Omega^2 S^3}$, which is in particular $1$-connected. Hence, both $C$ and $D$ are $2$-connected.\\
\end{proof}

\subsection{Deleting the spheres: Proof of Theorem \ref{thmB}} \label{subsec32} 

In this subsection, we aim to prove Theorem \ref{thmB} stated in the introduction, following a method of Steinebrunner \cite[Theorem B]{Steinebrunner_2020}. In Subsection \ref{3.1bis}, we showed that the difference between the morphism spaces in $\sCobtwo(X)$ and $\sCobtwooo(X)$ lies in the spherical components of the cobordisms. Let $\mathrm{H}$ denote $\mathrm{SO}$ or $\mathrm{SG}$. It suggests to separate the morphisms in $\sCobtwoH(X)$ which correspond to disjoint union of spheres $S^2$ from cobordisms $W$ which do not have any connected component equivalent to $S^2$. In this spirit, we define a reduced category \textcolor{black}{$\mathrm{Cob}^{\mathrm{H,red}}_2(X)$ }and a reduction functor $$\mathrm{red^H}(X) : \sCobtwoH(X) \rightarrow \mathrm{Cob}^{\mathrm{H,red}}_2(X)$$ such that $\mathrm{Cob}^{\mathrm{H,red}}_2(X)$ and $\sCobtwoH(X)$ have the same objects, but the morphisms of $\mathrm{Cob}^{\mathrm{H,red}}_2(X)$ do not have spherical components. We will apply Quillen's Theorem B to study the fiber of the induced map $$\B \sCobtwoH(X) \rightarrow \B \mathrm{Cob}^{\mathrm{H,red}}_2(X).$$ We start with defining reduced morphisms:

\begin{definition}
Let $(A,\tau,f)$ be an element of $\psi_{2,k}^{\mathrm{SO}}(n,1,X)$. We say that $A$ is \textcolor{black}{reduced} if no connected component of $A$ is levelwise diffeomorphic to $\Delta^k \times S^2$. If $A$ is not reduced, we define its reduction $(\mathrm{red}^{\mathrm{SO}}(A),\tau, f_{|\mathrm{red}^{\mathrm{SO}}(A)}) \in \psi_{2,k}^{\mathrm{SO}}(n,1,X)$, where  $\mathrm{red}^{\mathrm{SO}}(A)$ is obtained from $A$ by restricting to the connected components of $A$ which are not diffeomorphic to $\Delta^k \times S^2$.

Similarly, if $(A,\tau, \epsilon, f, f_{\epsilon}) \in \psi_{2,k}^{\mathrm{SO, tub}}(n,1,X)$, we define its \textcolor{black}{reduction} $$\mathrm{red}^{\mathrm{SO,tub}}(A,\tau, \epsilon, f,f_{\epsilon})$$ to be the tuple $$(\mathrm{red}^{\mathrm{SO}}(A),\tau, \epsilon, f_{|\mathrm{red}^{\mathrm{SO}}(A)}, (f_{\epsilon})_{|\mathrm{red}^{\mathrm{SO}}(A)}).$$
Lastly, let $(A,\tau, f)$ be an element of $\psi_{2,k}^{SG}(n,1,X)$. We say that $A$ is \textcolor{black}{reduced} if no connected component of $A$ is levelwise diffeomorphic to $\Delta^k \times S^2 \times \bR^{n-2}$.  If $A$ is not reduced, we define its reduction $(\mathrm{red}^{\mathrm{SG}}(A),\tau, f_{|\mathrm{red}^{\mathrm{SG}}(A)}) \in \psi_{2,k}^{\mathrm{SG}}(n,1,X)$, where  $\mathrm{red}^{\mathrm{SG}}(A)$ is obtained from $A$ by restricting to the connected components of $A$ which are not diffeomorphic to $\Delta^k \times S^2\times \bR^{n-2}$.
\end{definition}

Let $\mathrm{H}$ denote $\mathrm{SO}$ or $\mathrm{SG}$.

\begin{definition}
Let $\psi_{2,k}^{\mathrm{H,red}}(n,1,X)$ be the subset of reduced elements of $\psi_{2,k}^{\mathrm{H}}(n,1,X)$. These form a levelwise sub-simplicial set $\psi_{2,\bullet}^{\mathrm{H,red}}(n,1,X)$ of $\psi_{2,\bullet}^{\mathrm{H}}(n,1,X)$. Letting $n$ go to $\infty$, we get a sub-simplicial set $\psi_{2,\bullet}^{\mathrm{H,red}}(1,X)$ of $\psi_{2,\bullet}^{\mathrm{H}}(1,X)$.\\
Sending an element $(A,\tau,f)$ to its reduction $(\mathrm{red}^{\mathrm{H}}(A),f_{|\mathrm{red}^{\mathrm{H}}(A)})$ defines a map of simplicial sets $$\mathrm{red}^{\mathrm{H}} : \psi_{2,\bullet}^{\mathrm{H}}(1,X) \rightarrow \psi_{2,\bullet}^{\mathrm{H,red}}(1,X).$$
Similarly, we can define a simplicial set $\psi_{2,\bullet}^{\mathrm{SO,tub,red}}$ and a reduction map $$\psi_{2,\bullet}^{\mathrm{SO,tub}} \rightarrow \psi_{2,\bullet}^{\mathrm{SO,tub,red}}.$$
\end{definition}

We now define the reduced cobordism category $\mathrm{Cob}^{\mathrm{H,red}}_2(X)$ as follows. 

\begin{definition}\label{def-3.7-red-cob}
    Let \textcolor{black}{$\mathrm{Cob}_
    2^{\mathrm{H},\mathrm{red}}(X)$} be the simplicial category with:
    \begin{itemize}
    \item its simplicial set of objects is $\psi_{1}^{\mathrm{H}}(X)$, as in $\sCobtwoH(X)$;
     \item its simplicial set of morphisms is $\psi_{2}^{\mathrm{H},\mathrm{red}}(1,X)$;
        \item the source and target maps are as in $\sCobtwoH(X)$;
        \item the composition $(A,\tau, f)\circ (A,\tau',f')$ is given by $$(\mathrm{red}^{\mathrm{H}}(A\cup_{A_{\tau}} (A'+\tau e_1)),\tau+\tau', (f\cup_{|A_{\tau}} f')_{|\mathrm{red}^{\mathrm{H}}(A\cup_{A_{\tau}} (A'+\tau e_1))}).$$ In other words, we first do the composition in $\sCobtwoH(X)$ and then take its reduction.
    \end{itemize}
    There is a natural reduction functor \textcolor{black}{$$\mathrm{red}^{\mathrm{H}}(X) : \sCobtwoH(X) \rightarrow \mathrm{Cob}_
    2^{\mathrm{H},\mathrm{red}}(X)$$} given by the identity on the objects and sending every morphism to its reduction. 
\end{definition}

 Similarly, we define a reduced cobordism category of surfaces with the data of a tubular neighborhood $\mathrm{Cob}_2^{\mathrm{SO,tub,red}}(X)$, as well as a reduction functor: $$\mathrm{red}^{\mathrm{SO,tub}}(X) : \mathrm{Cob}_2^{\mathrm{SO,tub}}(X) \rightarrow \mathrm{Cob}_2^{\mathrm{SO,tub,red}}(X).$$ As in subsection \ref{2.3}, there is a zigzag of functors $$\mathrm{Cob}_2^{\mathrm{SO,red}}(X) \leftarrow \mathrm{Cob}_2^{\mathrm{SO,tub, red}}(X) \rightarrow \mathrm{Cob}_2^{\mathrm{SG,red}}(X).$$
The forgetful functor $\mathrm{Cob}_2^{\mathrm{SO,tub, red}}(X) \rightarrow \mathrm{Cob}_2^{\mathrm{SO,red}}(X) $ induces an equivalence on the nerves, since it does on the non-reduced categories. We prove in the lemma below that the functor $\mathrm{Cob}_2^{\mathrm{SO,tub, red}}(X) \rightarrow \mathrm{Cob}_2^{\mathrm{SG,red}}(X)$ induces an equivalence on the nerves too.

\begin{lem} \label{same-red}
    The induced map $$\B \mathrm{Cob}_2^{\mathrm{SO,tub, red}}(X) \rightarrow \B \mathrm{Cob}_2^{\mathrm{SG,red}}(X)$$ is a weak equivalence of spaces.
\end{lem}

\begin{proof}
    According to Subsection \ref{3.1bis}, if $W$ is reduced, the map $\B \Diff_{\partial}^+(W) \rightarrow \B \haut_{\partial}^+(W)$ is an equivalence. Consequently, the map $$\Map_{\partial}(W,X)\sslash\Diff_{\partial}^+(W) \rightarrow \Map_{\partial}(W,X)\sslash\haut_{\partial}^+(W)
    $$ is an equivalence when $W$ is reduced. Connected components of morphism spaces in $$\mathrm{Cob}_2^{\mathrm{SO,tub, red}}(X)$$ are equivalent to $\Map_{\partial}(W,X)\sslash\Diff_{\partial}^+(W),$ where $W$ is a reduced surface. On the other hand, connected components of the morphism spaces in $\mathrm{Cob}_2^{\mathrm{SG, red}}(X)$ are equivalent to $$\Map_{\partial}(W,X)\sslash\haut_{\partial}^+(W)$$ where $W$ is a reduced cobordism. Consequently, the functor $\mathrm{Cob}_2^{\mathrm{SO,tub, red}}(X) \rightarrow \mathrm{Cob}_2^{\mathrm{SG,red}}(X)$ induces an equivalence on objects and on morphism spaces, hence on geometric realization.
\end{proof}

On the other hand, we define a subcategory of $\sCobtwoH(X)$ which only contains the spherical morphisms:

\begin{definition}
    Let \textcolor{black}{$\mathrm{Sph^H}(X)$} denote the submonoid of $\sCobtwoH(X)(\emptyset, \emptyset)$ consisting of morphisms $(A,\tau,f) \in \mathrm{Mor}(\sCobtwoH(X))$ such that $A$ is homotopy equivalent to a finite disjoint union of $S^2$. 
    The morphism $u(X) : \B \sCobtwooo(X) \rightarrow \B \sCobtwo(X)$ restricts to a morphism $\B \mathrm{Sph^{SO}}(X) \rightarrow \B \mathrm{Sph^{SG}}(X)$.
\end{definition}

In the next Subsection \ref{subsec33}, we will prove that the homotopy fiber of the reduction map $$\B\mathrm{red}^{\mathrm{H}}(X) : \B\sCobtwoH(X) \rightarrow \B\mathrm{Cob}_
    2^{\mathrm{H},\mathrm{red}}(X)$$ is determined by the spherical morphisms. We state below the result, which we will prove in Subsection \ref{subsec33}:
    
\begin{thm} \label{thm-fiber-reduction}
    The homotopy fiber of $\B\mathrm{red}^{\mathrm{H}}(X) : \B\sCobtwoH(X) \rightarrow \B\mathrm{Cob}_
    2^{\mathrm{H},\mathrm{red}}(X)$ is $\B \mathrm{Sph}^{\mathrm{H}}(X).$
\end{thm} 

In Subsection \ref{2.4}, we showed that the functor $\B \sCobtwo(-)$ factors through $\mathrm{Sp}$. In particular, the functor $\Omega_{\emptyset} \B \sCobtwo(-)$ does as well. We now construct a non-connective delooping $$\mathrm{PH}(2,-) :\mathcal{S} \rightarrow \mathrm{Sp}$$ of $\Omega_{\emptyset}\B \sCobtwo(-)$, such that it is compatible with the non-connective delooping $(\mathrm{MTSO}(2) \otimes \susplus -)$ of $\Omega_{\emptyset} \B \sCobtwooo(-)$ given by the Galatius-Madsen-Tillmann-Weiss Theorem.  

In \cite{GMTW} and \cite{grw9}, they construct an equivalence $\alpha : \B \sCobtwooo(X) \rightarrow \Omega^{\infty-1} (\mathrm{MTSO}(2) \otimes \susplus X)$, called the scanning map. The adjoint of the composite $$\alpha : \Map(S^2,X)\sslash\mathrm{SO}(3) \rightarrow \Omega_{\emptyset} \B \sCobtwooo(X) \rightarrow \lop (\mathrm{MTSO}(2) \otimes \susplus X)$$ gives a map of spectra $$\mathrm{PT}_{S^2}^{\Diff}(X)  : \susplus \Map(S^2,X)//\mathrm{SO}(3) \rightarrow \lop (\mathrm{MTSO}(2) \otimes \susplus X).$$ We say the latter is a parametrized Pontryagin-Thom construction map, see Section \ref{4} for more details. 

\begin{definition}
    Let \textcolor{black}{$\mathrm{PH}(2,X)$} be the spectrum defined by the following pushout
\[\begin{tikzcd}
	{\Sigma^{\infty}_+ \mathrm{Map}(S^2,X)\sslash\mathrm{SO}(3)} & {\mathrm{MTSO}(2)\otimes \susplus X} \\
	{\Sigma^{\infty}_+ \mathrm{Map}(S^2,X)\sslash\mathrm{haut}^+(S^2)} & \textcolor{black}{{\mathrm{PH}(2,X)}}
	\arrow[""{name=0, anchor=center, inner sep=0}, "{\mathrm{PT}^{\mathrm{Diff}}_{S^2}(X)}"', from=1-1, to=1-2]
	\arrow[from=1-1, to=2-1]
	\arrow[from=1-2, to=2-2]
	\arrow[from=2-1, to=2-2]
	\arrow["\ulcorner"{anchor=center, pos=0.125, rotate=180}, draw=none, from=2-2, to=0]
\end{tikzcd}\]
This square is natural in $X$ and defines a functor $\mathrm{PH}(2,-) :\mathcal{S} \rightarrow \mathrm{Sp}$.
\end{definition}

We now reformulate \cref{thmB}, before giving its proof: 

\begin{thm}[Theorem \ref{thmB} in the Introduction]
    The square \eqref{main-square-thmB} is homotopy cartesian. Moreover, the pushout square \[\begin{tikzcd}
	{\Sigma^{\infty}_+ \mathrm{Map}(S^2,X)\sslash\mathrm{SO}(3)} & {\mathrm{MTSO}(2)\otimes \susplus X} \\
	{\Sigma^{\infty}_+ \mathrm{Map}(S^2,X)\sslash\mathrm{haut}^+(S^2)} & \textcolor{black}{{\mathrm{PH}(2,X)}}
	\arrow[""{name=0, anchor=center, inner sep=0}, "{\mathrm{PT}^{\mathrm{Diff}}_{S^2}(X)}"', from=1-1, to=1-2]
	\arrow[from=1-1, to=2-1]
	\arrow[from=1-2, to=2-2]
	\arrow[from=2-1, to=2-2]
	\arrow["\ulcorner"{anchor=center, pos=0.125, rotate=180}, draw=none, from=2-2, to=0]
\end{tikzcd}\] is equivalent after taking $\lop$ to the square \eqref{main-square-thmB}.
\end{thm}

\begin{proof}[Proof of Theorem \ref{thmB}]
The proof is in two parts. In a first part, we show that the square \eqref{main-square-thmB} is a homotopy pullback and in a second part, we show that the spectrum $\mathrm{PH}(2,X)$ is indeed a delooping of $\Omega_{\emptyset} \B \PCobdG(X)$. 

On the one hand, according to Theorem \ref{thm-fiber-reduction}, there is a map of fiber sequences
\begin{equation} \label{square-proof}
    \begin{tikzcd}
	{\mathrm{BSph}^{\mathrm{SO}}(X)} & {\mathrm{BCob}_2^{\mathrm{SO}}(X)} & {\mathrm{BCob}_2^{\mathrm{SO,red}}(X)} \\
	{\mathrm{BSph}^{\mathrm{SG}}(X)} & {\mathrm{BCob}_2^{\mathrm{SG}}(X)} & {\mathrm{BCob}_2^{\mathrm{SG,red}}(X)}
	\arrow[from=1-1, to=1-2]
	\arrow[from=1-1, to=2-1]
	\arrow[from=1-2, to=1-3]
	\arrow["{u(X)}", from=1-2, to=2-2]
	\arrow["{(a)}", from=1-3, to=2-3]
	\arrow[from=2-1, to=2-2]
	\arrow["{(b)}"{description}, draw=none, from=2-2, to=1-1]
	\arrow[from=2-2, to=2-3]
\end{tikzcd}.
\end{equation} According to Lemma \ref{same-red}, the right hand-side map labelled $(a)$ in the diagram \eqref{square-proof} is an equivalence, hence the left square labelled by $(b)$ is a homotopy pullback. It is still a homotopy pullback square after looping once. 

Concerning the second part of the statement, we recall that in \cref{main-prop-4.2}, we constructed a connective delooping $C^{\cSG}(2,-)$ and a natural equivalence $\B \sCobtwo(-) \Rightarrow \lop  C^{\cSG}(2,-)$. After looping once and taking connected truncations, we have a natural equivalence $$\Omega_{\emptyset} \B \sCobtwo(-) \Rightarrow \lop \tau_{\geq 0}(\Sigma^{-1} C^{\cSG}(2,-)).$$ Similarly, we have a natural equivalence $$\Omega_{\emptyset}\B \sCobtwooo(-) \Rightarrow \lop \tau_{\geq 0}(\mathrm{MTSO}(2) \otimes \susplus -).$$

It remains to identify the group-completion $\Omega_{\emptyset} \B \mathrm{Sph}^{\mathrm{H}}(X)$ with the free infinite loop space $\mathrm{Q}_+ \Map(S^2,X)\sslash\mathrm{\mathrm{H}}(S^2),$ where $\mathrm{H}, \mathrm{H}(S^2)$ respectively denote $\mathrm{SO}, \Diff^+(S^2)$ or $\mathrm{SG}, \haut^+(S^2)$. We note that $\mathrm{H}(S^2 \times \{1,...,n\})$ is equivalent to the wreath product monoid $(\mathrm{H}(S^2) \wr \Sigma_n)$. Then, we observe that we have natural equivalences $$\Map(S^2\times \{1,...,n\},X)\sslash \mathrm{H}(S^2 \times \{1,...,n\}) \simeq (\Map(S^2,X)\sslash \mathrm{H}(S^2))^n\sslash\Sigma_n.$$ According to \cite[Proposition $3.6$]{SEGAL1974293}, the group completion $\Omega \B (\bigsqcup_{j\geq 0} X^n\sslash\Sigma_n)$ is equivalent to the free infinite loop space $\lop \susplus X$. Moreover, it follows that the map $\Omega_{\emptyset}\B \mathrm{Sph}^{H}(X) \rightarrow \Omega_{\emptyset} \B \sCobtwoH(X)$ is equivalent to the free infinite loop-map $$\lop \susplus (\Map(S^2,X)//\mathrm{H}(S^2)) \rightarrow \Omega_{\emptyset} \B \sCobtwoH(X).$$

Finally, the adjoint $\mathrm{PT}_{S^2}^{\Diff} :\susplus (\Map(S^2,X)\sslash\mathrm{SO}(3)) \rightarrow \mathrm{MTSO}(2) \otimes \susplus X $ is also equivalent after taking $\lop$ to the free infinite loop map $\lop \susplus (\Map(S^2,X)\sslash\mathrm{SO}(3)) \rightarrow \lop (\mathrm{MTSO}(2)\otimes \susplus X)$. Hence, the map $$\lop \susplus (\Map(S^2,X) \sslash\mathrm{SO}(3)) \rightarrow \lop (\mathrm{MTSO}(2)\otimes \susplus X)$$ is equivalent to the map $\Omega_{\emptyset} \B \mathrm{Sph}^{SO}(X) \rightarrow \Omega_{\emptyset} \B \sCobtwooo(X)$.

Since $\susplus \Map(S^2,X)\sslash\mathrm{SO}(3)$ is connective, the map $$\susplus \Map(S^2,X)\sslash\mathrm{SO}(3) \rightarrow \mathrm{MTSO}(2)\otimes \susplus X$$ factors through the connective cover $\susplus \Map(S^2,X)\sslash\mathrm{SO}(3) \rightarrow \tau_{\geq 0}(\mathrm{MTSO}(2)\otimes \susplus X)$. 

The square \eqref{main-square-thmB} then deloops to the following commutative square of connective spectra:

\begin{equation} \label{square-spectra}
\begin{tikzcd}
	{\Sigma_+^{\infty} \mathrm{Map}(S^2,X)\sslash\mathrm{SO}(3)} & {\tau_{\geq 0}(\mathrm{MTSO}(2)\otimes \Sigma^{\infty}_+X) } \\
	{\Sigma_+^{\infty} \mathrm{Map}(S^2,X)\sslash\mathrm{haut}^+(S^2)} & {\tau_{\geq 0}( \Sigma^{-1}C^{\mathrm{SG}}(2,X))}
	\arrow[from=1-1, to=1-2]
	\arrow[from=1-1, to=2-1]
	\arrow["{\tilde{u}(X)}", from=1-2, to=2-2]
	\arrow[from=2-1, to=2-2]
\end{tikzcd}.
\end{equation}
The square \eqref{square-spectra} lives in $\mathrm{Sp}^{\geq 0}$ and becomes the pullback square \eqref{square-thmB} in $\mathcal{S}$ after taking $\Omega^{\infty}$. We can then conclude the square \eqref{square-spectra} is actually a pullback square, equivalently a pushout square, in $\mathrm{Sp}^{\geq 0}$, hence $\mathrm{Sp}$.

Let $D^{\mathrm{SG}}(2,X)$ denote the following pushout:

\[\begin{tikzcd}
	{\tau_{\geq 0}( \mathrm{MTSO}(2)\otimes \susplus X)} & {\mathrm{MTSO}(2)\otimes \susplus X} \\
	{\tau_{\geq 0}(\Sigma^{-1}C^{\mathrm{SG}}(2,X))} & {D^{\mathrm{SG}}(2,X)}
	\arrow[from=1-1, to=1-2]
	\arrow[from=1-1, to=2-1]
	\arrow[from=1-2, to=2-2]
	\arrow[from=2-1, to=2-2]
\end{tikzcd}.\]
It becomes a pullback square after applying $\Omega^{\infty}$. Hence, the map $$\lop \tau_{\geq 0}(\Sigma^{-1}C^{\mathrm{SG}}(2,X))\rightarrow \lop \mathrm{D^{SG}}(2,X)$$ is an equivalence since the top map is. 

We can combine both squares into the following diagram: 

\[\begin{tikzcd}
	{\Sigma^{\infty}_+ \mathrm{Map}(S^2,X)\sslash\mathrm{SO}(3)} & {\tau_{\geq 0}( \mathrm{MTSO}(2)\otimes \susplus X)} & {\mathrm{MTSO}(2)\otimes \susplus X} \\
	{\Sigma^{\infty}_+ \mathrm{Map}(S^2,X)\sslash\mathrm{haut}^+(S^2)} & {\tau_{\geq 0}(\Sigma^{-1}C^{\mathrm{SG}}(2,X))} & {D^{\mathrm{SG}}(2,X)}
	\arrow[""{name=0, anchor=center, inner sep=0}, from=1-1, to=1-2]
	\arrow[from=1-1, to=2-1]
	\arrow[""{name=1, anchor=center, inner sep=0}, from=1-2, to=1-3]
	\arrow[from=1-2, to=2-2]
	\arrow[from=1-3, to=2-3]
	\arrow[from=2-1, to=2-2]
	\arrow[from=2-2, to=2-3]
	\arrow["\ulcorner"{anchor=center, pos=0.125, rotate=180}, draw=none, from=2-2, to=0]
	\arrow["\ulcorner"{anchor=center, pos=0.125, rotate=180}, draw=none, from=2-3, to=1]
\end{tikzcd}.\]
Since the left and right squares are pushouts, we conclude the outer square is also a pushout. According to our previous discussion, the top composite map is equivalent to $\mathrm{PT}^{\Diff}_{S^2}(X)$. Finally, by definition, $\mathrm{PH}(2,X)$ is equivalent to $D^{\mathrm{SG}}(2,X)$. Moreover after taking $\lop$, the outer square is equivalent to the square \eqref{main-square-thmB}, which concludes the proof.
\end{proof}

In the proof of \cref{thmB}, we actually constructed a natural transformation 
$$\Omega_{\emptyset} \B \sCobtwo(-) \Rightarrow \lop \mathrm{PH}(2,-)$$ which factors as
\[\begin{tikzcd}
	{\Omega_{\emptyset} \B \sCobtwo(-)} & {\lop \tau_{\geq 0} (\Sigma^{-1}\mathrm{C}^{\mathrm{SG}}(2,-))} & {\lop \mathrm{PH}(2,-)}
	\arrow[from=1-1, to=1-2]
	\arrow[from=1-2, to=1-3]
\end{tikzcd}.\] 

In the following corollary, we show this natural equivalence deloops: 

\begin{cor} \label{cor-delooping-nerve}
    The equivalence $$\Omega_{\emptyset} \B \sCobtwo(-) \Rightarrow \lop \mathrm{PH}(2,-)$$ delooops to a natural equivalence $$\eta : \B \sCobtwo(-) \Rightarrow \lop \Sigma \mathrm{PH}(2,-).$$
\end{cor}

\begin{proof}
    Since $\Omega_{\emptyset} \B\sCobtwo(X)$ is equivalent to $\lop \mathrm{PH}(2,X)$, it suffices to show that $$\pi_0(\B \sCobtwo(X)) \cong \pi_{-1}(\mathrm{PH}(2,X)).$$ Since the homotopy categories $\mathrm{h}\sCobtwo(X)$ and $\mathrm{h}\sCobtwooo(X)$ are equivalent, there is a group isomorphism $\pi_0(\B \sCobtwooo(X)) \rightarrow \pi_0(\B\sCobtwo(X)).$ It suffices then to show that $\pi_{-1}(\mathrm{PH}(2,X))$ is isomorphic to $\pi_{-1}(\mathrm{MTSO}(2)\otimes \susplus X)$. 

    According to Theorem \ref{thmB}, the fiber of $\mathrm{MTSO}(2) \otimes \susplus X \rightarrow \mathrm{PH}(2,X)$ is equivalent to $ \Sigma^{-1} \sus C(X)$, where $C(X)$ is the homotopy cofiber of $$\Map(S^2,X)\sslash\mathrm{SO}(3) \rightarrow \Map(S^2,X)\sslash\haut^+(S^2).$$ 

    As in Lemma \ref{comparison-lemma-cofibfib}, we can use the homotopy excision theorem to deduce $C(X)$ is $2$-connected, hence $\Sigma^{-1} \sus C(X)$ is connective, which concludes the proof. 
\end{proof}

\subsection{Quillen's Theorem B for the reduction functor} \label{subsec33}

In this subsection, we aim to prove \cref{thm-fiber-reduction}. In other words, we want to show the homotopy fiber of the map $\B \mathrm{red^{SG}}(X) :\B\sCobtwo(X) \rightarrow \B\cobred(X)$ is equivalent to $\B \mathrm{Sph^{SG}}(X)$ for every space $X$, and similarly that the homotopy fiber of $\B \sCobtwooo(X) \rightarrow \B \mathrm{Cob}_2^{\mathrm{SO,red}}(X)$ is $\B \mathrm{Sph}^{\mathrm{SO}}(X)$. 

The main ingredient for studying the homotopy fiber of the geometric realization of the functor $ \mathrm{red^{SG}} : \sCobtwo \rightarrow \mathrm{Cob_2^{SG,red}}$ is Quillen's Theorem B. It gives a series of conditions on a functor so that the homotopy fiber of the geometric realization is equivalent to the classifying space of the \emph{genuine fiber} of the functor, which we now describe. 

For readability, assume $X$ is a point. As observed in \cref{rk-non-unit}, the category $\sCobtwo$ does not admit strict identity morphisms. However, it admits what we call \emph{weak units}. Let $U$ be a $0$-simplex of the simplicial set of objects of $\sCobtwo$. Then, the pair $(U\times \bR, \tau)$ is an endomorphism of $U$. Precomposing a morphism $(W,\tau')$ with $(U\times \bR,\tau)$ does not give back exactly $(W,\tau')$, but the resulting thickening $(U\times \bR \cup_{W_0} W )$ is equivalent to $W$. We say that such cylindrical endomorphisms $(U\times \bR,\tau)$ are weak units of the object $U$. As $0$-simplices, both $U$ and $(U\times \bR,\tau)$ generate constant subsimplicial sets of $\mathrm{Ob}(\sCobtwo)$ and $\mathrm{Mor}(\sCobtwo)$, which we also denote by $U$ and $(U \times \bR,\tau)$. We consider the pair $(U,(U\times \bR,\tau))$. We now describe the genuine fiber $(\mathrm{red^{SG}})^{-1}(U,(U\times \bR,\tau))$ of the functor $\mathrm{red^{SG}}$ at $(U,(U\times \bR,\tau))$. It is the subcategory of $\sCobtwo$ with objects $V$ such that $\mathrm{red^{SG}}(V)$ is $U$ and morphisms $(W,\tau')$ such that $\mathrm{red^{SG}}(W,\tau')$ is equivalent to $(U\times \bR,\tau)$. Since $(U\times \bR,\tau) \circ (U\times \bR,\tau)$ is equivalent to $(U\times \bR,\tau)$, this indeed defines a subcategory $(\mathrm{red^{SG}})^{-1}(U,(U\times \bR,\tau))$ of $\sCobtwo$.

We now define the fiber $(\mathrm{red^{SG}})^{-1}(U,(U\times \bR,\tau))$. By definition of the reduction functor, the objects of $(\mathrm{red^{SG}})^{-1}(U,(U\times \bR,\tau))$ are just given by the constant simplicial set $U$. The morphisms are pairs $(W,\tau')$ such that $\mathrm{red^{SG}}(W,\tau')$ is equivalent to $(U\times \bR, \tau)$. In fact, $\tau'$ can be any real number and $W$ is equivalent to the disjoint union of the cylinder $U\times \bR$ and a finite disjoint union of spherical components. In particular, in the case $U=\emptyset$, the genuine fiber $(\mathrm{red^{SG}})^{-1}(\emptyset,(\emptyset,\tau))$ is exactly $\sphG$. 

In this subsection, we use a version of Quillen's Theorem B for topological categories (in other words, categories internal to topological spaces $\mathrm{Top}$), as proved by Steinebrunner in \cite[Theorem $A$]{Steinebrunner_2022}. Before that, we explain how to obtain a topological category from a simplicial category:

\begin{definition}
    Let $\mathcal{C}$ be a category internal to $\mathrm{sSet}$, with simplicial sets of objects $\mathrm{Ob}(\mathcal{C})$, morphisms $\mathrm{Mor}(\mathcal{C})$. We can define a category $|\mathcal{C}|$ internal to $\mathrm{Top}$ by taking the space of objects to be $|\mathrm{Ob}(\mathcal{C})|$ and morphisms $|\mathrm{Mor}(\mathcal{C})|$.
\end{definition}
 
Most of the proof of \cref{thm-fiber-reduction} is adapted from \cite[Theorem B]{Steinebrunner_2022} for smooth cobordism categories. In particular, in what follows, we focus on the Poincaré cobordism category. For readability, we may work with the cobordism category over the point $\sCobtwo$. Without too much work, we can generalize the following to the cobordism category $\sCobtwo(X)$ over $X$. The following lemma describes equivalences in the categories $\sCobtwo$ and $\cobred$:
  
\begin{lem} \label{units-3.17}
    Let $(W,\tau)$ be a morphism in $\sCobtwo$ from $W_0$ to $W_{\tau}$. The precomposition map $$- \circ (W,\tau) : \sCobtwo(W_{\tau},B) \rightarrow \sCobtwo(W_0,B)$$ and the postcomposition map $$ (W,\tau) \circ - : \sCobtwo(B,W_0) \rightarrow \sCobtwo(B,W_{\tau})$$ are equivalences for all objects $B$ if and only if $W$ is diffeomorphic to a cylinder $W_0\times \bR$. If this is the case, we say $(W, \tau)$ is an \emph{\textcolor{Black}{equivalence}}.

    The functor $\mathrm{red^{SG}} :\sCobtwo \rightarrow \cobred$ sends equivalences to equivalences. 
\end{lem}

\begin{proof}
    We show that if precomposition and postcomposition with $(W,\tau)$ is an equivalence, then $W$ is equivalent to a cylinder. The other direction is immediate. Take $B$ to be $W_0$ and the precomposition map $- \circ (W,\tau) : \sCobtwo(W_{\tau},W_0) \rightarrow \sCobtwo(W_0,W_0)$. Then there exists a morphism $(V,\tau')$ from $W_{\tau}$ to $W_0$ such that $W \cup_{W_{\tau}} (V+ \tau e_1)$ is equivalent to $W_0 \times I$. On the other hand, by assumption, the postcomposition map $(W,\tau) \circ- : \sCobtwo(W_{\tau}, W_0) \rightarrow \sCobtwo(W_{\tau},W_{\tau})$ is also an equivalence. Then, there exists a morphism $(V',\tau'')$ from $W_0$ to $W_{\tau}$, such that $W \cup_{W_0} (V'+\tau'e_1)$ is equivalent to $W_{\tau} \times I$. It follows from these two points that $W$ is a cylinder.
\end{proof} 

We now show the following proposition: 

\begin{prop} \label{cat-fibrant}
    The source-target map $(s,t) : \mathrm{Mor}(\sCobtwo) \rightarrow  \mathrm{Ob}(\sCobtwo) \times \mathrm{Ob}(\sCobtwo) $ is a Kan fibration. 
\end{prop}

\begin{proof}
For $U$ a manifold with boundary, let $\mathrm{Sub}(U,\bR^n)$ denote the simplicial set $$S_{\bullet} \Emb(U,\bR^n)/S_{\bullet} \Diff(U).$$ There is a map $\partial : \mathrm{Sub}(U,\bR^n) \rightarrow \mathrm{Sub}(\partial U \times \bR,\bR^n)$, which corresponds to taking a collar of the boundary.

To see why $\sCobtwo$ is fibrant, it suffices to show the map $\partial : \mathrm{Sub}(U,\bR^n) \rightarrow \mathrm{Sub}(\partial U \times \bR, \bR^n)$ is a Kan fibration, where $U$ is an open $n$-dimensional thickening (potentially with boundary) of a Poincaré pair $(P,Q)$.

Let $U$ be such a manifold. For $k\geq 2, j=0,1,\ldots, i-1,i+1,\ldots, k$, let $B_j$ be a $(k-1)$-simplex of $\mathrm{Sub}(U,\bR^n)$, such that $d_l(B_j)= d_{j-1}(B_l)$ for $l<j$. Let $A$ be a $k$-simplex of $\mathrm{Sub}(\partial U \times \bR,\bR^n)$ such that $\partial B_j =d_j(A)$ for all $j$. We wish to construct a $k$-simplex $B$ of $\mathrm{Sub}_{\partial}$ such that its collared boundary is $A$.

Let $r : \Delta^k \rightarrow \Lambda_i^k$ be a retraction of horn inclusion. For each $\sigma$, there is a unique straight path $\gamma(\sigma)$ from $\sigma$ to $r(\sigma)$. 

For each $\sigma\in \Delta^k$, let $B_{\sigma}$ be $B_{r(\sigma)}\cup_{\partial B_{r(\sigma)}} A_{\gamma(\sigma)}$, where $A_{\gamma(\sigma)}$ denotes the restriction of the family $A$ to the path $\gamma(\sigma)$. The parametrized collection $(B_{\sigma})$ defines a $k$-simplex of $\mathrm{Sub}(U,\bR^n)$. Up to rescaling the collars and $A_{\gamma(\sigma)}$, the collar of $B$ is $A$ and for $\sigma \in \Lambda_i^k$, $B_{\sigma}$ coincides with $(B_j)_{\sigma}$ for some $j$, which concludes the proof.

\end{proof}

\begin{cor} \label{cor-cat-fibrant}
    The map $(s,t) : \mathrm{Mor}(\cobred) \rightarrow \mathrm{Ob}(\cobred) \times \mathrm{Ob}(\cobred)$ is a Kan fibration.
\end{cor}

\begin{proof}
    The map $\partial : \mathrm{Sub}(U;\bR^n) \rightarrow \mathrm{Sub}(\partial U\times \bR, \bR^n)$ is a Kan fibration for every open thickening $U$ of a Poincaré pair, including ones with no spherical components.
\end{proof}

We now prove the following decomposition lemma:

\begin{lem} \label{lem-prod-clos-red}
Let $P$ and $P'$ be two objects in $\sCobtwo$. Projecting on reduced and spherical components of the morphisms yields an equivalence: $$\sCobtwo(P,P') \rightarrow  \cobred(P,P') \times \mathrm{Sph}^{\mathrm{SG}}.$$   
\end{lem}

\begin{proof}
    The proof follows from the subsequent observation: the map induced by projection on components
    $$\B \haut_{\partial}(W \sqcup S^2 \sqcup \ldots \sqcup S^2) \rightarrow \Bhaut_{\partial}(W) \times \Bhaut(S^2 \sqcup \ldots \sqcup S^2)$$ is an equivalence whenever $W$ is a reduced surface. 
\end{proof}

We now show the reduction functor is a \textcolor{Black}{local fibration}, which is equivalent to the following proposition:

\begin{prop} \label{local-fib-red} 
    The reduction map on objects  $\mathrm{Ob}(\sCobtwo) \rightarrow \mathrm{Ob}(\cobred)$ and morphisms $\mathrm{Mor}(\sCobtwo) \rightarrow \mathrm{Mor}(\cobred)$ are Kan fibrations.
\end{prop}

\begin{proof}
    On objects, the reduction functor is the identity hence it is a Kan fibration. Let us now show it for morphisms. Let $B_j$ be $(k-1)$-simplices in the morphisms simplicial set $\mathrm{Mor}(\sCobtwo)$, for $j=1,\ldots, i-1,i+1,\ldots k$. Assume the $B_j$ have compatible faces, i.e. $d_l(B_j)=d_{j-1}(B_l)$ for $l<j$. Let $A$ be a $k$-simplex of $\mathrm{Mor}(\cobred)$ such that the $j$-th face of $A$, $d_j(A)$, is equivalent to the reduction $\mathrm{red}^{\mathrm{SG}}(B_j)$. Geometrically, $A$ is a $\Delta^k$-parametrized family of reduced morphisms, and each $B_j$ is a $\Delta^{k-1}$-parametrized families of non-reduced morphisms. We wish to construct a $k$-simplex $B$ of $\mathrm{Mor}(\sCobtwo)$ such that its reduction is equal to $A$, and its face $d_j(B)$ is $B_j$ for all $j\neq i$. Since the simplicial set $\mathrm{Mor}(\sCobtwo)$ is obtained as a filtered colimits of simplicial sets $\mathrm{Mor}(\sCobtwo(n))$, we can assume $(B_j)_j, A$ are respectively $(k-1), k$-parametrized subsets of $\Delta^{k-1}\times \bR^{n}, \Delta^k \times \bR^n$. 

    Up to slightly modifying the definition of $\sCobtwo$, we assume that the connected components of $A$ (resp. $B_j$) are at distance at least $1$ from each others in $\Delta^k \times \bR^n$.

     We take for each $j$, the restriction $\mathrm{sph}(B_j)$ to the spherical components of $B_j$. In particular, $\mathrm{sph}(B_j)$ is a $(k-1)$-simplex of the monoid $\mathrm{Sph}^{\mathrm{SG}}$. According to Lemma \ref{lem-2.2-0.2}, the simplicial sets $\mathrm{Sub}(U,\bR^n)$ are Kan complexes. Consequently, we can find a $k$-simplex $S$ of $\mathrm{Sph}^{\mathrm{SG}}$ such that $d_j(S) = B_j$.

     Let $\lambda : \Delta^k \rightarrow [0,\infty) $ be the continuous function such that:
    $\lambda(\sigma) = d(\sigma, \Lambda_i^k)$, where $d$ is a metric on $\Delta^k$. In particular, for $\sigma \in \Lambda_i^k$, $\lambda(\sigma)=0$. Since the $(B_j)_{\sigma}$ and the $S_{\sigma}$ depend continuously in $\sigma$, we can assume that the $S_{\sigma}$ is at distance at most $\frac{1}{2}$ from  $A_{\sigma}$ for $\sigma$ close enough to $\Lambda_i^k$, i.e. if $\lambda(\sigma)<\epsilon$, for some $\epsilon>0$.

    We would like to construct the $k$-simplex $B$ by taking the union of $A$ with $S$, however they may not be disjoint. Instead, we embed $A$, $B_j$ and $S$ in $\bR^{n+1}$, up to thickening them by crossing with an interval in the orthogonal direction.
    
    Since the connected components of $(B_j)_{\sigma}$ are at distance at least $1$ from each other, in particular the components of $\mathrm{Sph}(B_j)$ are at distance at least $1$ from $\mathrm{red^{SG}}(B_j)$. 
    
    We now construct the $k$-simplex $B \subset \Delta^k \times \bR^{n+1}$. We define $B_{\sigma}$ for $\sigma\in \Delta^k$ as follows:
        $$B_{\sigma} = A_{\sigma} \times (-1,1) \cup S_{\sigma} \times (-1+\frac{2}{\epsilon}\lambda(\sigma), 1+\frac{2}{\epsilon} \lambda(\sigma)).$$ When $\sigma \in \Lambda_i^k$, then $\lambda(\sigma)=0$ and $S_{\sigma} =(B_j)_{\sigma}$ for some $j$. When $\lambda(\sigma) \geq \epsilon$, the $A_{\sigma}$ and $S_{\sigma}$ have been made disjoint. When $\lambda(\sigma)< \epsilon$, the $S_{\sigma}$ and $A_{\sigma}$ are at distance at most $\frac{1}{2}$ by assumption, hence the $A_{\sigma}$ and $S_{\sigma}$ are also disjoint. 
            
        In the filtered colimit described in Subsection \ref{2.3}, the simplices $A$ and $A \times (-1,1)$ are identified. Moreover, since we are only adding spherical components, the reduction of $B$ and the reduction of $A$ agree, which concludes the proof. 
\end{proof}

From now on, we may confuse the simplicial categories $\sCobtwoH(X), \mathrm{Cob}_2^{\mathrm{SG,red}}(X), \mathrm{Sph^H}(X)$ with their associated topological categories $|\sCobtwoH(X)|, |\mathrm{Cob}_2^{\mathrm{SG,red}}(X)|, |\mathrm{Sph^H}(X)|$. We can now conclude the proof of Theorem \ref{thm-fiber-reduction}:

\begin{proof}[Proof of Theorem \ref{thm-fiber-reduction}]
    We start by discussing the Poincaré case. 

    Let $\mathcal{C}$ be the subcategory of $\cobred$ with one object $\emptyset$ and morphisms are pairs $(\emptyset, \tau)$ where $\tau>0$. By definition, the reduction functor $\sCobtwo \rightarrow \cobred$ restricts to a functor $\sphG \rightarrow \mathcal{C}$ on $\sphG$. After taking the nerve, we get the following commutative diagram:   
    \begin{equation} \label{square-quillen-B-proof}
\begin{tikzcd}
	{\mathrm{BSph^{SG}}} & {\mathrm{BCob}_2^{\mathrm{SG}}} \\
	{\mathrm{B}\mathcal{C}} & {\mathrm{BCob}_2^{\mathrm{SG,red}}}
	\arrow[from=1-1, to=1-2]
	\arrow[from=1-1, to=2-1]
	\arrow[from=1-2, to=2-2]
	\arrow[from=2-1, to=2-2]
\end{tikzcd}.
\end{equation}

Since $\mathcal{C}$ is equivalent to the terminal category, $\B \mathcal{C}$ is contractible. In particular, if we manage to prove the square \eqref{square-quillen-B-proof} is Cartesian, then we can conclude the proof of \cref{thm-fiber-reduction}. 

We note that the category $\mathcal{C} \times_{\cobred} \sCobtwo $ has one object $\emptyset$ and its space of morphisms is the space of pairs $(W,\tau)$ such that the reduction of $W$ is empty. In other words, the category $\mathcal{C} \times_{\cobred} \sCobtwo $ is equivalent to $\sphG$. According to \cite[Definition $1.2$]{Steinebrunner_2022}, if the functor $\mathrm{red^{SG}}$ is a \emph{realization fibration}, then the square \eqref{square-quillen-B-proof} is Cartesian.

The geometric realization of a Kan fibration is a Serre fibration. It then follows from \cref{cor-cat-fibrant} that the category $\cobred$ is \emph{fibrant} in the sense of \cite[Definition $5.6$]{Steinebrunner_2022}. It follows from \cref{local-fib-red} that the functor $\mathrm{red^{SG}} : \sCobtwo \rightarrow \cobred$ is a \emph{local fibration} in the sense of \cite[Definition $5.6$]{Steinebrunner_2022}.

According to \cref{units-3.17}, equivalences in the categories $\sCobtwo$ and $\cobred$ are morphisms $(W,\tau)$ such that $W$  is equivalent to a cylinder. Since a cylinder is already reduced, the reduction functor $\mathrm{red}^{\mathrm{SG}}$ takes equivalences to equivalences. On the other hand, an endomorphism $(W,\tau)$ of an object $U$ is said to be a \emph{weak unit} in the sense of \cite[Definition $5.4$]{Steinebrunner_2022} if it is an equivalence and if $(W,\tau) \circ (W,\tau)$ is equivalent to $(W,\tau)$. The latter condition is satisfied when $W$ is a cylinder. In particular, every object $U$ of $\sCobtwo$ and $\cobred$ admits weak units. The categories $\sCobtwo$ and $\cobred$ are then \emph{weakly unital} in the sense of \cite[Definition $5.4$]{Steinebrunner_2022}. The reduction functor $\mathrm{red^{SG}}$ sends weak units to weak units, hence it is \emph{weakly unital} in the sense of \cite[Definition $5.4$]{Steinebrunner_2022}.

According to \cite[Theorem $A$]{Steinebrunner_2022}, if $\sCobtwo$ and $\cobred$ are \emph{weakly unital}, $\cobred$ is \emph{fibrant}, the functor  $\mathrm{red^{SG}} : \sCobtwo \rightarrow \cobred$ is \emph{weakly unital}, a \emph{local fibration}, and \emph{locally Cartesian} and \emph{locally coCartesian}, then the functor $\mathrm{red^{SG}}$ is a \emph{realization fibration} in the sense of \cite[Definition $1.2$]{Steinebrunner_2022}. Thus, it remains to show that $\mathrm{red^{SG}}$ is indeed \emph{locally Cartesian} and \emph{coCartesian} in the sense of \cite[Definition $5.8$]{Steinebrunner_2022}. 

We start by showing $\mathrm{red^{SG}}$ is \emph{locally Cartesian}. Let $(A, \tau)$ be a morphism from $A_0$ to $A_{\tau}$ such that $A$ is equivalent to a reduced surface. It can equivalently be seen as a morphism in $\cobred$ or $\sCobtwo$ since $A$ has no spherical component. According to Proposition \ref{lem-prod-clos-red}, the reduction map $\mathrm{red}^{\mathrm{SG}}(P,P') : \sCobtwo(P,P') \rightarrow \cobred(P,P')$ is equivalent to the projection map $\cobred(P,P') \times \mathrm{Sph^{SG}} \rightarrow \cobred(P,P')$, hence the homotopy fiber of $\mathrm{red}(P,P')$ is equivalent to $\mathrm{Sph^{SG}}$.

    Let $(A_{0}\times \bR, 1)$ be an equivalence of $A_0$. We now consider the following diagram:
\[\begin{tikzcd}
	{\sCobtwo(A_0,A_0)} & {} & {\sCobtwo(A_0,A_{\tau})} \\
	{\cobred(A_0,A_0)} && {\cobred(A_0,A_{\tau})}
	\arrow["{(A,\tau)\circ -}"', from=1-1, to=1-3]
	\arrow["{\mathrm{red}^{\mathrm{SG}}(A_0,A_0)}"', from=1-1, to=2-1]
	\arrow["{\mathrm{red}^{\mathrm{SG}}(A_0,A_{\tau})}", from=1-3, to=2-3]
	\arrow["{(A,\tau)\circ -}"', from=2-1, to=2-3]
\end{tikzcd}.\] Since $A$ is equivalent to a reduced surface, the postcomposition map $(A,\tau)\circ -$ induces an equivalence $$(A,\tau)\circ - : \mathrm{hofib}(\mathrm{red^{SG}}(A_0,A_0))_{(A_0\times \bR,1)} \rightarrow \mathrm{hofib}(\mathrm{red^{SG}}(A_0,A_{\tau}))_{(A,\tau)\circ (A_0\times \bR,1)}.$$ This shows that $(A,\tau)$ is locally $\mathrm{red^{\mathrm{SG}}}$-Cartesian, as in \cite[Definition $5.8$]{Steinebrunner_2022}. We can find a reduced morphism $(A,\tau)$ between any pair of objects of $\cobred$, hence we can conclude that the functor $\mathrm{red^{SG}}$ is locally Cartesian, as in \cite[Definition $5.8$]{Steinebrunner_2022}. The opposite functor $(\mathrm{red}^{\mathrm{SG}})^{\mathrm{op}} : (\sCobtwo)^{\mathrm{op}} \rightarrow (\cobred)^{\mathrm{op}}$ is again locally Cartesian since it only reverses the cobordisms, hence $\mathrm{red}^{\mathrm{SG}}$ is also locally coCartesian.

Finally, the smooth case follows by adapting the steps above or adapting the proof, without too much work, of \cite[Theorem B]{Steinebrunner_2020}. 

\end{proof}    
    
We finish this subsection with the following remarks:

\begin{rk} \label{remark-d-reduced}
    A priori, we could write fiber sequences similar to the one in Propositon \ref{thm-fiber-reduction} in arbitrary dimension $d$. Moreover, instead of only deleting spherical components, we could choose to delete all endomorphisms of the empty object in $\PCobdG(X)$. Let $\mathrm{Cob}_d^{\mathrm{(S)G,red}}(X)$ denote the cobordism category obtained from $\PCobdG(X)$ by deleting all connected components in the morphisms which are equivalent to $d$-dimensional Poincaré complexes. There is again a reduction functor $\PCobdG(X) \rightarrow \mathrm{Cob}_d^{\mathrm{\cSG,red}}(X)$.
    Then, we could prove in the same vein as Proposition \ref{thm-fiber-reduction}, that there is a fiber sequence $$\B (\PCobdG(X)(\emptyset,\emptyset)) \rightarrow \B \PCobdG(X) \rightarrow \B \mathrm{Cob}_d^{\mathrm{\cSG,red}}(X).$$

    However, there is no hope in generalizing the formula from Theorem \ref{thmB} to higher dimensions. Indeed, in dimension $2$, we were able to compare the reduced cobordism categories $\mathrm{Cob}_2^{\mathrm{SO,red}}(X)$ and $\cobred(X)$. We even showed they are equivalent, which is most probably far from being the case in higher dimension.
\end{rk}

\begin{rk} \label{segalspace} 
    We showed in \cref{cat-fibrant} that the category $\sCobtwo(X)$ is fibrant, in other words the map $$\mathrm{Mor}(\sCobtwo(X)) \rightarrow \mathrm{Ob}(\sCobtwo(X))^2$$ is a Kan fibration. We deduce that the semi-simplicial space $\mathrm{N}_{\bullet} |\sCobtwo(X)|$, i.e. the levelwise geometric realization of the nerve $\mathrm{N}_{\bullet} \sCobtwo(X)$, is a semi-Segal space. 
    
    Moreover, equivalences in $\sCobtwo(X)$ are exactly the weak units, hence the nerve $\mathrm{N}_{\bullet} |\sCobtwo(X)|$ is a complete semi-Segal space.
\end{rk}

\begin{rk} 
    It seems like the fiber sequence from Proposition \ref{thm-fiber-reduction} could be deduced from the more general result in \cite[Observation $5.19$]{barkan2024equifiberedapproachinftyproperads}. Let $\mathcal{P}$ be an infinity-properad. Roughly speaking, it is a symmetric monoidal $\infty$-category, such that its spaces of objects and morphisms are freely generated by a suspace of \emph{connected objects} and \emph{connected morphisms}. Let $\mathcal{P}_0$ be the space of endomorphisms of the unit of $\mathcal{P}$. Let $\bar{\mathcal{P}}$ denote the cofiber $\mathcal{P}/\mathcal{P}_0$, in $\mathrm{Cat}_{\infty}^{\otimes}$, the $\infty$-category of symmetric monoidal $\infty$-categories. Morally, $\bar{\mathcal{P}}$ is obtained from $\mathcal{P}$ by removing all endomorphisms of the unit. It is shown in \cite[Observation $5.19$]{barkan2024equifiberedapproachinftyproperads} that there is a fiber sequence on geometric realization $$|\mathcal{P}_0| \rightarrow |\mathcal{P}| \rightarrow |\bar{\mathcal{P}}|.$$ Moreover, they show that $\Omega |\mathcal{P}_0|$ is equivalent to $Q_+ \mathcal{P}(\emptyset,\emptyset)$. Roughly speaking, $\mathcal{P}(\emptyset,\emptyset)$ is the space of \emph{connected generators} of $\mathcal{P}_0$.

    Using the $\Gamma$-structure from Subsection \ref{4.1} and Remark \ref{segalspace}, we can say that $\PCobdG$ and $\PCobdO$ are symmetric monoidal $\infty$-category. In both categories, all objects are obtained by taking finite disjoint unions on the subspace of objects which are connected, and similarly for morphisms. We could then argue that this makes $\PCobdG, \PCobdO$ into $\infty$-properads. Taking the reduction $\overline{\PCobdG}$ would correspond to taking the reduction category $\mathrm{Cob}_d^{\mathrm{\cSG,red}}$ described in Remark \ref{remark-d-reduced}. 
\end{rk}

\section{Parametrized Pontryagin-Thom Construction for Fibrations} \label{4}

We wish in this section to generalize the parametrized Pontryagin-Thom construction for manifold bundles from \cite{GMTW} to maps with homotopy fiber a Poincaré duality space, without relying on pointset models. More precisely, let $P$ be a Poincaré complex and let $p : E\rightarrow B $ be a map with homotopy fiber $P$. In \cref{dualizingcomplexfamilies}, we define a map which assembles the family of Spivak fibrations of the fibers of $p$ into a spherical fibration $\mathrm{D}_p$ over the total space $E$. Then, we show we can put together the family of the Pontryagin-Thom maps of the fibers into a single map $$ \mathrm{PT}_p: \susplus B \rightarrow \mathrm{Th}(\mathrm{D}_p).$$ 

On the other hand, the map of \cite{GMTW} is defined for smooth bundles $\pi : E \rightarrow B$ with fiber a smooth manifold $M$. If $T\pi$ denotes the vertical tangent bundle of $\pi$, i.e. the family of the tangent bundles of the fibers, then taking the family of Pontryagin-Thom maps of the fibers gives a map $$ \susplus B \rightarrow \Th(-T\pi),$$ where $-T\pi$ is the stable inverse or equivalently the family of the stable normal vector bundles of the fibers. In particular, we see that the map of \cite{GMTW} coincides with our map $\mathrm{PT}_{\pi}$, where we forget the bundle structure on $\pi$. 

The construction of \cite{GMTW} relies on pointset models of $\B \Diff(M)$. On the other hand, our approach uses the yoga of six functors on parametrized spectra. In particular, we claim no originality to the material presented. We rather adapt the language of relative dualizing objects and twisted norm maps, as defined in \cite{cnossen2023twistedambidexterityequivarianthomotopy}, to families of Poincaré complexes.

\subsection{Dualizing objects and Spherical Fibrations} \label{4.1}

In this subsection we work with the infinity-category of spaces $\mathcal{S}$. Most of the material presented below is explained in Appendix A of \cite{land2021reducibilitylowdimensionalpoincare}.
For $X$ a space, let $\mathrm{Sp}^X$ denote the category $\mathrm{Fun}(X^{\mathrm{op}}, \mathrm{Sp})$ of parametrized spectra over $X$. For any map $f: X \rightarrow Y$, the pullback functor $f^* : \mathrm{Sp}^Y \rightarrow \mathrm{Sp}^X$ admits a left adjoint $f_!$ and a right adjoint $f_*$, both obtained by taking respectively the left and right Kan extension. The category $\mathrm{Sp}^X$ is symmetric monoidal where we denote the product by $\otimes$. It admits an internal hom, which we denote by $\mathrm{hom}_X$, adjoint to the tensor product. In others words, $\mathrm{hom}_X(\mathcal{F},-)$ is right adjoint to $(\mathcal{F} \otimes -)$ for any object $\mathcal{F}$. The proposition below gives an important relation between the functors $f_!$ and $f^*$:

\begin{prop}
    Let $f: X \rightarrow Y$. The functor $f^*$ is symmetric monoidal. Moreover the left adjoint $f_!$ satisfies the following formula, called the projection formula, for any $\mathcal{F}, \mathcal{G} \in \mathrm{Sp}^X, \mathrm{Sp}^Y$:
    $$f_!(\mathcal{F}) \otimes \mathcal{G} \simeq f_!(\mathcal{F} \otimes f^*(\mathcal{G})).$$
\end{prop}

Let $r : X \rightarrow \star$ denote the unique map to the point. Then the functor $r_! : \mathrm{Sp}^X \rightarrow \mathrm{Sp}$ corresponds to taking the colimit over $X$ while the functor $r_* : \mathrm{Sp}^X \rightarrow \mathrm{Sp}$ corresponds to taking the limit over $X$.

We can now reformulate spherical fibrations over $X$ as well as their Thom spectra in terms of parametrized spectra:

\begin{definition}
    Let $\mathcal{F}$ be an object of $\mathrm{Sp}^X$. We say $\xi \in \mathrm{Sp}^X$ is a \emph{spherical fibration over $X$} if it is pointwise invertible or equivalently an element of $\mathrm{Fun}(X^{\mathrm{op}}, \mathrm{Pic}(\bS))$. Let $\bS_X$ denote the trivial spherical fibration over $X$.

    If $X$ is connected and $\xi$ is spherical fibration over $X$, let $\mathrm{rk}(\xi)$ be the integer defined by post-composing the functor $\xi  : X^{\mathrm{op}} \rightarrow \mathrm{Pic}(\bS)$ with the rank map $\mathrm{Pic}(\bS) \rightarrow \bZ$.

    If $r : X \rightarrow \star$ and $\xi$ is a spherical fibration, then we say $r_!(\xi)$ is the Thom spectrum of $\xi$. 
\end{definition}
The following proposition gives relations between the pullbacks and 
pushforwards functors induced by a pullback square: 

\begin{prop}
    Let $$\xymatrix{E' \ar[r]^{g} \ar[d]^{q} & E \ar[d]^p \\ B' \ar[r]^f & B}$$ be a pullback square. Then the functors satisfy the Beck-Chevalley isomorphisms:
    \begin{equation} \label{Beck-Chev}
\begin{split}
q_!g^* & \simeq f^*p_!\\
 q_*g^* & \simeq f^*p_*\\
 g_!q^* & \simeq p^*f_!\\
 g_*q^* &\simeq p^*f_*. 
\end{split}
\end{equation}
\end{prop}

Let $p : E\rightarrow B$ be a map. In the following definition, we associate to $p$ a parametrized spectrum $\mathrm{D}_p$ over $E$, which is used afterwards to compare the pushforward functor $p_*$ and the lower shriek functor $p_!$: 

\begin{definition}[Definition $3.5$ in \cite{cnossen2023twistedambidexterityequivarianthomotopy}]
    Let $p : E \rightarrow B$ be a map and consider the pullback square 
\begin{equation} \label{square-dualizing-map} \begin{tikzcd}
	{E\times_B E} & E \\
	E & B
	\arrow["{\pi_1}"', from=1-1, to=1-2]
	\arrow["{\pi_2}"', from=1-1, to=2-1]
	\arrow["p", from=1-2, to=2-2]
	\arrow["p", from=2-1, to=2-2]
\end{tikzcd}.
\end{equation}
    The dualizing complex of $p$ is defined as $\pi_{1_*} \Delta_!(\bS_E)$ where $\Delta : E \rightarrow E\times_B E$ is the diagonal map.
 \end{definition}

Playing with the Beck-Chevalley isomorphisms induced by the square \eqref{square-dualizing-map}, the counit of the adjunction $(\pi_1^*,\pi_{1_*})$,  as well as the projection formula gives the following natural transformation: 
\begin{equation} \label{Beck-Chev-appli}
\begin{split}
p^*p_!(\mathrm{D}_p \otimes -) & \simeq \pi_{2_!} \pi_1^* (- \otimes \pi_{1_*}\Delta_! \bS_E) \\
 & \simeq \pi_{2_!} (\pi_{1}^*(-) \otimes \pi_1^* \pi_{1_*} \Delta_!(\bS_E))\\
 & \longrightarrow \pi_{2_!} (\pi_{1}^*(-) \otimes \Delta_!(\bS_E))\\
 &\simeq \pi_{2_!} \Delta_! (\Delta^* \pi_{1}^*(-) \otimes \bS_E)\\
 & \simeq \text{id}(-). 
\end{split}
\end{equation} 

The adjoint of the natural transformation \eqref{Beck-Chev-appli} gives a natural transformation $$p_!(\mathrm{D}_p \otimes -) \rightarrow p_*(-),$$ called the \emph{twisted norm map} in \cite{cnossen2023twistedambidexterityequivarianthomotopy}. The following proposition from \cite{cnossen2023twistedambidexterityequivarianthomotopy} says when exactly this natural transformation is an equivalence:

\begin{prop}[Corollary $3.14$ in \cite{cnossen2023twistedambidexterityequivarianthomotopy}] \label{prop-5.1.5}
    Let $p : E \rightarrow B$ be a map such that the fiber is a compact space (in the $\infty$-categorical sense, e.g. a finitely dominated space). Then the natural transformation \eqref{Beck-Chev-appli} $$p_!(\mathrm{D}_p \otimes -) \rightarrow p_*(-)$$ is an equivalence.
\end{prop}

We now apply the previous notions to the terminal map $r : X \rightarrow \star$.

\begin{definition}\cite[Definition A.4]{land2021reducibilitylowdimensionalpoincare} \label{dualizing-spec-PD}
    Let $X$ be a compact space. Let $\mathrm{D}_X$ denote the dualizing complex $\mathrm{D}_r$ with $r : X \rightarrow \star$. According to Proposition \ref{prop-5.1.5}, the natural transformation \eqref{Beck-Chev-appli} $$r_!(\mathrm{D}_X \otimes -) \rightarrow r_*(-)$$ is an equivalence. We say that $\mathrm{D}_X$ is the \emph{dualizing complex} of $X$. Combining the unit of the adjunction $(r^*,r_*)$ with the equivalence \eqref{Beck-Chev-appli} gives a canonical map: $$\bS \rightarrow r_!(\mathrm{D}_X).$$ We call it the \emph{Pontryagin-Thom map} of $X$.
\end{definition} 

We now give an equivalent characterization of Poincaré complexes in terms of parametrized spectra as formulated in \cite{land2021reducibilitylowdimensionalpoincare} and originally due to \cite{Klein_2007}:

\begin{definition}[Definition $A.7$ in \cite{land2021reducibilitylowdimensionalpoincare}]\
  Let $X$ be a finite space. Then, $X$ is a Poincaré complex in the sense of \cref{definition-2.1.1} if and only if its dualizing complex $\mathrm{D}_X$ is a spherical fibration. Its dimension is defined as $-\mathrm{rk}(\mathrm{D}_X)$, where $\mathrm{rk}(\mathrm{D}_X)$ denotes the rank of the spherical fibration $\mathrm{D}_X$.
\end{definition}

We observe that if $X$ is a Poincaré complex, then $\mathrm{D}_X$ corresponds to the Spivak fibration from \cref{thm-2.1.3} and the Pontryagin-Thom map $\bS \rightarrow r_!(\mathrm{D}_X)$ is the same as in \cref{thm-2.1.3}. Moreover, according to \cite[Lemma A.6]{land2021reducibilitylowdimensionalpoincare}, if $X$ is compact, then $r_!(\mathrm{D}_X)$ is equivalent to the Spanier-Whitehead dual $D(X_+)$. While the Pontryagin-Thom map $\bS \rightarrow r_!(\mathrm{D}_X)$ is equivalent to the map $D(\star_+) \rightarrow D(X_+)$ induced by $r: X\rightarrow \star$.

The following proposition asserts that taking the dualizing complex actually preserves pullbacks. As we have not found a proof of this fact in the literature, we give one below: 

\begin{prop} \label{dua-pull}
    Let  $$\xymatrix{E' \ar[r]^{g} \ar[d]^{q} & E \ar[d]^p \\ B' \ar[r]^f & B}$$ be a pullback square. Then there is an equivalence $$\mathrm{D}_q \simeq g^* \mathrm{D}_p.$$
\end{prop}

\begin{proof}
    Each face of the following cube is Cartesian: 
\[\begin{tikzcd}
	{E'\times_{B'}E'} && {E\times_B E} \\
	& {E'} && E \\
	{E'} && E \\
	& {B'} && B
	\arrow["H"{description, pos=0.3}, from=1-1, to=1-3]
	\arrow["{\pi'_1}"{description}, from=1-1, to=2-2]
	\arrow["{\pi'_1}"{description, pos=0.7}, from=1-1, to=3-1]
	\arrow["{\pi_1}"{description}, from=1-3, to=2-4]
	\arrow["{\pi_1}"{description, pos=0.7}, dashed, from=1-3, to=3-3]
	\arrow["g"{description, pos=0.3}, from=2-2, to=2-4]
	\arrow["q"{description, pos=0.6}, from=2-2, to=4-2]
	\arrow["p"{description, pos=0.6}, from=2-4, to=4-4]
	\arrow["g"{description, pos=0.3}, dashed, from=3-1, to=3-3]
	\arrow["q"{description}, from=3-1, to=4-2]
	\arrow["p"{description}, from=3-3, to=4-4]
	\arrow["f"{description}, from=4-2, to=4-4]
\end{tikzcd}.\] On the other hand, consider the following commutative diagram:

\begin{equation} \label{bigsquarediagonal}
\begin{tikzcd}
	{E'} && {E'\times_{B'}E'} && {E'} \\
	E && {E\times_B E} && E
	\arrow["{\Delta'}"{description}, from=1-1, to=1-3]
	\arrow["g"{description}, from=1-1, to=2-1]
	\arrow["{\pi'_1}"{description}, from=1-3, to=1-5]
	\arrow["H"{description}, from=1-3, to=2-3]
	\arrow["g"{description}, from=1-5, to=2-5]
	\arrow["\Delta"{description}, from=2-1, to=2-3]
	\arrow["{\pi_1}"{description}, from=2-3, to=2-5]
\end{tikzcd}.
\end{equation} The composite top and bottom maps are respectively equivalent to the identity on $E'$ and $E$. Consequently, the outer square is Cartesian. Since the right-hand square is also Cartesian, according to the pullback pasting law, the left square  $$\xymatrix{ E' \ar[d]^{\Delta'} \ar[r]^g & E \ar[d]^{\Delta} \\ E'\times_{B'} E' \ar[r]^{H}  & E\times_{B} E}$$ is once again a pullback. We can now conclude by applying the Beck-Chevalley isomorphisms \eqref{Beck-Chev} to both the left and right squares in the diagram \eqref{bigsquarediagonal}:
\begin{equation*}
\begin{split}
\mathrm{D}_q & \simeq \pi'_{1_*} \Delta'_! (\bS_{E'}) \\
 & \simeq \pi'_{1_*} \Delta'_! g^*(\bS_{E})\\
 & \simeq \pi'_{1_*} \Ho^* \Delta_!(\bS_E)\\
 &\simeq g^* \pi_{1_*} \Delta_!(\bS_E)\\
 & \simeq g^*\mathrm{D}_p. 
\end{split}
\end{equation*}

\end{proof}

As a corollary let $X$ be the fiber of a map $p : E \rightarrow B$ and assume it is compact. We can write the following pullback square: 
\begin{equation} \label{square-param-PT-5.1.9}\begin{tikzcd}
	X & E \\
	{\star} & B
	\arrow["j"', from=1-1, to=1-2]
	\arrow["r"', from=1-1, to=2-1]
	\arrow["p"', from=1-2, to=2-2]
	\arrow["{\mathrm{ev}}"', from=2-1, to=2-2]
\end{tikzcd},
\end{equation} where $\mathrm{ev} : \star \rightarrow B$ denotes taking a point in $B$. The pullback $j^*\mathrm{D}_p$ of $\mathrm{D}_p$ along the inclusion of the fiber is then equivalent to $\mathrm{D}_X$. The dualizing object $\mathrm{D}_p$ can then be seen as a family over $E$ of the dualizing complexes of the fibers.

On the other hand, the unit of the adjunction $(p^*,p_*)$ gives a map $\bS_B \rightarrow p_*p^* (\mathbb{S}_B) \simeq p_!(\mathrm{D}_p)$. Applying $\mathrm{ev}^*$ gives a map $\mathrm{ev}^*(\bS_B) \simeq  \bS \rightarrow \mathrm{ev}^*p_!(\mathrm{D}_p)$. Beck-Chevalley isomorphisms applied to the square \eqref{square-param-PT-5.1.9} identify $\mathrm{ev}^*p_!(\mathrm{D}_p)$ with $r_!j^*(\mathrm{D}_p)$, which is equivalent to $r_!(\mathrm{D}_X)$. Consequently, the map $\mathrm{ev}^*(\bS_B) \simeq  \bS \rightarrow \mathrm{ev}^*p_!(\mathrm{D}_p)$ is equivalent to the Pontryagin-Thom collapse map $\bS \rightarrow r_!(\mathrm{D}_X)$ of the fiber $X$. Broadly speaking, the unit map $\bS_B \rightarrow p_!(\mathrm{D}_p)$ is the family of the Pontryagin-Thom maps of the fibers.

We can now defined the parametrized Pontryagin-Thom map announced in the introduction of this \cref{4}.

\begin{definition} \label{dualizingcomplexfamilies}
    Let $p : E \rightarrow B$ be a map whose fiber is a Poincaré complex $X$. Then the dualizing complex $\mathrm{D}_p$ of $p$ is a spherical fibration over $E$ such that $j^*\mathrm{D}_p$ is equivalent to $\mathrm{D}_X$ where $j$ is as in the square \eqref{square-param-PT-5.1.9}. 

    Let $r : B \rightarrow \star$. Applying $r_!$ to the unit map $\bS_B \rightarrow p_!(\mathrm{D}_p)$ gives a map: $$ \mathrm{PT}_p : \susplus B \rightarrow r_!(\mathrm{D}_p).$$ We say it is the \emph{parametrized Pontryagin-Thom map} of $p$. 
\end{definition} 

Let $f: X \rightarrow Y$ be a space, and let $\xi  \in \mathrm{Sp}^Y$ be a spherical fibration over $Y$. Let $r^X, r^Y$ be the terminal maps. Using the relation $r^X= r^Y \circ f$ and the counit map of the adjunction $(f_!,f^*)$, we get a map of Thom spectra $r^X_!(f^*\xi) \rightarrow r^Y_!(\xi)$. We now show that the parametrized Pontryagin-Thom construction is compatible with pullback:

\begin{cor} \label{comm-PT}
 Let 
\[\begin{tikzcd}
	{E'} & E \\
	{B'} & B
	\arrow["g", from=1-1, to=1-2]
	\arrow["q"', from=1-1, to=2-1]
	\arrow["p", from=1-2, to=2-2]
	\arrow["f"', from=2-1, to=2-2]
\end{tikzcd}\] be a pullback square. Assume the fiber of $p$ is a Poincaré complex. Then, the following diagram is \emph{commutative}:

\[\begin{tikzcd}
	{\susplus B'} & {(r^{B'})_!(D_q)} \\
	{\susplus B} & {(r^B)_!(D_p)}
	\arrow["{\mathrm{PT}_q}", from=1-1, to=1-2]
	\arrow["f"', from=1-1, to=2-1]
	\arrow[from=1-2, to=2-2]
	\arrow["{\mathrm{PT}_p}"', from=2-1, to=2-2]
\end{tikzcd}.\]
\end{cor}

\begin{proof}
    The unit of the adjunction $(p^*,p_*)$ gives a map $\bS_B \rightarrow p_*p^*(\bS^B)$. The counit of the adjunction $(f_!,f^*)$ gives a natural transformation $f_!f^* \rightarrow \mathrm{id}$. Consequently, the following diagram commutes: 
\[\begin{tikzcd}
	{f_!f^*(\bS_B)} & {\bS_B} \\
	{f_!f^*(p_*p^*(\bS_B))} & {p_*p^*(\bS_B)}
	\arrow[from=1-1, to=1-2]
	\arrow[from=1-1, to=2-1]
	\arrow[from=1-2, to=2-2]
	\arrow[from=2-1, to=2-2]
\end{tikzcd}.\] We now apply $(r^B)_!$ to it: 
\[\begin{tikzcd}
	{(r^B)_!(f_!f^*(\bS_B))} & {(r^B)_!(\bS_B)} \\
	{(r^B)_!(f_!f^*(p_*p^*(\bS_B)))} & {(r^B)_!(p_*p^*(\bS_B))}
	\arrow[from=1-1, to=1-2]
	\arrow[from=1-1, to=2-1]
	\arrow[from=1-2, to=2-2]
	\arrow[from=2-1, to=2-2]
\end{tikzcd}.\] By construction, we identify the right-hand vertical map with $\mathrm{PT}_p$. The Beck-Chevalley isomorphisms as well as \cref{dua-pull} give equivalences $f^*p_!(\mathrm{D}_p) \simeq q_!g^*(\mathrm{D}_p) \simeq q_!(\mathrm{D}_q)$. Since $r^B \circ f \circ q =r^{E'}$, we conclude the left-hand side vertical map is equivalent to $\mathrm{PT}_q$. 
\end{proof}

As a corollary, we also recover the following result due to Gottlieb (see \cite{442c80bf-efe8-3420-b1f6-cc10866773ac}) on Poincaré complexes: 

\begin{cor}
      Let $p : E \rightarrow B$ be a map with fiber $P$ a Poincaré duality space of dimension $d$. If $B$ is a Poincaré duality space of dimension $m$ and $E$ is finite, then $E$ is also a Poincaré duality space of dimension $m+p$ and its dualizing complex satisfies the relation $$\mathrm{D}_E \simeq \mathrm{D}_p \otimes p^*\mathrm{D}_B.$$
\end{cor}

\begin{proof}
     Let $r^E : E \rightarrow *$, $r^B : B \rightarrow *$. 
    Then we have the following chain of natural equivalences: 
    \begin{equation} 
\begin{split}
r^E_*(-) & \simeq r^B_* p_*(-) \\
 & \simeq r^B_! (\mathrm{D}_B \otimes p_*(-))\\
 & \simeq r^B_! (\mathrm{D}_B \otimes p_!(\mathrm{D}_p \otimes -))\\
 &\simeq r^B_! p_! (p^*\mathrm{D}_B \otimes \mathrm{D}_p \otimes -)\\
 & \simeq r ^E_!(p^*\mathrm{D}_B \otimes \mathrm{D}_p \otimes -). 
\end{split}
\end{equation}
This identifies the dualizing complex $\mathrm{D}_E$ of $E$ with $p^*\mathrm{D}_B\otimes \mathrm{D}_p$. Since $B$ and $P$ are Poincaré complexes, both $\mathrm{D}_B$ and $\mathrm{D}_p$ are spherical fibrations, hence $\mathrm{D}_E$ is a spherical fibration.
\end{proof}

We now informally explain why the parametrized Pontryagin-Thom map from \cite{GMTW} and our map from \cref{dualizingcomplexfamilies} are equivalent. Let $\pi : E \rightarrow B$ be a smooth bundle with fiber $M$. Let $\nu \pi$ denote the stable inverse of the bundle vertical tangent bundle $T\pi$ over $E$. The bundle $\nu \pi$ is then the family of the stable normal bundles of the fibers of $\pi$. The underlying spherical fibration of $\nu \pi$ is the family of the underlying spherical fibrations of the normal bundles of the fibers of $\pi$, hence it is the family of the dualizing complexes $\mathrm{D}_M$ of the fibers. Then, we can deduce the underlying spherical fibration of $\nu \pi$ is equivalent to the spherical fibration $\mathrm{D}_{\pi}$ over $E$. We now briefly sketch the construction of the parametrized Pontryagin-Thom map $\susplus B \rightarrow \Th(\nu \pi)$. Up to filtering $B$ by finite CW-complexes, we pick an embedding $e :E \hookrightarrow B \times \bR^N$. The bundle $\nu \pi$ is then equivalent to the stabilization of the normal bundle $\nu_e$ of $e$. We then take the collapse map $(B\times \bR^N)^+ \rightarrow \nu E^+$, where $+$ denotes the one-point compactification. We then identify $(B\times \bR^N)^+$ with $\Sigma^N_+B$ and the one point compactification of a tubular neighborhood $\nu E$ of $e$ with $\Th(\nu e)$. On each fiber, we get a collapse map $S^N \rightarrow \Th(i^*\nu e)$, where $i : M \hookrightarrow E$ denotes the inclusion of the fiber, and $i^* \nu e$ is a model of the normal bundle of $M$.

This discussion is summarized in the following corollary: 

\begin{cor}
    Let $M$ be a smooth $d$-dimensional manifold and let $\pi : E \rightarrow B$ be a bundle with fiber $M$. Let $T\pi$ denote the fiberwise tangent bundle of $M$ over $E$ and let $\nu \pi$ denote its stable inverse $-T\pi$. Then the underlying spherical fibration of $\nu \pi$ is equivalent to the spherical fibration $D\pi$. Moreover, the parametrized Pontryagin-Thom map or scanning map from \cite{GMTW} is equivalent to $\mathrm{PT}_{\pi} : \susplus B \rightarrow r_!(\mathrm{D}\pi)$.
\end{cor}

We conclude with some remark on notations. In section \ref{5}, we use the parametrized Pontryagin-Thom construction for maps $P\sslash\mathrm{H} \rightarrow \mathrm{BH}$, where $\mathrm{H}$ is a submonoid of $\haut(P)$.

\begin{notation} \label{notation-4.1.12}
    Let $P$ be a Poincaré duality space and take a monoid map $\mathrm{H} \rightarrow \haut(P)$. We consider the map $p : P\sslash\text{H} \rightarrow \mathrm{BH}$ with fiber $P$. We denote by $\mathrm{D}_P^H$ the parametrized spectrum $\mathrm{D}_p$. We may also write $\nu_P^H$ for the spherical fibration and $\mathrm{PT}_P^H$ instead of $\mathrm{PT}_p$.
\end{notation}

\definecolor{mygreen}{RGB}{0,100,0}

\section{Goodwillie Calculus and Quotients of Stable Mapping Spaces} \label{5}

This section aims to give the necessary tools and results from Goodwillie calculus for the proof of \cref{thmA} in Subsection \ref{6.2}. Let $\mathcal{S}$ be the $\infty$-category of spaces. 
In what follows, by category, limit, and colimit we actually mean infinity-category, homotopy limit and homotopy colimit.

In Section \ref{3}, \ref{thmB}, we constructed a functor $\mathrm{PH}(2,-) : \mathcal{S} \rightarrow \mathrm{Sp}$ and showed it is a delooping of $\B \sCobtwo(-)$. In subsection \ref{5.2}, we see that post-composing with $\lop$ commutes with most operations we describe below. Hence, it is reasonable to restrict to functors from $\mathcal{S}$ to $\mathrm{Sp}$. 

A functor $E: \mathcal{S} \rightarrow \mathrm{Sp}$ is excisive if it preserves pushouts. Using a result of Lurie (\cref{Lurie-6151}), we give in Subsection \ref{5.1}, \cref{class_excisive}, a general decomposition formula for an excisive functor $E: \mathcal{S} \rightarrow \mathrm{Sp}$. More precisely, we show every excisive functor $E$ preserving filtered colimits is equivalent to the pullback of a diagram 
\[\begin{tikzcd}
	& {A\otimes (\susplus -)} \\
	B & A
	\arrow[from=1-2, to=2-2]
	\arrow["f"{description}, from=2-1, to=2-2]
\end{tikzcd},\] where the map $A\otimes (\susplus -) \rightarrow\mathrm{A}$ is induced by the terminal map $r : - \rightarrow \star$ and $f :\mathrm{A} \rightarrow B$ is a map in $\mathrm{Sp}$.

Goodwillie showed in \cite{Goodwillie_2003} that any functor $\mathrm{F}: \mathcal{S} \rightarrow \mathrm{Sp}$ can be approximated by an excisive functor $\mathrm{P}_1F$ via a map $p_1 \mathrm{F} : \mathrm{F} \Rightarrow \mathrm{P}_1 F$. Let $F$ be a filtered-colimit preserving functor and let $\eta : \mathrm{F} \Rightarrow\mathrm{A}\otimes (\susplus -)$ be a natural transformation. Let $E$ be the functor such that the square in the following diagram is a pullback: 
\[\begin{tikzcd}[ampersand replacement=\&,cramped]
	\textcolor{mypurple}{F} \\
	\& \textcolor{mypurple}{E} \& {A\otimes (\susplus -)} \& {.} \\
	\& {F(\star)} \&\mathrm{A}
	\arrow[draw={mypurple}, Rightarrow, from=1-1, to=2-2]
	\arrow["\eta"{description}, curve={height=-12pt}, Rightarrow, from=1-1, to=2-3]
	\arrow["{(-\circ r)}"{description}, curve={height=6pt}, from=1-1, to=3-2]
	\arrow[from=2-2, to=2-3]
	\arrow[from=2-2, to=3-2]
	\arrow["\lrcorner"{anchor=center, pos=0.125}, draw=none, from=2-2, to=3-3]
	\arrow["{(-\circ r)}", from=2-3, to=3-3]
	\arrow["{\eta(\star)}"', from=3-2, to=3-3]
\end{tikzcd}\]
Then the map $\eta$ as well as the terminal map $F \rightarrow F(\star)$ induce a natural transformation $\mathrm{F} \Rightarrow\mathrm{E}$ represented by the purple map in the diagram. \cref{5.15} in Subsection \ref{5.2} gives a necessary condition on $\eta$ for the purple map to be equivalent to the first approximation map $p_1 \mathrm{F} : \mathrm{F} \Rightarrow \mathrm{P}_1 F$. 

Finally, in Subsection \ref{5.3}, we use the recipe from Subsections \ref{5.1} and \ref{5.2} to compute the first polynomial approximation of the functor $\FPH : X \mapsto \susplus \Map(P,X)\sslash H$, where $P$ is a Poincaré complex and $H$ is a grouplike submonoid of $\haut(P)$. We begin with constructing in \cref{constr-parpt} a natural transformation $$\PTPH(-) : \FPH(-) \Rightarrow r_!(\DPH) \otimes (\susplus -)$$ such that at the point, the map $\FPH(\star) \rightarrow r_!(\DPH)$ is equivalent to the parametrized Pontryagin-Thom map $\PTPH$ from \cref{4.1}. Let $\mathrm{E}_{\mathrm{P}}^{\mathrm{H}}$ be the pullback of the cospan 
\[\begin{tikzcd}[ampersand replacement=\&,cramped]
	\& {r_!(\DPH)\otimes (\susplus -)} \\
	{\FPH(\star)} \& {r_!(\DPH)} 
	\arrow[from=1-2, to=2-2]
	\arrow[from=2-1, to=2-2]
\end{tikzcd}.\]
Next, in \cref{5.15}, we show the first approximation map $\FPH \Rightarrow \mathrm{P}_1 \FPH$ is equivalent to the natural transformation $\FPH \Rightarrow \mathrm{E}_{\mathrm{P}}^{\mathrm{H}}$ induced by $\eta$ and the terminal map $\FPH \rightarrow \FPH(\star)$, shown in purple in the following diagram:

\[\begin{tikzcd}[ampersand replacement=\&,cramped]
	\textcolor{mypurple}{{\FPH}} \\
	\& \textcolor{mypurple}{\mathrm{E}_{\mathrm{P}}^{\mathrm{H}}} \& {r_!(\DPH)\otimes (\susplus -)} \& {.} \\
	\& {\FPH(\star)\simeq \Sigma^{\infty}_+ \mathrm{BH}} \& {r_!(\DPH)}
	\arrow[draw={mypurple}, Rightarrow, from=1-1, to=2-2]
	\arrow["{\PTPH}", curve={height=-12pt}, Rightarrow, from=1-1, to=2-3]
	\arrow["{(-\circ r)}"{description}, curve={height=12pt}, from=1-1, to=3-2]
	\arrow[from=2-2, to=2-3]
	\arrow[from=2-2, to=3-2]
	\arrow["{(-\circ r)}", from=2-3, to=3-3]
	\arrow["{\PTPH}", from=3-2, to=3-3]
\end{tikzcd}\]

\subsection{Excisive functors} \label{5.1}

In this subsection, we define excisive functors and give a partial classification of excisive functors $\mathcal{S} \rightarrow \mathrm{Sp}$, following \cite[Chapter 6]{HA} 

\begin{definition} 
Let  $\mathcal{D}$ denote either $\mathcal{S}, \mathcal{S}_*$ or $\mathrm{Sp}$. A functor $E : \mathcal{S}  \rightarrow \mathcal{D}$ is excisive if it takes a pushout square to a pullback square in $\mathcal{D}$.
\end{definition}

\begin{rk} 
Let $A$ be a spectrum. Then the functor $X \mapsto A\otimes \susplus X$ is excisive. Moreover, the homotopy groups $\pi_*(A\otimes \susplus X)$ define a generalized homology theory on spaces. Intuitively, the excisive condition is like a generalization of Mayer-Vietoris. 
\end{rk}

The categories $\mathcal{S},\mathcal{S}_*$ and $\mathrm{Sp}$ all admit a terminal object. We consider the following class of functors, called reduced:

\begin{definition}
   A functor $\mathrm{F} : \mathcal{S} \rightarrow \mathcal{D}$ is reduced if it preserves the terminal object. A functor  $E : \mathcal{S} \rightarrow \mathcal{D}$ is linear or homogeneous if it is excisive and reduced.
\end{definition}

We now state the main proposition we wish to show in this subsection: 

\begin{prop} \label{class_excisive}
    Let $E : \mathcal{S} \rightarrow \mathrm{Sp}$ be an excisive functor commuting with filtered colimits. There exists a map of spectra $ f : B \rightarrow A$ such that the functor $E$ is equivalent to the pullback functor of the diagram 
\[\begin{tikzcd}
	& {A\otimes (\susplus -)} \\
	B & A
	\arrow["{(-\circ r)}", from=1-2, to=2-2]
	\arrow["g"{description}, from=2-1, to=2-2]
\end{tikzcd}.\]
    
\end{prop}

Proposition \ref{class_excisive} is a corollary of the subsequent theorem due to Lurie in \cite{HA}. Let $\textit{Fin}$ denote the category of finite sets and $\textit{Fin}^{\leq n}$ denote the category of finite sets with cardinal less or equal than $n$. In particular $\textit{Fin}^{\leq 1}$ is the category with two objects $\emptyset$ and $\{1\}$. The only non-identity morphism is the inclusion. Any excisive functor $E :\mathcal{S} \rightarrow \mathrm{Sp}$ preserving filtered colimits is actually determined by its values at $\emptyset$ and $\star$ as follows:

\begin{thm}[Theorem $6.1.5.1$ in \cite{HA}] \label{Lurie-6151}
    Let $E : \mathcal{S} \rightarrow \mathrm{Sp}$ be a functor. The following conditions are equivalent: 
   \begin{enumerate}
       \item The functor $E$ is excisive and commutes with filtered colimits.
       \item The functor $E$ is a left Kan extension of $E_{|\mathrm{N}(\text{Fin}^{\leq 1})}$. 
   \end{enumerate}
\end{thm}

We now prove Proposition \ref{class_excisive}:

\begin{proof}[Proof of Proposition \ref{class_excisive}]
   Let $B$ be $E(\star)$ and let $f$ be the map $E(\emptyset \rightarrow \star)$. Then let $A$ be the pushout of the following diagram:

\[\begin{tikzcd}
	{E(\emptyset)} & \star \\
	{E(\star)} & A
	\arrow[from=1-1, to=1-2]
	\arrow["f", from=1-1, to=2-1]
	\arrow[from=1-2, to=2-2]
	\arrow["g"{description}, from=2-1, to=2-2]
\end{tikzcd}.\] Let $g$ be the map $E(\star) \rightarrow\mathrm{A}$ coming from the pushout.

Let $D$ be the functor given by the pullback 
\[\begin{tikzcd}
	D & {A\otimes (\susplus -)} \\
	{E(\star)} & A
	\arrow[from=1-1, to=1-2]
	\arrow[from=1-1, to=2-1]
	\arrow["{(-\circ r)}", from=1-2, to=2-2]
	\arrow["g"{description}, from=2-1, to=2-2]
\end{tikzcd}.\] Then $D$ is again excisive since pullbacks commute with pullbacks in $\mathrm{Sp}$. It still commutes with filtered colimits since filtered colimits commute with finite limits in $\mathrm{Sp}$. 

At the point, $D(\star)$ is equivalent to $E(\star)$. Evaluating at $\emptyset$ gives the pullback square 
\[\begin{tikzcd}
	{D(\emptyset)} & \star \\
	{E(\star)} & A
	\arrow[from=1-1, to=1-2]
	\arrow[from=1-1, to=2-1]
	\arrow[from=1-2, to=2-2]
	\arrow["g"{description}, from=2-1, to=2-2]
\end{tikzcd}.\] 
Then $D(\emptyset)$ is the fiber of the map $g :\mathrm{E}(\star) \rightarrow\mathrm{A}$, hence the map $D(\emptyset) \rightarrow D(\star)$ is equivalent to the map $f :\mathrm{E}(\emptyset) \rightarrow\mathrm{E}(\star)$.

Both functors $E$ and $D$ satisfies the assumptions of \cref{Lurie-6151} and their restrictions to $N(\textit{Fin}^{\leq 1})$ agree, hence according to \cref{Lurie-6151} both functors $E$ and $D$ agree.
\end{proof}

We now give two corollaries of Proposition \ref{class_excisive}. The first one characterizes natural equivalences between excisive functors: 

\begin{cor} \label{cor 4.8}
    Let $E :\mathcal{S} \rightarrow \mathrm{Sp}$ be an excisive functor preserving filtered colimits and let $A$ be a spectrum. Let $\eta :\mathrm{E} \Rightarrow\mathrm{A} \otimes (\susplus -)$ be a natural transformation such that the square $$\xymatrix{E(\emptyset) \ar[r]^{\eta(\emptyset)} \ar[d]_{E(\emptyset\rightarrow \star)} & \star \ar[d]\\\mathrm{E}(\star) \ar[r]^{\eta(\star)} &\mathrm{A}}$$ is coCartesian. Then the induced natural transformation $\alpha$ from $E$ to the pullback of the cospan 
\begin{equation}\label{square-pullback-section51}\begin{tikzcd}
	& {A\otimes (\susplus -)} \\
	{E(\star)} & A
	\arrow[from=1-2, to=2-2]
	\arrow["\eta(\star)"{description}, from=2-1, to=2-2]
\end{tikzcd}\end{equation} is an equivalence of functors.
    \end{cor}

\begin{proof} Let $D$ be the functor obtained as the pullback of the span \eqref{square-pullback-section51}. For an object $X$, the map $\alpha(X)$ is such that this diagram commutes: 

\[\begin{tikzcd}
	{E(X)} \\
	& {D(X)} & {A\otimes \susplus X} \\
	& {E(\star)} & A
	\arrow["{\alpha(X)}"{description}, from=1-1, to=2-2]
	\arrow["{\eta(X)}"{description}, curve={height=-12pt}, from=1-1, to=2-3]
	\arrow["{E(r)}"{description}, curve={height=12pt}, from=1-1, to=3-2]
	\arrow[from=2-2, to=2-3]
	\arrow[from=2-2, to=3-2]
	\arrow[from=2-3, to=3-3]
	\arrow["{\eta(\star)}"{description}, from=3-2, to=3-3]
\end{tikzcd}.\]
At the point, $A \otimes \susplus \star \rightarrow\mathrm{A}$ is an equivalence, hence $\alpha(\star)$ is an equivalence.

We evaluate the diagram at $\emptyset$:
$$\xymatrix{E(\emptyset) \ar@/^1pc/[rr]^{\eta(\emptyset)} \ar[r]_{\alpha(\emptyset)} \ar[d]_{E(\emptyset \rightarrow \star)} & D(\emptyset) \ar[r] \ar[d] & \star \ar[d] \\\mathrm{E}(\star) \ar@{=}[r] &\mathrm{E}(\star) \ar[r]_{\eta(\star)} &\mathrm{A}}.$$
By definition, the square on the right handside is coCartesian. By assumption, the total square is coCartesian. The left handside square is then coCartesian according to the homotopy pullback pasting law. 
Consequently $\alpha(\emptyset)$ is an equivalence and $D(\emptyset) \rightarrow\mathrm{E}(\star)$ is equivalent to $E(\emptyset \rightarrow \star)$.

The restriction of functors $\alpha_{|\text{Fin}^{\leq 1}} : E_{|\text{Fin}^{\leq 1}} \Rightarrow D_{|\text{Fin}^{\leq 1}}$ is then an equivalence. Since both functors $E$ and $D$ are left Kan extensions of their restrictions $E_{|\text{Fin}^{\leq 1}}, D_{|\text{Fin}^{\leq 1}}$, $\alpha:\mathrm{E} \Rightarrow D$ is an equivalence of functors as well according to \cref{Lurie-6151}.
\end{proof}

We now introduce two notations:

\begin{notation}
    Let $\Sigma_{\mathcal{S}} : \mathcal{S} \rightarrow \mathcal{S}$ denote the functor taking a space $X$ to the pushout: 
\[\begin{tikzcd}
	X & \star \\
	\star & {\Sigma_{\mathcal{S}}X}
	\arrow[from=1-1, to=1-2]
	\arrow[from=1-1, to=2-1]
	\arrow[from=1-2, to=2-2]
	\arrow[from=2-1, to=2-2]
\end{tikzcd}.\] We note that by definition $\Sigma_{\mathcal{S}} \emptyset \simeq S^0$. 

Secondly, if $X$ is unbased, let $\widetilde{\Sigma}^{\infty} X$ denote the homotopy fiber of the map $\susplus X \rightarrow \susplus \star \simeq \bS$ induced by $r: X\rightarrow \star$. We note that $\widetilde{\Sigma}^{\infty} \emptyset \simeq \Sigma^{-1} \bS$.
\end{notation}

We end this subsection with characterizing homogeneous functors:

\begin{cor} \label{cor-homogeneous-functors}
    Let $E : \mathcal{S} \rightarrow \mathrm{Sp}$ be a homogeneous functor commuting with filtered colimits. There exists a spectrum $A$ such that $E$ is equivalent to the functor $A\otimes (\widetilde{\Sigma}^{\infty} -)$.
\end{cor}

\begin{proof}
Since $E(\star) \simeq \star$, according to Proposition \ref{class_excisive}, there exists a spectrum $A$ such that $E$ is equivalent to the pullback of the cospan:

\[\begin{tikzcd}
	& {A\otimes (\susplus -)} \\
	\star & A
	\arrow[from=1-2, to=2-2]
	\arrow[from=2-1, to=2-2]
\end{tikzcd},\] hence the statement.
\end{proof}

\subsection{Polynomial approximation and the first derivative} \label{5.2}

Let $\mathcal{D}$ denote either $\mathcal{S}$ or $\mathrm{Sp}$. All these categories are complete, cocomplete and admit a terminal object. Moreover, finite limits and filtered colimits commute in $\mathcal{D}$. Let $\mathrm{Fun}(\mathcal{S}, \mathcal{D})$ denote the $\infty$-category of functors and $\mathrm{Exc}_1(\mathcal{S}, \mathcal{D})$ denote the full $\infty$-category of excisive functors.  The following theorem, originally due to Goodwillie in \cite{Goodwillie_2003}, states that every functor $F : \mathcal{S} \rightarrow \mathcal{D}$ can be approximated by an excisive functor: 
\begin{thm}[Theorem $6.1.1.10$ in \cite{HA}] \label{P_1^X_F}
    For an object $X$ of $\mathcal{S}$, the composite functor $$\mathrm{Exc}_1(\mathcal{S}_{/X},\mathcal{D}) \hookrightarrow \mathrm{Fun}(\mathcal{S}_{/X}, \mathcal{D}) \hookrightarrow \mathrm{Fun}(\mathcal{S}, \mathcal{D})$$ admits a left adjoint denoted by $$\mathrm{P}_1^X.$$ For a functor $F \in \mathrm{Fun}(\mathcal{S}, \mathcal{D})$, we say that $\mathrm{P}_1^X F$ is the first polynomial approximation at $X$ of $F$.
\end{thm}

The unit of the adjunction gives a natural transformation $ p_1^X\mathrm{F} : F \rightarrow \mathrm{P}_1^X F$. Moreover the functor $\mathrm{P}_1^X F$ satisfies the following universal property. Any natural transformation $\mathrm{F} \Rightarrow\mathrm{E}$ in $\mathrm{Fun}(\mathcal{S}_{/X}, \mathcal{D})$ where $E$ is an excisive functor factors as $\mathrm{F} \Rightarrow \mathrm{P}_1^X \mathrm{F} \Rightarrow\mathrm{E}$. This makes $\mathrm{P}_1^X F$ the best possible approximation of $F$ at the object $X$ by an excisive functor. At the terminal object $X$ of $\mathcal{S}_X$, we can show $P_1^X F(X)$ is actually equivalent to $F(X)$:

\begin{prop}[Prop $1.17$ in \cite{Goodwillie_2003}]
  Let $F \in \mathrm{Fun}(\mathcal{S},\mathcal{D})$ and $X$ be an object of $\mathcal{S}$. The linear approximation of $F$ at $X$ is the homogeneous functor $\mathrm{D}_1^X \mathrm{F} : \mathcal{S}_{/X} \rightarrow \mathcal{D}$ defined at each object $Y$ by: $$\mathrm{D}_1^X F(Y)=\mathrm{fib}(\mathrm{P}_1^X F(Y) \rightarrow F(X)).$$  
\end{prop}

Let $\mathrm{Exc}_1^{\mathrm{red}}(\mathcal{S}, \mathcal{D})$ denote the subcategory of $\mathrm{Exc}_1(\mathcal{S},\mathcal{D})$ of reduced excisive functors. 

\begin{notation}
 Since $\mathcal{S}$ admits a final object $\star$ and $\mathcal{S}_{/\star}$ is equivalent to $\mathcal{S}$, we write $\mathrm{P}_1 F$ instead of $\mathrm{P}_1^*F$, and similarly we can denote by $\mathrm{D}_1 F$ its linearization $\mathrm{D}_1^* F$.
\end{notation}

We now restrict to functors from spaces to spectra, since most functors we are interested factor through $\mathrm{Sp}$. In the Corollary \ref{cor-homogeneous-functors}, we characterized homogeneous functors $\mathcal{S} \rightarrow \mathrm{Sp}$, and showed they are entirely characterized by a spectrum:

\begin{definition}
     Let $F : \mathcal{S} \rightarrow \mathrm{Sp}$ be a functor preserving filtered colimits. In particular, $\mathrm{D}_1F$ is a homogeneous functor preserving filtered colimits. There exists a spectrum $\partial_1 F$, called the \emph{first derivative of $F$ at the point} such that $\mathrm{D}_1F$ is equivalent to the functor $\partial_1F \otimes (\widetilde{\Sigma}^{\infty} -)$.
\end{definition}

We now aim to show the following Proposition. Given a functor $F : \mathcal{S} \rightarrow \mathrm{Sp}$ and a natural transformation $\eta : F \Rightarrow A \otimes (\susplus -)$, it provides a necessary condition for the derivative $\partial_1 F$ to be equivalent to $A$:

\begin{prop} \label{5.15}
    Let $F : \mathcal{S} \rightarrow \mathrm{Sp}$ be a functor preserving filtered colimits such that $F(\emptyset)$ is contractible. Let $\eta : F \Rightarrow\mathrm{A} \otimes (\susplus -)$ be a natural transformation, where $A \in \mathrm{Sp}$. If the map $\eta(\star)  : F(\star) \rightarrow\mathrm{A}$ is equivalent to the inclusion in the filtered colimit \[F(\star) \rightarrow \colim_n \Sigma^{-n+1} \mathrm{fib}(F(\Sigma^n_{\mathcal{S}} \emptyset) \rightarrow F(\star)), \] then $A$ is equivalent to $\partial_1 F$ and $\mathrm{P}_1F$ is equivalent to the pullback of the cospan $$\xymatrix{ &\mathrm{A} \otimes (\susplus -) \ar[d]^{- \rightarrow \star} \\ F(\star) \ar[r]^{\eta(\star)} &\mathrm{A} }.$$ The natural transformation $\mathrm{p}_1 F : F \Rightarrow \mathrm{P}_1F$ is induced by $\eta$ and the terminal transformation $F \Rightarrow F(\star)$.
\end{prop}

Before showing \cref{5.15}, we recall below an explicit formula due to Goodwillie to compute the first polynomial approximation $\mathrm{P}_1F$. For readability assume $F$ is reduced. Let $F :\mathcal{S} \rightarrow \mathrm{Sp}$ be a functor. For any object $X$ there is a pushout square:
$$\xymatrix{X \ar[d] \ar[r] & \star \ar[d] \\ \star \ar[r] & \Sigma_{\mathcal{S}} X} $$ which gives a commutative square in $\mathcal{C}$ after applying $F$:

\[\begin{tikzcd}
	{F(X)} & \star \\
	\star & {F(\Sigma_{\mathcal{S}}X)}
	\arrow[from=1-1, to=1-2]
	\arrow[from=1-1, to=2-1]
	\arrow[from=1-2, to=2-2]
	\arrow[from=2-1, to=2-2]
\end{tikzcd}.\] 
By universal property of pushouts there is a map $$F(X) \rightarrow \Sigma^{-1} F(\Sigma_{\mathcal{S}} X).$$ We have the following result:

\begin{prop}[Theorem $1.8$ in \cite{Goodwillie_2003}]\label{good_approx} 
    Let $F : \mathcal{S} \rightarrow \mathrm{Sp}$ be a reduced functor. The inclusion of the $0$-th stage in the colimit induces a natural transformation 
    \[F \Rightarrow \colim_{n\rightarrow \infty} \Sigma^{-n} \circ F\circ \Sigma_{\mathcal{S}}\] 
    equivalent to the first polynomial approximation map $$p_1F : F \Rightarrow \mathrm{P}_1 F.$$
\end{prop}

We now prove Proposition \ref{5.15}:

\begin{proof}
    The natural transformation $\eta : F \rightarrow A\otimes (\susplus -)$ factors through a natural transformation $\mathrm{P}_1 F \rightarrow A \otimes (\susplus -)$. Since $F$ preserves filtered colimits, $\mathrm{P}_1 F$ does as well. As we are now in the situation of \cref{cor 4.8}, it suffices to show the square 
\begin{equation} \label{square-A-partial}
    \begin{tikzcd}
	{\mathrm{P}_1 F(\emptyset)} & \star \\
	{\mathrm{P}_1F(\star)\simeq F(\star)} & A
	\arrow[from=1-1, to=1-2]
	\arrow[from=1-1, to=2-1]
	\arrow[from=1-2, to=2-2]
	\arrow["{\eta(\star)}"', from=2-1, to=2-2]
\end{tikzcd}
\end{equation} is coCartesian.

    Let $\tilde{F}$ denote the reduced functor $\mathrm{fib}(F(-) \rightarrow F(\star))$. The homogeneous approximation $\mathrm{D}_1 F$ is then equivalent to $\mathrm{P}_1 \tilde{F}$. By definition, $$\mathrm{cofib}(\mathrm{P}_1F(\emptyset) \rightarrow F(\star)) \simeq \Sigma \mathrm{D}_1 F(\emptyset) \simeq \Sigma \mathrm{P}_1 \tilde{F}(\emptyset).$$ We also observe that $\tilde{F}(\emptyset)=\Sigma^{-1} F(\star)$, since $F(\emptyset)$ is contractible by assumption.

    On the other hand, since $\tilde{F}$ is reduced, we can apply \cref{good_approx}. At emptyset, the first approximation map \[\mathrm{p}_1 \tilde{F} : \tilde{F}(\emptyset)\rightarrow \mathrm{P}_1 \tilde{F}(\emptyset)\] is equivalent to the inclusion in the colimit \[\tilde{F}(\emptyset) \rightarrow \colim_{n \rightarrow \infty} \Sigma^{-n} \tilde{F}(\Sigma^n_{\mathcal{S}} \emptyset).\] After applying suspension, we see the map $\Sigma \tilde{F}(\emptyset) \simeq F(\star) \rightarrow \Sigma \mathrm{P}_1 \tilde{F}(\emptyset)$ is equivalent to the map $\eta(\star)$. We conclude that $\mathrm{cofib}(\mathrm{P}_1 F(\emptyset) \rightarrow F(\star))$ is equivalent to $A$. Finally, the square \eqref{square-A-partial} is coCartesian and \cref{cor 4.8} concludes the proof. 
\end{proof}

A direct corollary of Proposition \ref{5.15} gives a criterion to determine whether a functor $F : \mathcal{S} \rightarrow \mathrm{Sp}$ is excisive or not:

\begin{cor} \label{crit-ex}
    Let $F : \mathcal{S} \rightarrow \mathrm{Sp}, \eta$ be as in \cref{5.15}. If the map $\eta(\star)$ is not an equivalence then the functor $F$ is not excisive.
\end{cor}

\begin{proof}
    According to \cref{5.15}, $\mathrm{P}_1F(\emptyset)$ is equivalent to $\mathrm{fib}(\eta(\star))$. Since by assumption the map $\eta(\star)$ is not an equivalence, it follows that $\mathrm{P}_1 F(\emptyset)$ is not contractible. However, if $F$ were excisive then $\mathrm{P}_1F(\emptyset)$ would be equivalent to $F(\emptyset)$ which is contractible by assumption. 
    \end{proof}

In the last part of this subsection, let $\mathcal{D}$ denote either $\mathcal{S}$ or $\mathrm{Sp}$ again. We discuss properties of the functor $\mathrm{P}_1$. The categories $\mathcal{S}$ and $\mathcal{D}$ being complete and cocomplete, we can compute colimits and limits in $\mathrm{Fun}(\mathcal{S},\mathcal{D})$ pointwise. Furthermore, because pullbacks commute with limits and with filtered colimits in $\mathcal{D}$, the category $\mathrm{Exc}_1(\mathcal{S},\mathcal{D})$ is closed under pullbacks and filtered colimits. Using the functor $\mathrm{P}_1$ is a left adjoint and the usual rules for commuting limits and colimits, we get the following proposition:  

\begin{prop}[Propositions $1.7$, $1.18$ in \cite{Goodwillie_2003} and Remark $6.1.1.32$ in \cite{HA} ] \label{P_1-commutes}
    The functors $\Pone^* : \mathrm{Fun}(\mathcal{S},\mathcal{D}) \rightarrow \mathrm{Exc}_1(\mathcal{S},\mathcal{D})$ and $\mathrm{D}_1^*: \mathrm{Fun}(\mathcal{S},\mathcal{D}) \rightarrow \mathrm{Exc}_1^{\mathrm{red}}(\mathcal{S},\mathcal{D})$ commute with: \begin{enumerate}
        \item finite limits, in particular with fiber sequences;
        \item filtered colimits;
        \item all colimits if $\mathcal{D}$ is $\mathrm{Sp}$;
        \item the functor $\Omega^{\infty} : \mathrm{Sp} \rightarrow \mathcal{S}$;
        \item the functors $\Sigma^{-1}, \Sigma : \mathrm{Sp} \rightarrow \mathrm{Sp}$.
    \end{enumerate}
\end{prop}

We conclude with the following remark:

\begin{rk}
    In particular, taking homotopy orbits $(-)\sslash G$ is the same as taking the colimit $\mathrm{colim}_{\B G}(-),$ where $G$ is a grouplike monoid. Consequently, we can consider a functor which can be written as $F\sslash G,$ where $F : \mathcal{S} \rightarrow \mathcal{D}$. Then, according to \cref{P_1-commutes}, $\mathrm{P}_1 (F\sslash G) \simeq (\mathrm{P}_1 F)\sslash G$.
\end{rk}

\subsection{Parametrized Pontryagin-Thom construction as a best approximation map} \label{5.3}

Let $P$ be a path-connected Poincaré complex. Let $H \rightarrow \haut(P)$ be a monoid map. Then the monoid $H$ acts by precomposition on $\Map(P,X)$. Let $\FPH : \mathcal{S} \rightarrow \mathrm{Sp}$ be the functor defined on objects by $\FPH(X)= \susplus (\Map(P,X)\sslash H)$. The functor $\FPH$ is defined on morphisms by post-composition on the spaces of maps. Unless otherwise specified we may simply write $\mathrm{F_P}$ instead of $\mathrm{F_P^{\{\text{id}_P\}}}$. We note that  $\FPH$ preserves filtered colimits, since finite spaces are compact objects in $\mathcal{S}$.  In this subsection, we aim to compute the first polynomial approximation $\Pone \FPH$ of $\FPH$.

Let $p_P^H$ be the $P$-fibration $P\sslash H \rightarrow BH$. We denote by $\DPH$ the dualizing complex of $p_P^H$ as defined in \cref{dualizingcomplexfamilies}.  We denote by $\PTPH$ the parametrized Pontryagin-Thom construction map $\mathrm{PT}_{p_P^H}$. Using Section \ref{4}, we construct in \cref{constr-parpt} a natural transformation  $$\PTPH(-) : \FPH(-) \Rightarrow  (r_!(\DPH) \otimes (\susplus -))$$ such that the map $\PTPH(\star)$ is equivalent to the parametrized Pontryagin-Thom construction map $\PTPH$ from \cref{notation-4.1.12}.

The terminal map induces a map $ r_!(\DPH) \otimes \susplus X \rightarrow  r_!(\DPH)$. At the point there is the parametrized Pontryagin-Thom map $\PTPH  : \FPH(\star) \simeq  \susplus \mathrm{BH} \rightarrow  r_!(\DPH)$. Let $\mathrm{E}_{\mathrm{P}}^{\mathrm{H}}$ be the excisive functor given by the pullback $$\xymatrix{\mathrm{E}_{\mathrm{P}}^{\mathrm{H}}(-) \ar[r] \ar[d] &   (r_!(\DPH) \otimes (\susplus -)) \ar[d]\\  \susplus \mathrm{BH} \ar[r]^{\PTPH} &  r_!(\DPH)}.$$ The natural transformation $\PTPH(-)$ and the natural transformation $\FPH \Rightarrow \FPH(\star)$ induce a natural transformation $\FPH \Rightarrow \mathrm{E}_{\mathrm{P}}^{\mathrm{H}}$.

In this subsection, we aim to prove the following proposition: 

\begin{prop} \label{main-prop-P1}
    The natural transformation $$\FPH \Rightarrow \mathrm{E}_{\mathrm{P}}^{\mathrm{H}}$$ induced by the parametrized Pontryagin-Thom construction is equivalent to the first polynomial approximation map $\mathrm{p}_1 \FPH: \FPH \Rightarrow \Pone \FPH$. In particular, the first derivative of $\FPH$ is equivalent to the Thom spectrum $r_!(\mathrm{D_P^H})$.
\end{prop}

We start with constructing the natural transformation $\PTPH(-)$. 

\begin{constr} \label{constr-parpt}
    For each $X$, we consider the map $\mathrm{p_P^H}(X) :( P\times \Map(P,X))\sslash H \rightarrow \Map(P,X)\sslash H$ with fiber $P$. We constructed in \cref{dualizingcomplexfamilies} a parametrized Pontryagin-Thom map $\susplus \Map(P,X)\sslash\mathrm{H} \rightarrow r_!(\mathrm{D_{p^H(X)}})$. We can show this construction is actually natural in $X$. Indeed, for $f : X \rightarrow Y$ a map, there is a pullback square: 
 \[\begin{tikzcd}
	{(\Map(P,X)\times P)\sslash H} & {(\Map(P,Y)\times P)\sslash H} \\
	{\Map(P,X)\sslash H} & {\Map(P,Y)\sslash H}
	\arrow[from=1-1, to=1-2]
	\arrow[from=1-1, to=2-1]
	\arrow[from=1-2, to=2-2]
	\arrow[from=2-1, to=2-2]
\end{tikzcd}.\]
By applying \cref{comm-PT}, we observe that the $\mathrm{PT}_{p_{P}^H(X)} : \susplus \Map(P,X)\sslash H \rightarrow r_!(D_{p_P^H(X)})$ give a natural transformation $$\FPH \Rightarrow r_!(D_{p_P^H(-)}).$$

We now construct a natural transformation $r_!(\mathrm{D_{p^H_P(-)}}) \Rightarrow r_!(\DPH) \otimes (\susplus -)$. Let $\pi_X$ denote the projection $(P\times \Map(P,X))\sslash H \rightarrow P\sslash H$. In particular, according to \cref{dua-pull}, the dualizing complex $D_{p^H_P(X)}$ is equivalent to $\pi_X^* \DPH$. On the other hand, evaluation gives a map $\mathrm{ev}_X : (P\times \Map(P,X))\sslash H \rightarrow X$. Then, the map $\pi_X$ factors as:  
\[\begin{tikzcd}
	{(P\times \Map(P,X))\sslash H} & {P\sslash H \times X} & {P\sslash H}
	\arrow["{\pi_X \times \mathrm{ev}_X}"', from=1-1, to=1-2]
	\arrow["{\pi_X}"{description}, curve={height=18pt}, from=1-1, to=1-3]
	\arrow[from=1-2, to=1-3]
\end{tikzcd}.\] We can then say $$\mathrm{D}_{\mathrm{p_P^H}(X)} \simeq (\pi_X \times \mathrm{ev}_X)^*(\DPH \times X),$$  where $\DPH \times X$ is the product fibration on $P\sslash H \times X$. The map $(\pi_X \times \mathrm{ev}_X)$ induces a map of Thom spectra $$r_!(\mathrm{D}_{\mathrm{p_P^H}(X)}) \rightarrow r_!(\DPH \times X) \simeq r_!(\DPH) \otimes (\susplus X).$$ Since this construction is natural in $X$, we get a natural transformation $$r_!(\mathrm{D}_{\mathrm{p_P^H}(-)}) \Rightarrow r_!(\DPH) \otimes (\susplus -).$$

Composing the two natural transformations $\FPH \Rightarrow r_!(D_{p_P^H(-)})$ and $r_!(\DPH) \otimes (\susplus -)$ gives a natural transformation $$\FPH \Rightarrow r_!(\DPH) \otimes (\susplus -).$$
\end{constr}

\begin{definition}
    Let $\PTPH(-) : \FPH(-) \rightarrow r_!(\DPH)\otimes (\susplus -)$ denote the natural transformation from \cref{constr-parpt}. By definition, at the point, we get the parametrized Pontryagin-Thom map $\PTPH : \susplus \mathrm{BH} \rightarrow r_!(\DPH)$.
\end{definition}

In what follows, let $\tilde{\mathrm{F}}_P^H$ be the reduction $\mathrm{fib}(\FPH(-) \rightarrow \FPH(\star))$ of $\FPH$. Let $\tilde{\mathrm{F}}_\mathrm{P}$ be the reduction $\tilde{\mathrm{F}}_{\mathrm{P}}^{\mathrm{id_P}}$ of $\mathrm{F_P}$. Then, the functor $\tilde{\mathrm{F}}_P$ is equivalent to the functor $ \widesus (\Map(P,-))$. Similarly the functor $\tilde{\mathrm{F}}_P^H$ is equivalent to the functor $ \widesus (\Map(P,-)\sslash H)$. 

Before proving Proposition \ref{main-prop-P1}, we begin with treating the case $H=\{ \text{id}_P\}$, where we write $\mathrm{E_P}$ instead of $\mathrm{E_P^{\{\mathrm{id}\}}}$. The derivatives of the functor $\mathrm{F_P}$ were already computed by Goodwillie and Arone in \cite{Goodwillie_2003} and \cite{arone2019goodwillie}. We give here a less general proof of the computation of $\partial_1 F_P$ than the one written in \cite{Goodwillie_2003} or \cite{arone2019goodwillie}.

\begin{prop} \label{P_1_F_P}
    The natural transformation $\mathrm{F_P} \Rightarrow \mathrm{E_P}$ is equivalent to the first approximation map $\mathrm{p}_1 \mathrm{F_P} : \mathrm{F_P} \Rightarrow \Pone \mathrm{F_P}$. In particular, the first derivative $\partial_1 \mathrm{F_P}$ is equivalent to $\mathrm{D}(P_+)$.
\end{prop}

\begin{proof}
According to \cref{5.15}, it suffices to show the map $$\mathrm{PT_P} : \mathrm{F_P}(\star)\simeq \bS \rightarrow \mathrm{D}(P_+)$$ is equivalent to the map:
\[\bS \rightarrow \colim_n \Sigma^{-n+1} \tilde{\mathrm{F}}(\Sigma^n_{\mathcal{S}} \emptyset) \simeq \colim_n \Sigma^{-n+1} \widesus \Map(P,\mathrm{S}^{n-1}).\]

The spectrum $\mathrm{D}(P_+)$ can be presented as a sequential spectrum $(\Map_*(P_+,S^n))_n$. Let $d$ be the dimension of $P$. The maps  $(\Sigma \Map_*(P_+,S^n) \rightarrow \Map_*(P_+,S^{n+1}))$ induce maps of spectra $\Sigma^{-n} \sus \Map_*(P_+,S^n) \rightarrow \Sigma^{-n+1} \sus \Map_*(P_+,S^{n+1})$. On the other hand, there are induced maps $$f_n : \Sigma^{-n}\sus \Map_*(P_+,\Sigma^n X) \rightarrow \mathrm{D}(P_+).$$ The maps $f_n$ assemble into a map \[\colim_n \Sigma^{-n} \sus \Map_*(P_+,S^n) \rightarrow D(P_+),\] which is an equivalence.

Secondly, the parametrized Pontryagin-Thom map $\bS \rightarrow D(P_+)$ is equivalent to the map $D(*_+) \simeq \bS \rightarrow D(P_+)$ induced by the terminal map. In particular, at the level of sequential spectra the map $\mathrm{PT}_P$ is induced by suspending the map $$\xymatrix{\Map_*(*_+,S^0) \simeq S^0 \ar[r]^{\simeq} & \Map_*(P_+,S^0) \simeq S^0}.$$ We see in fact $\mathrm{PT}_P$ is equivalent to the map $f_0 : \sus \Map_*(P_+,S^0) \rightarrow D(P_+)$, which is equivalent to the inclusion in the colimit \[\sus \Map_*(P_+,S^0) \rightarrow \colim_n \Sigma^{-n} \sus \Map_*(P_+,\mathrm{S}^{\mathrm{n}})\] Using $\widesus X \simeq \sus X$ for $X$ a based space and properties of shifts, we see that $\mathrm{PT}_P$ is equivalent to the inclusion in the colimit \[\bS \rightarrow \colim_n \Sigma^{-n+1} \tilde{F}_P(\Sigma^n_{\mathcal{S}} \emptyset) \simeq \colim_n \Sigma^{-n+1} \widesus \Map(P,\mathrm{S}^{n-1}).\]
\end{proof}

Since taking homotopy orbits is a colimit, \cref{P_1-commutes} gives a formula for $\mathrm{P}_1 \FPH$: 

\begin{cor} \label{cor-4-20}
    The natural transformation $F_P \Rightarrow E_P$ is $H$-equivariant for every $H$ a grouplike submonoid of $\haut(P)$. The induced natural transformation $\FPH \Rightarrow E_P\sslash H$ is equivalent to the first approximation map $$p_1 : \FPH \Rightarrow \mathrm{P}_1 \FPH.$$ In particular there is an equivalence $\partial_1 \FPH \simeq D(P_+)\sslash H$ where $H$ acts by precomposition. 
\end{cor}

We now prove \cref{main-prop-P1}:

\begin{proof}[Proof of Proposition \ref{main-prop-P1}]
    The case $H =\{\mathrm{id}_P\}$ was already treated in \cref{P_1_F_P}.
    
    There is a natural transformation $\PTPH :\FPH \Rightarrow r_!(\DPH) \otimes (\susplus -)$ which factors through a natural transformation $\mathrm{P}_1 \PTPH : \mathrm{P}_1 \FPH \Rightarrow r_!(\DPH) \otimes (\susplus -)$. At the point,  $(\mathrm{P}_1 \FPH)(\star)$ is equivalent to $\susplus \mathrm{BH}$ and $\mathrm{P}_1 \PTPH(\star)$ is equivalent to the map $\PTPH(\star)$.
    
    After evaluating the functor $\mathrm{E}_{\mathrm{P}}^{\mathrm{H}}$ at $\emptyset$, we obtain the following commutative square:
    \begin{equation} \label{squareproof}
        \begin{gathered}
             \xymatrix{\mathrm{P}_1 \FPH(\emptyset) \ar[r] \ar[d]^f & \star \ar[d]\\ \susplus \mathrm{BH} \ar[r]^{\PTPH} & r_!(\DPH)}
        \end{gathered}.
    \end{equation}
   
    The map $f : \mathrm{P}_1 \FPH(\emptyset) \rightarrow \susplus \mathrm{BH}$ is induced by the initial map.
    
    The category $\mathrm{Sp}^{BH}$ of parametrized spectra over $BH$ is equivalent to the category of spectra with an action of $H$. Taking homotopy orbits $-\sslash H$ is equivalent to taking the colimit functor $r_!= \text{colim}_{BH}$.
    
    According to \cref{cor-4-20}, the map $f$ is obtained as the homotopy quotient by $H$ of the map $g : \mathrm{P}_1 F_P(\emptyset) \rightarrow F_P(\star)$, where we recall $F_P(\star)\simeq \bS$. Since $\bS$ is the sphere spectrum with a trivial $H$-action, it corresponds to the constant spectrum $\bS_{BH}$ in $\mathrm{Sp}^{BH}$. Since $H$ acts on $\mathrm{P}_1 F_P(\emptyset)$, $\mathrm{P}_1F_P(\emptyset)$ is an object of $\mathrm{Sp}^{BH}$. The map $g$ can then be seen as a map $g: \mathrm{P}_1 F_P(\emptyset) \rightarrow \bS_{BH}$ in $\mathrm{Sp}^{BH}$. 
    
    According to \cref{dualizingcomplexfamilies}, the map $\PTPH$ is obtained from applying $r_!(-)$ to the unit map $$\epsilon :\bS_{BH} \rightarrow \mathrm{p}_!(\DPH)$$ in $\mathrm{Sp}^{BH}$. We then have the following commutative diagram in $\mathrm{Sp}^{BH}$:
    \begin{equation} \label{diag-SpBH}
        \begin{gathered}
            \xymatrix{\mathrm{P}_1F_P(\emptyset) \ar[r] \ar[d]^{g} & \star \ar[d] \\ \bS_{BH} \ar[r]^{\epsilon} & \mathrm{p}_!(\DPH)}
        \end{gathered}.
    \end{equation}
    
    We evaluate the maps $\epsilon$ and $g$ at each point $x\in \mathrm{BH}$. By construction, according to \cref{dualizingcomplexfamilies}, the map $\epsilon_x : \bS \rightarrow (\mathrm{p}_!(\DPH))_x \simeq D(P_+)$ is equivalent to the Pontryagin-Thom map $\mathrm{PT}_P : \bS \rightarrow D(P_+)$ for $P$. At each point $x\in \mathrm{BH}$, the map $g_x$ is equivalent to the map $\mathrm{P}_1 F_P(\emptyset) \rightarrow \bS\simeq F_P(\star)$ induced by the initial map. At each point $x$, we then have a commutative diagram: 
    $$\xymatrix{\mathrm{P}_1F_P(\emptyset) \ar[r] \ar[d]^{g_x} & \star \ar[d] \\ \bS \ar[r]^{\epsilon_x} &D(P_+)}.$$ According to \cref{P_1_F_P} and \cref{cor 4.8}, this square is coCartesian. Since the square \eqref{diag-SpBH} is pointwise coCartesian, we can conclude the square \eqref{diag-SpBH} is coCartesian in $\mathrm{Sp}^{BH}$.
    
    The square \eqref{squareproof} is obtained from the square \eqref{diag-SpBH} by applying the functor $r_!=\text{colim}_{BH}$. Since colimits commute with colimits, the square $$ \xymatrix{\Pone \FPH(\emptyset) \ar[r] \ar[d]^f & \star \ar[d]\\ \susplus \mathrm{BH} \ar[r]^{\PTPH} & r_!(\DPH)}$$ is also coCartesian. Applying \cref{cor 4.8} allows to conclude the natural transformation $$\Pone \FPH \Rightarrow \mathrm{E}_{\mathrm{P}}^{\mathrm{H}}$$ is an equivalence of functors. 
\end{proof}

We end this section with the following corollary, answering the question of whether $\FPH$ is excisive.
\begin{cor} 
    Assume $P$ is connected. The functor $\FPH$ is excisive if and only if $P$ is a point.
\end{cor}

\begin{proof}
    The functor $F_{\star}$ is the functor $\susplus  : \mathcal{S} \rightarrow \mathrm{Sp}$, which is excisive. On the other hand, assume $P$ is connected and not contractible. Then $P$ is of dimension greater than $1$. The Thom spectrum $r_!(\DPH)$ has negative homotopy groups, since the spherical fibration $\DPH$ has negative rank. Then, the Pontryagin-Thom map $\susplus \mathrm{BH} \rightarrow r_!(\DPH)$ is not an equivalence, because $\susplus \mathrm{BH}$ is a connective spectra. We conclude thanks to \cref{crit-ex} that $\FPH$ is not excisive.
    \end{proof}

\section{The best excisive approximation of \texorpdfstring{$\B \sCobtwo(-)$}{B sCobtwo(-)}} \label{6}

In this section, we aim to prove the two remaining theorems announced in the introduction: \cref{thmA} and \cref{thmC}.

\subsection{Proof of Theorem \ref{thmA}} \label{6.2}

We start with proving \cref{thmA}. In the first part of this subsection, we use results from \cref{4} and Subsection \ref{5.3} to construct the natural transformation announced in the introduction $$\alpha(-) : \B \sCobtwo (-) \Rightarrow \Omega^{\infty} (\Sigma \Th(\nuparStwo) \otimes (\susplus -)).$$ According to \cref{cor-delooping-nerve}, there is an equivalence of functors $\eta : \B \sCobtwo(-) \Rightarrow \lop \Sigma \mathrm{PH}(2,-)$. On the other hand, the functor $\mathrm{PH}(2,-) : \mathcal{S} \rightarrow \mathrm{Sp}$ is given by the following pushout square: 

\begin{equation} \label{6.1-square-pushout}
\begin{tikzcd}
	{\susplus (\Map(S^2,-)\sslash\mathrm{SO}(3))} & {\mathrm{MTSO}(2)\otimes (\susplus -)} \\
	{\susplus (\Map(S^2,-)\sslash\haut^+(S^2))} & {\mathrm{PH}(2,-)}
	\arrow[""{pos=0.6}, from=1-1, to=1-2]
	\arrow[from=1-1, to=2-1]
	\arrow[from=1-2, to=2-2]
	\arrow[from=2-1, to=2-2]
\end{tikzcd}.
\end{equation} According to \cite{GMTW}, the top map is the composite of the natural transformation $$\mathrm{PT}_{S^2}^{\Diff} : \susplus \Map(S^2,-)\sslash\mathrm{SO}(3) \Rightarrow \mathrm{Th}(\nu_{S^2}^{\Diff}) \otimes (\susplus -)$$ with a natural transformation $f \otimes (\susplus -) : \Th(\nu_{S^2}^{\Diff}) \otimes (\susplus - ) \Rightarrow \mathrm{MTSO}(2) \otimes (\susplus -)$ induced  by a map $f : \Th(\nu_{S^2}^{\Diff}) \rightarrow \mathrm{MTSO}(2)$. The bundle $\nu_{S^2}^{\Diff}$ is stably inverse of the vertical tangent bundle $T\pi$ of the universal bundle $S^2\sslash\mathrm{SO}(3) \rightarrow \mathrm{BSO}(3)$. Thus, the bundle $(-T\pi)$ is pulled back from the stable inverse of the universal $2$-dimensional vector bundle $\gamma_2$ over $\mathrm{BSO}(2)$. The map $f : \Th(\nu_{S^2}^{\Diff}) \rightarrow \mathrm{MTSO}(2)$ is then the induced map on Thom spectra. 

In the following lemma, we show that the spectra $\Th(\nu_{S^2}^{\Diff})$ and $\mathrm{MTSO}(2)$ are actually equivalent: 

\begin{lem} \label{MTSO2S^2}
    The map $f : \Th(\nu_{S^2}^{\Diff}) \rightarrow \mathrm{MTSO}(2)$ is an equivalence.
\end{lem}

\begin{proof}
    To begin with, we note that the universal $S^2$-bundle is given up to homotopy by 
\[\begin{tikzcd}
	{S^2} & {\mathrm{BSO}(2)} & {\mathrm{BSO}(3)}
	\arrow["j", from=1-1, to=1-2]
	\arrow["\pi", from=1-2, to=1-3]
\end{tikzcd}.\] Let $T\pi$ denote the vertical tangent bundle of dimension $2$ on the total space $\mathrm{BSO}(2)$. It is classified by a map $T\pi : \text{BSO}(2) \rightarrow \text{BSO}(2)$. If $T\pi$ is homotopic to the identity, then $T\pi$ is the universal bundle $\gamma_2$ over $\mathrm{BSO}(2)$. The statement would follow from $(-T\pi)$ and $(-\gamma_2)$ being equivalent. 

    The pullback of $T\pi$ along the inclusion of the fiber $j : S^2 \rightarrow \text{BSO}(2)$ is isomorphic to the tangent bundle of $S^2$ classified by a map $TS^2 : S^2 \rightarrow \text{BSO}(2)$. At the level of classifying maps for bundles there is an equivalence $T\pi \circ j \simeq TS^2$. 
    
    We have an isomorphism $[S^2,\text{BSO}(2)] \cong \Ho^2(S^2, \bZ) \cong \bZ$. The bundle $TS^2$ is classified by its Euler class $e(TS^2) \in \Ho^2(S^2,\bZ)$. The class $e(TS^2)$ is given by $2.u$ where $u$ is a generator of $\Ho^2(S^2,\bZ)$.
  
    Unwinding the long exact sequence of homotopy groups for $\pi$ we recover that the map induced on $\pi_2$ by $j$ is the multiplication by $2$, hence $j : S^2 \rightarrow \text{BSO}(2)$ corresponds as well to $2.u$ in $\Ho^2(S^2,\bZ)$. We conclude $T\pi$ is homotopic to the identity.
\end{proof}

To construct the natural transformation $\alpha : \B\sCobtwo(-) \Rightarrow \lop (\Sigma \Th(\nuparStwo) \otimes \susplus -)$, it suffices to define a natural transformation $\mathrm{PH}(2,-) \Rightarrow \Th(\nuparStwo) \otimes (\susplus -)$. 

The bundle $\pi : \mathrm{BSO}(2) \rightarrow \mathrm{BSO}(3)$ is a pullback of the universal fibration $\mathrm{Bhaut}^+_*(S^2) \rightarrow \mathrm{Bhaut}^+(S^2)$. According to \cref{dua-pull}, the parametrized Spivak fibration $\nu_{S^2}^{\Diff}$ of $\pi$ is pulled back from $\nuparStwo$ along the map $\mathrm{BSO}(3) \rightarrow \mathrm{Bhaut}^+(S^2)$. The latter induces then a map $\Th(\nu_{S^2}^{\Diff}) \rightarrow \mathrm{Th}(\nuparStwo)$, hence a map $j :  \mathrm{MTSO}(2) \rightarrow \Th(\nuparStwo)$ according to \cref{MTSO2S^2}.

We then have two natural transformations $j : \mathrm{MTSO}(2) \otimes (\susplus -) \Rightarrow \Th(\nuparStwo) \otimes (\susplus -)$ and $\mathrm{PT}_{S^2}^{\haut} : \susplus (\Map(S^2,-)\sslash\mathrm{haut}(S^2)) \Rightarrow \Th(\nuparStwo) \otimes (\susplus -)$ from \cref{constr-parpt}. In the following lemma, we show they are both compatible with the maps out of $\susplus (\Map(S^2,-)\sslash\mathrm{SO}(3))$ in the square \eqref{6.1-square-pushout}:

\begin{lem}
    The square of natural transformations 
\[\begin{tikzcd}
	{\susplus (\Map(S^2,-)\sslash\mathrm{SO}(3))} & {\mathrm{MTSO}(2)\otimes (\susplus-)} \\
	{\susplus (\Map(S^2,-)\sslash\mathrm{haut}^+(S^2)} & {\mathrm{Th}(\nuparStwo)\otimes (\susplus -)}
	\arrow["{\mathrm{PT}_{S^2}^{\Diff}}", from=1-1, to=1-2]
	\arrow[from=1-1, to=2-1]
	\arrow["j", from=1-2, to=2-2]
	\arrow["{\mathrm{PT}_{S^2}^{\haut}}"', from=2-1, to=2-2]
\end{tikzcd}\] is commutative.
\end{lem}

\begin{proof}
    As in \cref{comm-PT}, this square is the map of parametrized Pontryagin-Thom constructions induced by the diagram 
\[\begin{tikzcd}
	{\mathrm{BSO}(2)} & {\mathrm{Bhaut}_*^+(S^2)} \\
	{\mathrm{BSO}(3)} & {\mathrm{Bhaut}^+(S^2)}
	\arrow[from=1-1, to=1-2]
	\arrow[from=1-1, to=2-1]
	\arrow[from=1-2, to=2-2]
	\arrow[from=2-1, to=2-2]
\end{tikzcd}.\] The conclusion follows from \cref{comm-PT}.
\end{proof}

Since the square \eqref{6.1-square-pushout} is a pushout in $\mathrm{Sp}$, we can construct a natural transformation $\mathrm{PH}(2,-) \Rightarrow \Th(\nuparStwo) \otimes (\susplus -)$ as follows: 

\begin{constr} \label{constr-6.1-nattrans}
    Let $\beta : \mathrm{PH}(2,-) \Rightarrow \Th(\nuparStwo) \otimes (\susplus -)$ be the natural transformation induced by the following diagram 
\[\begin{tikzcd}
	{\susplus (\Map(S^2,-)\sslash\mathrm{SO}(3))} & {\mathrm{MTSO}(2)\otimes (\susplus -)} \\
	{\susplus (\Map(S^2,-)\sslash\haut^+(S^2))} & \textcolor{mypurple}{{\mathrm{PH}(2,-)}} \\
	&& \textcolor{mypurple}{{\Th(\nuparStwo) \otimes (\susplus -)}}
	\arrow[from=1-1, to=1-2]
	\arrow[from=1-1, to=2-1]
	\arrow[from=1-2, to=2-2]
	\arrow["j"{description},curve={height=-12pt}, Rightarrow, from=1-2, to=3-3]
	\arrow[from=2-1, to=2-2]
	\arrow["\mathrm{PT}^{\haut}_{S^2}"{description},curve={height=12pt}, Rightarrow, from=2-1, to=3-3]
	\arrow["\beta"{description}, color=mypurple, Rightarrow, from=2-2, to=3-3]
\end{tikzcd}.\] 
Let $\alpha$ denote the natural transformation $$\lop \Sigma \beta \circ \eta  : \B \sCobtwo(-) \Rightarrow \lop \Sigma \mathrm{PH}(2,-) \Rightarrow \lop (\Sigma\Th(\nuparStwo) \otimes (\susplus -)).$$
\end{constr}

We now have all the elements to establish the proof of \cref{thmA}.

\begin{proof}[Proof of \cref{thmA}]
    According to \cref{cor-delooping-nerve}, the natural transformation $\eta : \B \sCobtwo(-) \Rightarrow \lop \Sigma \mathrm{PH}(2,-)$ is an equivalence. Then, according to \cref{P_1-commutes}, $P_1B \sCobtwo(-)$ is equivalent to $\lop \Sigma \circ \mathrm{P}_1 \mathrm{PH}(2,-)$, hence it suffices to compute $\mathrm{P}_1 \mathrm{PH}(2,-)$. 

   Let $D$ be the $1$-excisive functor given by the pullback $$\xymatrix{D \ar[r] \ar[d] & \Th(\nuparStwo) \otimes (\susplus -) \ar[d]^{(-\rightarrow \star)} \\ \mathrm{PH}(2,\star) \simeq \mathrm{PH}(2) \ar[r]^{\beta(\star)}  & \Th(\nuparStwo)},$$ where the right vertical map is induced by the terminal map. The natural transformation $\beta : \mathrm{PH}(2,-) \Rightarrow \mathrm{Th}(\nuparStwo) \otimes (\susplus -)$ factors through a natural transformation $ \Pone \mathrm{PH}(2,-) \Rightarrow \Th(\nuparStwo) \otimes (\susplus -)$. We recall that $\Pone \mathrm{PH}(2,\star) \simeq \mathrm{PH}(2)$.

   We now consider the following square: 
\begin{equation} \label{proof-thmC-suare}\begin{tikzcd}
	{\Pone\mathrm{PH}(2,\emptyset)} & \star \\
	{\mathrm{PH}(2)} & {\Th(\nuparStwo)}
	\arrow[from=1-1, to=1-2]
	\arrow[from=1-1, to=2-1]
	\arrow[from=1-2, to=2-2]
	\arrow["{\beta(\star)}", from=2-1, to=2-2]
\end{tikzcd},
\end{equation} where the map $\mathrm{P}_1 \mathrm{PH}(2,\emptyset) \rightarrow \mathrm{PH}(2)$ is induced by the map $\emptyset \rightarrow \star$. According to \cref{cor 4.8}, if the square \eqref{proof-thmC-suare} is coCartesian, then the natural transformation $\Pone \mathrm{PH}(2,-) \Rightarrow D$ is an equivalence.
   
   The functor $\mathrm{PH}(2,-)$ is obtained as a pushout of functors $F_{S^2}^{\Diff}(-), F_{S^2}^{\haut}(-)$ as in Subsection \ref{5.3} and $\mathrm{MTSO}(2) \otimes (\susplus -)$. According to \cref{P_1-commutes}, taking $P_1$ commutes with colimits. The first polynomial approximation of $\mathrm{PH}(2,-)$ is then determined by the following pushout 
\[\begin{tikzcd}
	{\mathrm{P}_1 \mathrm{F}_{S^2}^{\Diff}(-)} & {\mathrm{MTSO}(2)\otimes (\susplus -)} \\
	{\mathrm{P}_1 \mathrm{F}_{S^2}^{\haut}(-)} & {\mathrm{P}_1 \mathrm{PH}(2,-)}
	\arrow[from=1-1, to=1-2]
	\arrow[from=1-1, to=2-1]
	\arrow[from=1-2, to=2-2]
	\arrow[from=2-1, to=2-2]
\end{tikzcd}.\]
 Let $H$ denote either $\Diff^+(S^2)$ or $\haut^+(S^2)$. At the point, $F_{S^2}^H(\star)$ and $\mathrm{MTSO}(2) \otimes (\susplus \star)$ are respectively equivalent to $\susplus \mathrm{BH}$ and $\mathrm{MTSO}(2)$. According to \cref{main-prop-P1}, the square 
\[\begin{tikzcd}
	{\mathrm{P}_1 \mathrm{F}_{S^2}^{\mathrm{H}}(\emptyset)} & \star \\
	{\susplus \mathrm{BH}} & \textcolor{mygreen}{{\mathrm{Th}(\nu_{S^2}^{\mathrm{H}})}}
	\arrow[from=1-1, to=1-2]
	\arrow[from=1-1, to=2-1]
	\arrow[from=1-2, to=2-2]
	\arrow["{\mathrm{PT}_{S^2}^{\mathrm{H}}}"', from=2-1, to=2-2]
\end{tikzcd}\] is coCartesian. Moreover, according to \cref{MTSO2S^2}, $\Th(\nu_{S^2}^{\Diff})$ is equivalent to $\mathrm{MTSO}(2)$.

Let $\textcolor{mygreen}{A}$ be the cofiber of the map $\mathrm{P}_1 \mathrm{PH}(2,\emptyset) \rightarrow \mathrm{PH}(2)$.

The map $\mathrm{P}_1 \mathrm{PH}(2,\emptyset) \rightarrow \mathrm{PH}(2)$ is induced by a map of pushout diagrams. This is represented on the left cube of the following diagram \eqref{big-diagram-proof-THMA}. The right hand-side is obtained by taking the cofibers of the horizontal maps $\mathrm{P}_1 \mathrm{F}_{S^2}^{\mathrm{H}}(\emptyset) \rightarrow \susplus \mathrm{BH}$, $\star \rightarrow \mathrm{MTSO}(2)$ and $\mathrm{P}_1 \mathrm{PH}(2,\emptyset) \rightarrow \mathrm{PH}(2)$. Since both squares labelled $(1)$ and $(2)$ in the diagram \eqref{big-diagram-proof-THMA} are coCartesian and taking cofibers commutes with pushouts, we deduce that the square labelled by $\textcolor{mygreen}{(3)}$ is a pushout square.

\begin{equation} \label{big-diagram-proof-THMA} 
\begin{tikzcd}
	{\mathrm{P}_1 \mathrm{F}_{S^2}^{\Diff}(\emptyset)} && {\susplus \mathrm{BSO}(3)} && \textcolor{mygreen}{\mathrm{MTSO}(2)} \\
	& \star && {\mathrm{MTSO}(2)} && \textcolor{mygreen}{\mathrm{MTSO}(2)} \\
	{\mathrm{P}_1 \mathrm{F}_{S^2}^{\haut}(\emptyset)} && {\susplus \mathrm{Bhaut}^+(S^2)} && \textcolor{mygreen}{\mathrm{Th}(\nuparStwo)} \\
	& {\mathrm{P}_1 \mathrm{PH}(2,\emptyset)} && {\mathrm{PH}(2)} && \textcolor{mygreen}{A}
	\arrow[from=1-1, to=1-3]
	\arrow[from=1-1, to=2-2]
	\arrow[from=1-1, to=3-1]
	\arrow["{(1)}"{description}, draw=none, from=1-1, to=4-2]
	\arrow["{\mathrm{PT}_{S^2}^{\Diff}}"{description}, from=1-3, to=1-5]
	\arrow["{\mathrm{PT}_{S^2}^{\Diff}}"{description}, from=1-3, to=2-4]
	\arrow[dashed, from=1-3, to=3-3]
	\arrow["{(2)}"{description}, draw=none, from=1-3, to=4-4]
	\arrow["\simeq"{description}, color={mygreen}, from=1-5, to=2-6]
	\arrow["j"{description, pos=0.3}, color={mygreen}, dashed, from=1-5, to=3-5]
	\arrow["{(3)}"{description}, color={mygreen}, draw=none, from=1-5, to=4-6]
	\arrow[from=2-2, to=2-4]
	\arrow[from=2-2, to=4-2]
	\arrow["\simeq"{description, pos=0.2}, from=2-4, to=2-6]
	\arrow[from=2-4, to=4-4]
	\arrow["j"{description},color=mygreen,from=2-6, to=4-6]
	\arrow[dashed, from=3-1, to=3-3]
	\arrow[from=3-1, to=4-2]
	\arrow["{\mathrm{PT}^{\haut}_{S^2}}"{description, pos=0.7}, dashed, from=3-3, to=3-5]
	\arrow[from=3-3, to=4-4]
	\arrow[color=mygreen,from=3-5, to=4-6]
	\arrow[from=4-2, to=4-4]
	\arrow[from=4-4, to=4-6]
\end{tikzcd}
\end{equation}
We deduce from the square $\textcolor{mygreen}{(3)}$ being coCartesian that the map $\Th(\nuparStwo) \rightarrow A$ is an equivalence. We see in the diagram \eqref{big-diagram-proof-THMA} that the map $\mathrm{PH}(2) \rightarrow A$ is obtained as follows: 
\[\begin{tikzcd}
	{\susplus \mathrm{BSO}(3)} & {\mathrm{MTSO}(2)} & {\Th(\nuparStwo)} \\
	{\susplus \mathrm{Bhaut}^+(S^2)} & \textcolor{mygreen}{{\mathrm{PH}(2)}} \\
	{\Th(\nuparStwo)} && \textcolor{mygreen}{A}
	\arrow["{\mathrm{PT}_{S^2}^{\Diff}}", from=1-1, to=1-2]
	\arrow[from=1-1, to=2-1]
	\arrow["j", from=1-2, to=1-3]
	\arrow[from=1-2, to=2-2]
	\arrow["\simeq"{description}, from=1-3, to=3-3]
	\arrow[from=2-1, to=2-2]
	\arrow["{\mathrm{PT}_{S^2}^{\haut}}"', from=2-1, to=3-1]
	\arrow[color=mygreen, from=2-2, to=3-3]
	\arrow["\simeq"{description}, from=3-1, to=3-3]
\end{tikzcd}.\]

By \cref{constr-6.1-nattrans} of $\beta$, we conclude that the map $\mathrm{PH}(2) \rightarrow A$ is equivalent to $$\beta(\star) : \mathrm{PH}(2) \rightarrow \Th(\nuparStwo).$$ Consequently, the square \eqref{proof-thmC-suare} is coCartesian, which concludes the proof.

\end{proof}

\subsection{Proof of Theorem \ref{thmC}} \label{6.3}

In this final subsection, we prove \cref{thmC}. The upshot is that the map from Subsection \ref{6.2} $$\Omega \alpha(\star): \Omega_{\emptyset} \B \sCobtwo \rightarrow \Omega^{\infty} \Th(\nuparStwo)$$ is not a weak equivalence. According to \cref{constr-6.1-nattrans}, the map $\Omega \alpha(\star)$ comes from a map $\beta(\star) : \mathrm{PH}(2) \rightarrow \Th(\nuparStwo)$. 

More precisely, we first establish in \cref{thmC-spectra-version} a version of \cref{thmC} for the map of spectra $\beta(\star) : \mathrm{PH}(2) \rightarrow \Th(\nuparStwo)$. In other words, we show that the map $\beta(\star)  : \mathrm{PH}(2) \rightarrow \Th(\nuparStwo)$ induces an isomorphism on $\pi_*$ for $* \leq 0$ and we construct a nonzero class $\epsilon. U \in H^1(\Th(\nuparStwo),\Ftwo)$ such that $\beta(\star)^*(\epsilon.U)$ vanishes in $H^1(\mathrm{PH}(2),\Ftwo)$. We then conclude the proof of \cref{thmC} by propagating the latter result after taking $\lop(-)$.

Let $\mathbb{F}_2[\ldots]$ and $\Lambda[\ldots]$ respectively denote taking the polynomial algebra and the exterior algebra over $\Ftwo$. We start with describing the cohomology ring of $\B\haut^+_*(S^2)$ in the next proposition: 

\begin{prop} \label{cohomolo-bhaut-based}
      The $\Ftwo$-cohomology ring $$\Ho^*(\Bhaut_*^+(S^2), \Ftwo)$$ is isomorphic to $$\Ftwo[w_2] \otimes \Lambda[\epsilon, e_n]_{n\geq 2}$$ where: \begin{itemize}
        \item the classes $w_2$ is the second Stiefel-Whitney class;
        \item the class $\epsilon$ is in degree $3$;
        \item the classes $e_n$ are defined for $n\geq 2$ and have degree $2^n$.
    \end{itemize}
\end{prop}

The cohomology ring of $\B\haut(S^2)$ was already determined by Milgram in \cite{milgram-cohomo}: 
\begin{prop}[Theorem A in \cite{milgram-cohomo}] \label{milgram-bhaut}
    The $\Ftwo$-cohomology ring $$\Ho^*(\Bhaut(S^2), \Ftwo)$$ is isomorphic to $$\Ftwo[w_1,w_2,w_3] \otimes \Lambda[\epsilon, e_n]_{n\geq 2}$$ where: \begin{itemize}
        \item the classes $w_1, w_2,w_3$ are the Stiefel-Whitney classes in degree $1,2$ and $3$;
        \item the class $\epsilon$ is in degree $3$;
        \item the classes $e_n$ are defined for $n\geq 2$ and have degree $2^n$.
    \end{itemize}
\end{prop}

In particular, Stiefel-Whitney classes are also defined for spherical fibrations and the pullback map $$\Ho^*(\Bhaut(S^2), \Ftwo) \rightarrow \Ho^*(\text{BO}(3),\Ftwo)$$ is surjective. As a corollary, we briefly prove \cref{cohomolo-bhaut-based}.

\begin{proof}[Proof of Proposition \ref{cohomolo-bhaut-based}]
    Since $\B \haut^+(S^2)$ is simply connected, the class $w_1$ vanishes. In particular, we deduce $H^*(\B \haut^+(s^2),\Ftwo) \cong \Ftwo[w_2,w_3]\otimes \Lambda[\epsilon, e_n]_{n\geq 2}$.

    The homotopy fiber of $\pi : \Bhaut^+_*(S^2) \rightarrow \Bhaut^+(S^2)$ is $S^2$. According to the Thom-Gysin long exact sequence, there exists a class $c\in \Ho^3(\Bhaut^+(S^2))$ such that the following is exact:  
\[\begin{tikzcd}[row sep=tiny, column sep=small]
	\ldots & {\Ho^*(\Bhaut^+(S^2),\mathbb{F}_2)} & {\Ho^*(\Bhaut^+_*(S^2),\mathbb{F}_2)} & {\Ho^{*-2}(\Bhaut^+(S^2),\mathbb{F}_2)} & {} \\
	{} & {\Ho^{*+1}(\Bhaut^+(S^2),\mathbb{F}_2)} & {\Ho^{*+1}(\Bhaut^+_*(S^2),\mathbb{F}_2)} & \ldots & {}
	\arrow[from=1-1, to=1-2]
	\arrow["{\pi^*}", from=1-2, to=1-3]
	\arrow[from=1-3, to=1-4]
	\arrow[no head, from=1-4, to=1-5]
	\arrow["{- \smile c}", from=2-1, to=2-2]
	\arrow["{\pi^*}", from=2-2, to=2-3]
	\arrow[from=2-3, to=2-4]
\end{tikzcd}.\] 
Since $\pi$ induces an isomorphism on $H^2(-,\Ftwo)$, we deduce $c$ is non-zero. On the other hand, the class $w_3$ is an obstruction to a spherical fibration for having a section. Since $\B \haut_*^+(S^2)$ classifies spherical fibrations with a section, we deduce the class $w_3$ vanishes in $H^3(\B\haut^+_*(S^2))$. After inspection of the Thom-Gysin sequence, we deduce $c$ is $w_3$, which concludes the proof.
\end{proof}

Before finishing the proof of \cref{thmC}, we give an interpretation in the remark below of the class $\epsilon$:

\begin{rk}[Exotic characteristic classes]
    The characteristic classes coming from $\epsilon$ and $e_i$ in the cohomology $\Ho^*(\Bhaut^+(S^2),\Ftwo)$ vanish when evaluated onto vector bundles. We can speak of exotic characteristic classes. In \cite{HEIL1985269}, Heil constructs these classes via some secondary cohomology operations based on the Adem relations in the Steenrod algebra. There is a simpler interpretation of $\epsilon$ though, given by Gitler and Stasheff in \cite{GITLER1965257}. Let $f : B \rightarrow \Bhaut_*^+(S^2)$ be a spherical fibration. The first non-trivial obstruction $o_3(f)$ to lift $f$ to a vector bundle classified by $\tilde{f} : B \rightarrow \mathrm{BSO}(2)$ lives in $\Ho^3(B,\pi_2(\mathrm{fib}(\mathrm{BSO}(2) \rightarrow \Bhaut_*^+(S^2))).$ In particular, the obstruction $o_3(f)$ lives in $\Ho^3(B, \Ftwo)$. If $f$ is the identity of $\B \haut_*^+(S^2)$, $o_3$ is non-zero since $f$ classifies the universal fibration. It follows from Milgram's result and \cref{prop-approx-Bhaut} that the only non-zero class in $H^3(\B\haut_*^+(S^2),\Ftwo)$ is $\epsilon$.    

\end{rk}

We now wish to study the map $\beta(\star)$. For convenience, we write $\beta$ instead of $\beta(\star)$. The map $\beta$ fits in the following commutative diagram in $\mathrm{Sp}$: 
\[\begin{tikzcd}
	{\mathrm{MTSO}(2)} \\
	{\mathrm{PH}(2)} & {\Th(\nu^{\haut}_{S^2})} 
	\arrow["a"{description},from=1-1, to=2-1]
	\arrow["j"{description}, from=1-1, to=2-2]
	\arrow["{\beta}"{description}, from=2-1, to=2-2]
\end{tikzcd}.\] Here $j$ is as in \cref{constr-6.1-nattrans} and the map $a$ is as in the square \eqref{6.1-square-pushout}.

According to \cref{cohomolo-bhaut-based}, the cohomology group $H^3(\B \haut_*^+(S^2),\Ftwo)$ is generated by a class $\epsilon$. The Thom isomorphism produces a non-zero class $\epsilon. U \in H^1(\Th(\nuparStwo),\Ftwo)$, where $U\in H^{-2}(\Th(\nuparStwo))$ is the Thom class of $\nuparStwo$.

\begin{prop} \label{thmC-spectra-version}
    The map $\beta : \mathrm{PH}(2) \rightarrow \Th(\nuparStwo)$ \begin{itemize}
        \item is a rational equivalence;
        \item induces an isomorphism on $\pi_*$ for $*\leq 0$; 
    \end{itemize} however the class $\epsilon.U$ generates $H^1(\Th(\nuparStwo),\Ftwo)$, while $\beta^*(\epsilon.U)$ vanishes in $H^1(\mathrm{PH}(2),\Ftwo)$.
\end{prop}

The following lemma computes the connectivity of the map $a$: 

\begin{lem} \label{lemma-5.15}
    The map $a : \mathrm{MTSO}(2) \rightarrow \mathrm{PH}(2)$ is $2$-connected.
\end{lem}

\begin{proof}
    According to \cref{thmB}, the spectrum $\mathrm{PH}(2)$ is given as a pushout 
\[\begin{tikzcd}
	{\susplus \mathrm{BSO}(3)} & {\mathrm{MTSO}(2)} \\
	{\susplus \mathrm{Bhaut}^+(S^2)} & {\mathrm{PH}(2)}
	\arrow[from=1-1, to=1-2]
	\arrow["{\susplus \iota}"', from=1-1, to=2-1]
	\arrow["a", from=1-2, to=2-2]
	\arrow[from=2-1, to=2-2]
\end{tikzcd}.\] The cofiber of $\susplus \iota$ is equivalent to $\sus C$, where $C$ is as in \cref{comparison-lemma-cofibfib}. According to Freudenthal suspension theorem and \cref{comparison-lemma-cofibfib}, we deduce $\sus C$ is $2$-connected, hence the claim.
\end{proof}

We now show there is no difference rationally between $\mathrm{MTSO}(2)$ and $\Th(\nuparStwo)$:

\begin{lem} \label{lemma-5.16}
    The map $\beta : \mathrm{PH}(2)\rightarrow \Th(\nuparStwo)$ is a rational equivalence.  
\end{lem}

\begin{proof}
    On the one hand, the map $\iota : \text{BSO}(3) \rightarrow \Bhaut^+(S^2)$ is a rational equivalence. It follows from the square \eqref{6.1-square-pushout} being coCartesian, that the map $a : \mathrm{MTSO}(2) \rightarrow \mathrm{PH}(2)$ is a rational equivalence.
    On the other hand, the map $ \text{BSO}(2) \rightarrow \Bhaut^+_*(S^2)$ is also a rational equivalence. Hence, the induced map on Thom spectra $j : \mathrm{MTSO}(2) \rightarrow \Th(\nuparStwo)$ is also a rational equivalence. The map $\beta : \mathrm{PH}(2) \rightarrow \Th(\nuparStwo) $ is then a rational equivalence by a two-out-of-three argument.
\end{proof}

In the following lemma, we show the map $\beta : \mathrm{PH}(2) \rightarrow \Th(\nuparStwo)$ induces an isomorphism on homotopy groups in nonpositive degree:

\begin{lem} \label{lemma-5.17}
    The map $\beta : \mathrm{PH}(2) \rightarrow \Th(\nuparStwo)$ induces an isomorphism on $\pi_*$ for $* \leq 0$.
\end{lem}

\begin{proof} 
According to \cref{lemma-5.15}, the map $a : \mathrm{MTSO}(2) \rightarrow \mathrm{PH}(2)$ is $2$-connected. Consequently, the map $\beta: \mathrm{PH}(2) \rightarrow \Th(\nuparStwo)$ induces an isomorphism on $\pi_*$ for $*\leq 0$ if and only if $j : \mathrm{MTSO}(2) \rightarrow \Th(\nuparStwo)$ does. 

We now show the map $ j : \mathrm{MTSO}(2) \rightarrow \Th(\nuparStwo)$ induces an isomorphism on nonpositive homotopy groups $\pi_*$ for $*\leq 0$. We can apply the relative Atiyah-Hirzebruch spectral sequence to the map $ j : \mathrm{MTSO}(2) \rightarrow \Th(\nuparStwo)$. The $E^2$-page is given by the $\Ho_p(\Th(\nuparStwo), \mathrm{MTSO}(2),\pi_q(\bS))$ and converges to $\pi_{p+q}(\Th(\nuparStwo),\mathrm{MTSO}(2))$. The map $\text{BSO}(2) \rightarrow \Bhaut_*^+(S^2)$ is $2$-connected, hence the relative homology groups $$\Ho_*(\Bhaut_*^+(S^2),\text{BSO}(2))$$ vanish for $*\leq 2$. On the other hand, it follows from the Thom isomorphism that: $$\Ho_p(\Th(\nuparStwo), \mathrm{MTSO}(2),\pi_q(\bS))\cong \Ho_{p+2}(\Bhaut_*^+(S^2),\text{BSO}(2),\pi_q(\bS)).$$ Since $\bS$ is connective, we conclude the terms $E^2_{p,q}$ vanish in the spectral sequence for $p+q \leq 0$. As a consequence, the relative homotopy groups $\pi_*(\Th(\nuparStwo),\mathrm{MTSO}(2))$ vanish in degrees $*\leq 0$.

It remains to show the map $j$ is injective on $\pi_0$. In the $E^2$-page, since $$H_*(\Th(\nu_{S^2}^{\haut}), \mathrm{MTSO}(2))$$ vanishes for $*\leq 0$, there is only one non-zero element in the line $p+q=1$, given by $$H_1(\Th(\nuparStwo), \mathrm{MTSO}(2),\pi_0(\bS)).$$ We can deduce from Thom isomorphism and Hurewicz Theorem that $H_1(\Th(\nuparStwo), \mathrm{MTSO}(2),\bZ)$ is isomorphic to $\Ftwo$. According to \cite[Theorem $1.0.1$]{ebert2007lowdimensionalhomotopystablemapping}, the group $\pi_0(\mathrm{MTSO}(2))$ is isomorphic to $\bZ$.

The long exact sequence on homotopy groups of $j$ is then as follows:

\footnotesize {\[\begin{tikzcd}
	{} & {\pi_1(\Th(\nuparStwo),\mathrm{MTSO}(2))} & {\pi_0(\mathrm{MTSO}(2))} & {\pi_0(\Th(\nuparStwo))} & {\pi_0(\Th(\nuparStwo),\mathrm{MTSO}(2))} \\
	{} & \Ftwo & \bZ & \bZ & 0
	\arrow[from=1-1, to=1-2]
	\arrow[from=1-2, to=1-3]
	\arrow["\cong"{marking, allow upside down}, draw=none, from=1-2, to=2-2]
	\arrow[from=1-3, to=1-4]
	\arrow["\cong"{marking, allow upside down}, draw=none, from=1-3, to=2-3]
	\arrow[from=1-4, to=1-5]
	\arrow[from=1-4, to=2-4]
	\arrow["\cong"{marking, allow upside down}, draw=none, from=1-5, to=2-5]
	\arrow[from=2-1, to=2-2]
	\arrow["0",from=2-2, to=2-3]
	\arrow[from=2-3, to=2-4]
	\arrow[from=2-4, to=2-5]
\end{tikzcd}.\]} 

\normalsize After inspection, we infer $\pi_0(\mathrm{MTSO}(2)) \rightarrow \pi_0(\Th(\nuparStwo))$ is an isomorphism.

\end{proof}

We now have all the elements to prove \cref{thmC}: 

\begin{proof}[Proof of \cref{thmC}]
    According to \cref{lemma-5.16}, the map $\beta : \mathrm{PH}(2) \rightarrow \Th(\nuparStwo)$ is a rational equivalence. After taking $\Omega^{\infty}$, the map $$\Omega \alpha(\star) : \Omega^{\infty} \mathrm{PH}(2) \rightarrow \Omega^{\infty} \Th(\nuparStwo)$$ is also a rational equivalence.
    
    According to \cref{lemma-5.17}, the map $$\beta : \mathrm{PH}(2) \rightarrow \Th(\nuparStwo)$$ induces an isomorphism on $\pi_*$ for $*\leq 0$. Since $\Omega^{\infty}$ preserves homotopy groups, the map $\Omega \alpha(\star): \Omega_{\emptyset}\B\sCobtwo \rightarrow \Omega^{\infty} \Th(\nuparStwo)$ induces an isomorphism on $\pi_0$.
    
    To conclude, we need to show two things: the class $\kappa_{\epsilon} \in \Ho^1(\Omega^{\infty} \Th(\nuparStwo), \Ftwo)$ is non-zero and the class $(\Omega \alpha(\star))^*\kappa_{\epsilon}$ vanishes. 

    We start with the second point. It follows from the stronger statement: $\Ho^1(\Omega^{\infty} \mathrm{PH}(2), \Ftwo)$ is null. Let $\lop_0 -$ denote restricting to the path-component of a basepoint. Let $\tau_{\geq 1} : \mathrm{Sp} \rightarrow \mathrm{Sp}^{\geq 1}$ denote the truncation functor, such that $\pi_*(\tau_{\geq 1}(X))$ vanishes for nonpositive degrees and $\pi_*(\tau_{\geq 1}(X)) \cong \pi_*(X)$ for $* \geq 1$, where $X$ is a spectrum. It fits in a fiber sequence: $$\tau_{\geq 1} \mathrm{X} \rightarrow \mathrm{X} \rightarrow \tau_{<1}  \mathrm{X},$$ where $\tau_{<1} \mathrm{X}$ only remembers homotopy groups of $\mathrm{X}$ in degrees $*\leq 0$.

    In particular, we have an equivalence $\lop_0 \tau_{\geq 1} \mathrm{PH}(2) \rightarrow \lop_0 \mathrm{PH}(2)$. By Hurewicz theorem, we compute $$H_1(\lop_0\mathrm{PH}(2)) \cong H_1(\lop_0\tau_{\geq 1} \mathrm{PH}(2)) \cong \pi_1(\lop_0 \tau_{\geq 1} \mathrm{PH}(2)).$$ According to \cref{lemma-5.15}, the map $\mathrm{MTSO}(2) \rightarrow \mathrm{PH}(2)$ is $2$-connected. Thus, $\pi_1(\mathrm{PH}(2)) \cong \pi_1(\mathrm{MTSO}(2))$. We can conclude since the latter vanishes, according to \cite[Theorem $1.0.1$]{ebert2007lowdimensionalhomotopystablemapping}. 
    
    We now show the class $\kappa_{\epsilon}$ is not null. Similarly, we have an equivalence $\lop_0 \tau_{\geq 1} \Th(\nuparStwo) \rightarrow \lop_0 \Th(\nuparStwo)$. By naturality of the Hurewicz morphism, as well as naturality of the suspension morphism $\sigma_* : \Ho_*(\lop_0 -) \rightarrow \Ho_*(-),$ we get the following commutative diagram: 

    \[\begin{tikzcd}
	{\pi_1(\tau_{\geq 1}\Th(\nu^{\haut}_{S^2}))} & \textcolor{mypurple}{{\Ho_1(\tau_{\geq 1} \Th(\nu^{\haut}_{S^2}))}} & \textcolor{mypurple}{{\Ho_1(\Th(\nu^{\haut}_{S^2}))}} \\
	{\pi_1(\Omega^{\infty}_0\tau_{\geq 1}\Th(\nu^{\haut}_{S^2}))} & {\Ho_1(\Omega^{\infty}_0\tau_{\geq 1} \Th(\nu^{\haut}_{S^2}))} & {\Ho_1(\Omega^{\infty}_0\Th(\nu^{\haut}_{S^2}))}
	\arrow["{h_1}", from=1-1, to=1-2]
	\arrow[equals, from=1-1, to=2-1]
	\arrow["g", color=mypurple, from=1-2, to=1-3]
	\arrow["{h_1}", from=2-1, to=2-2]
	\arrow["{\sigma_*}", from=2-2, to=1-2]
	\arrow["\cong", from=2-2, to=2-3]
	\arrow["{\sigma_*}", from=2-3, to=1-3]
\end{tikzcd}.\]

The top and bottom maps $h_1$ are the Hurewicz morphisms. Since $\tau_{\geq 1}$ has its homotopy groups concentrated in degrees $*\geq 1$, both top and bottom $h_1$ are isomorphisms. By a two-out-of-three argument, we deduce $$\sigma_* : \Ho_1(\lop_0 \tau_{\geq 1} \Th(\nuparStwo)) \rightarrow \Ho_1(\tau_{\geq 1} \Th(\nuparStwo))$$ is an isomorphism. In order to show the right vertical map $$\sigma_* : \Ho_1(\lop_0 \Th(\nuparStwo)) \rightarrow \Ho_1(\Th(\nuparStwo))$$ is an isomorphism, it remains to show the morphism $$\textcolor{mypurple}{g : \Ho_1(\tau_{\geq 1} \Th(\nuparStwo)) \rightarrow \Ho_1(\Th(\nuparStwo))},$$ in purple in the diagram,
is an isomorphism.

After taking truncations, the map $j : \mathrm{MTSO}(2) \rightarrow \Th(\nuparStwo)$ induces a map of fiber sequences: 
\[\begin{tikzcd}
	{\tau_{\geq 1} \mathrm{MTSO}(2)} & {\mathrm{MTSO}(2)} & {\tau_{<1}\mathrm{MTSO}(2)} \\
	{\tau_{\geq 1} \Th(\nu^{\haut}_{S^2})} & {\Th(\nu^{\haut}_{S^2})} & {\tau_{<1}\Th(\nu^{\haut}_{S^2})} 
	\arrow[from=1-1, to=1-2]
	\arrow["{\tau_{\geq 1}j}", from=1-1, to=2-1]
	\arrow[from=1-2, to=1-3]
	\arrow["j", from=1-2, to=2-2]
	\arrow["{\tau_{<1}j}", from=1-3, to=2-3]
	\arrow[from=2-1, to=2-2]
	\arrow[from=2-2, to=2-3]
\end{tikzcd}.\]
Since $j : \mathrm{MTSO}(2) \rightarrow \Th(\nuparStwo)$ induces an isomorphism on $\pi_*$ for $*<1$, the map $$\tau_{<1} (j) : \tau_{<1} \mathrm{MTSO}(2) \rightarrow \tau_{<1}\Th(\nuparStwo)$$ is an equivalence. As a consequence, the left square is a pushout in $\mathrm{Sp}$. Let $B$ denote the cofiber of the map $j$. We write below the long exact sequence in homology groups: 

\small{\[\begin{tikzcd}
	{\Ho_1(\tau_{\geq 1}\mathrm{MTSO}(2)) } & {\Ho_1(\tau_{\geq 1}\Th(\nu^{\haut}_{S^2}))} & {\Ho_1(B)} & 0 & 0 \\
	{\Ho_1(\mathrm{MTSO}(2)) \cong 0} & {\Ho_1(\Th(\nu^{\haut}_{S^2}))} & {\Ho_1(B)} & {\Ho_0(\mathrm{MTSO}(2))} & {\Ho_0(\Th(\nu^{\haut}_{S^2}))}
	\arrow[from=1-1, to=1-2]
	\arrow[from=1-1, to=2-1]
	\arrow[from=1-2, to=1-3]
	\arrow["g",from=1-2, to=2-2]
	\arrow["0"', from=1-3, to=1-4]
	\arrow[equals, from=1-3, to=2-3]
	\arrow[from=1-4, to=1-5]
	\arrow[from=1-4, to=2-4]
	\arrow[from=1-5, to=2-5]
	\arrow[from=2-1, to=2-2]
	\arrow[from=2-2, to=2-3]
	\arrow["0"', from=2-3, to=2-4]
	\arrow["\cong"', from=2-4, to=2-5]
\end{tikzcd}\]}
\normalsize The top right-handside of the diagram is $0$ because the spectra $\tau_{\geq 1} \mathrm{MTSO}(2)$ and $\tau_{\geq 1} \Th(\nuparStwo)$ are $0$-connected. On the other hand, according to what we wrote above, $\Ho_1(\tau_{\geq 1} \mathrm{MTSO}(2))$ vanishes. Consequently, the middle maps in the diagram $\Ho_1(\tau_{\geq 1} \Th(\nuparStwo)) \rightarrow \Ho_1(B)$ and $\Ho_1( \Th(\nuparStwo)) \rightarrow \Ho_1(B)$ are isomorphisms and it follows that the map $$\textcolor{mypurple}{g: \Ho_1(\tau_{\geq 1} \Th(\nuparStwo)) \rightarrow \Ho_1(\Th(\nuparStwo))}$$ is an isomorphism. To conclude, the morphism $$\sigma_* : \Ho_1(\lop_0\Th(\nuparStwo),\Ftwo) \rightarrow \Ho_1(\Th(\nuparStwo),\Ftwo)$$ is an isomorphism. Since we work with field coefficients, the dual morphism $$\sigma^1 : \Ho^1(\Th(\nuparStwo), \Ftwo) \rightarrow \Ho^1(\lop_0 \Th(\nuparStwo),\Ftwo)$$ is also an isomorphism. Consequently, the element $\epsilon.U$ generating $\Ho^1(\Th(\nuparStwo),\Ftwo)$ is sent to a nonzero class $\kappa_{\epsilon}=\sigma^*(\epsilon.U) \in \Ho^1(\lop_0 \Th(\nuparStwo),\Ftwo),$ which concludes the proof.
\end{proof}

Before showing the functor $\B \sCobtwo(-)$ is not $1$-excisive, we give a remark below on the Madsen-Weiss theorem: 

\begin{rk}
    One initial motivation, as in \cite{GMTW}, for studying the homotopy type of the nerve of the cobordism category $\mathrm{Cob}_2^{\mathrm{SO}}$ is the cohomology of stable moduli space of surfaces \[\mathcal{M}_{\infty}=\hocolim_{g\rightarrow \infty} \B\Diff_{\partial}(\Sigma_{g,1}).\] The connection comes from the Madsen-Weiss Theorem, proven in \cite{GMTW}. The latter states that there is a map $$\mathcal{M}_{\infty} \rightarrow \lop_0 \mathrm{MTSO}(2),$$ such that it is a homology equivalence, or in other words, induces an isomorphism on homology. Here $\lop_0 \mathrm{MTSO}(2)$ denotes the restriction to the path-component of a basepoint.
    
    One can then wonder what happens if we replace diffeomorphisms with self-homotopy equivalences. In dimension $2$, according to Subsection \ref{3.1bis}, we deduce that $\mathcal{M}_{\infty}$ is equivalent to \[\hocolim_{g\rightarrow \infty} \Bhaut_{\partial}(\Sigma_{g,1}).\] However, according to \cref{thmC}, the group-completion $\Omega_{\emptyset} \B \sCobtwo$ is not homotopy or homology equivalent to $\Omega_{\emptyset} \B\mathrm{Cob}_2^{\mathrm{SO}}$. It suggests we may not have an analogue of Madsen-Weiss theorem for classifying spaces of self-homotopy equivalences of Poincaré complexes.
\end{rk}

Finally, we prove the functor $\B \sCobtwo(-)$ is not $1$-excisive. 

\begin{proof}[Proof of \cref{COR-thmC}]
According to \cref{thmC} the map $\Omega_{\emptyset} \B \sCobtwo \rightarrow \lop \Th(\nuparStwo)$ is $0$-connected and is not an equivalence. We deduce the map $$\B \sCobtwo \rightarrow \Omega^{\infty-1} \Th(\nuparStwo) \simeq \partial_1 \B \sCobtwo,$$ is not an equivalence. According to \cref{crit-ex}, we infer $\B \sCobtwo(-)$ is not $1$-excisive.
\end{proof}

\printbibliography

\end{document}